\documentclass[11pt]{article}
\usepackage{amsfonts}
\usepackage{eucal}
\usepackage{mathtools} 
\usepackage{amssymb,epsfig,latexsym}
\usepackage{enumitem}
\usepackage{graphicx,amsmath,amssymb,latexsym,psfrag}
\usepackage{amsthm}
\usepackage{graphics,graphicx}
\usepackage{pgf,tikz}
\usepackage{mathrsfs}
\usepackage{epstopdf}
\usetikzlibrary{arrows}
\DeclareMathOperator*{\Var}{\mathbb{V}ar}

\newcommand{\R}{\mathbb{R}}

\newcommand{\st}{{\ :\ }}
\renewcommand{\d}{{\mathrm d}}
\newcommand{\PP}{\mathbb{P}}

\numberwithin{equation}{section}

\newtheorem{Theorem}{Theorem}[section]
\newtheorem{Proposition}[Theorem]{Proposition}

\newtheorem{Corollary}[Theorem]{Corollary}
\newtheorem{Lemma}[Theorem]{Lemma}
\newtheorem{Example}[Theorem]{Example}
\theoremstyle{remark}
\newtheorem{remark}{Remark}[section]
\def\ind{{\rm 1\hspace{-0.90ex}1}}

\begin{document}
\allowdisplaybreaks
	
\title{\huge Poincar\'e and transportation cost inequalities for marked point processes}
\author{Ian Flint \thanks{Division of Mathematical Sciences,
Nanyang Technological University, SPMS-MAS-04-02, 21 Nanyang Link
Singapore 637371. e-mail: \tt{iflint@ntu.edu.sg}}\and Nicolas Privault\thanks{Division of Mathematical Sciences,
Nanyang Technological University, SPMS-MAS-05-43, 21 Nanyang Link
Singapore 637371. e-mail: \tt{nprivault@ntu.edu.sg}} \and Giovanni
Luca Torrisi\thanks{Istituto per le Applicazioni del Calcolo
"Mauro Picone", CNR, Via dei Taurini 19, 00185 Roma, Italy.
e-mail: \tt{giovanniluca.torrisi@cnr.it}}}
\date{}
\maketitle
\begin{abstract}
	In recent years, a number of functional inequalities have been derived for Poisson random measures, with a wide range of applications.
	In this paper, we prove that such inequalities can be extended to the setting of marked temporal point processes, under mild assumptions on their Papangelou conditional intensity.
	First, we derive a Poincar\'e inequality.
	Second, we prove two transportation cost inequalities.
	The first one refers to functionals of marked point processes with a Papangelou conditional intensity and is new even in the setting of Poisson random measures.
	The second one refers to the law of marked temporal point processes with a Papangelou conditional intensity, and extends a related inequality which is known to hold on a general Poisson space.
	Finally, we provide a variational representation of the Laplace transform of functionals of marked point processes with a Papangelou conditional intensity.
	The proofs make use of an extension of the Clark-Ocone formula to marked temporal point processes. Our results are shown to apply to classes of renewal, nonlinear Hawkes and Cox point
processes.
\end{abstract}

\section{Introduction}\label{sec:int}

Point processes with a Papangelou conditional intensity (\cite{daley2}, \cite{georgii}, \cite{lieshout}, \cite{matthes}, \cite{NZa}) constitute an important class of point process models, which generalizes the Poisson process.
Roughly speaking, the intuitive meaning of this notion of conditional intensity, denoted by $\pi_x(\omega)$, is that, for a suitable state space $S$ and reference measure $\sigma$ on $S$, $\pi_x(\omega)\,\sigma(\d x)$, $x\in S$,
is the conditional probability of having a particle in the infinitesimal region $\d x$ when the configuration $\omega$ is given outside $\d x$.

In this paper we provide several functional inequalities for marked temporal point processes having a Papangelou conditional intensity, with times in $\R_+$ and marks in a complete separable metric space $E$.

Our main achievements are
$(i)$ a Poincar\'e inequality for square-integrable functionals of marked point processes with a Papangelou conditional intensity (Theorem~\ref{prop:poincare});
$(ii)$ transportation cost inequalities for the law of functionals of marked point processes with a Papangelou conditional intensity
(Theorem \ref{prop:transport}) and for the law of the marked point process itself (Theorem \ref{thm:transportlawpp});
$(iii)$ a variational representation of the Laplace transform of functionals, bounded from above, of marked point processes with a Papangelou conditional intensity (Theorem~\ref{thm:laplace}).

The Poincar\'e inequality and variational representations of the Laplace transform for functionals on the Poisson space have attracted a lot of interest (see \cite{lastpenrosebook}, \cite{privault}, \cite{wu} for the Poincar\'e inequality and \cite{bdm}, \cite{zhang1} for variational representations of the Laplace transform).
Our results in this direction show that such functional relations hold in a non-Poissonian setting.
We emphasize that the Poincar\'e inequality proved in this paper concerns one-dimensional marked point processes and it holds under different conditions than the Poincar\'e inequality for Gibbs point processes provided in \cite{kondratiev} (we refer the reader to Remark \ref{rem:KL} for a more detailed discussion).
On the one hand, to the best of our knowledge, the transportation cost inequality for functionals of marked point processes with a Papangelou conditional intensity is new even in the Poisson setting. On the other hand, the transportation cost inequality for the law of the marked point process itself generalizes a related relation proved in \cite{MaWu}.

A key ingredient in the proofs is a new Clark-Ocone formula for square-integrable functionals of marked point processes with a Papangelou conditional intensity (Theorem~\ref{thm:clarkoperative}), which generalizes the corresponding formula in \cite{COFormula} in two directions.
First, we allow for point processes with values on an unbounded time interval and marks in a complete separable metric space.
Second, we prove the square integrability of the
Clark-Ocone integrand, which is crucial in the proofs of the functional inequalities mentioned above since it enables the application of the isometry formula provided by Proposition~\ref{prop:ExtendedIsometry}.

Deviation bounds,
which are a classical application of transportation cost inequalities,
are presented in Remarks \ref{re:deviationbound} and \ref{re:deviationboundomega}.
The variational representation of the Laplace transform for functionals of marked point processes with a Papangelou conditional intensity
can also be useful to derive large deviation principles for those functionals, along the lines of \cite{bdm} and more generally relying on the theory developed in \cite{dupuisellis}.
Further applications of our results to
various classes of non-Poissonian point processes such as renewal, nonlinear Hawkes and Cox are presented in Corollaries \ref{cor:PoincareRenewal}, \ref{cor:Hawkes}, \ref{cor:Cox},
\ref{cor:transportRenewal}, \ref{cor:transportHawkes}, \ref{cor:transportCox}, \ref{cor:transportRenewalForLawN}, \ref{cor:transportHawkesForLawN},
\ref{cor:transportCoxForLawN} and in Remark~\ref{re:variationalrep}.

The proof of the Poincar\'e inequality is based on the evaluation of the variance of the functional by a combination of the Clark-Ocone formula and the isometry formula for marked point processes with a stochastic intensity.
The proof of the transportation cost inequality exploits its characterization via exponential moments proved in \cite{gozlan}.
Such exponential moments are controlled by a stochastic convex inequality for functionals of marked point processes with a Papangelou conditional intensity (Proposition~\ref{lem:Theorem41Privault}), which is based on the Clark-Ocone formula and generalizes the corresponding result in \cite{kleinmaprivault}.
The proof of the variational representation of the Laplace transform
uses a localization argument to deal with the unbounded case,
  which is out of the reach of the techniques in \cite{zhang1},
  and relies on the Clark-Ocone formula
in order to take into account the non-Poissonian dynamics of the point process.
The proof of the Clark-Ocone formula is based on the representation theorem for square-integrable martingales,
on the use of an integration by parts formula for functionals of point processes with a Papangelou conditional intensity,
and on an isometry formula for point processes with stochastic intensity,
which allows us to identify the integrand appearing in the representation theorem.

We note once again that, in contrast to the corresponding result in \cite{COFormula}, the Clark-Ocone formula of Theorem~\ref{thm:clarkoperative} guarantees the square integrability of the integrand \eqref{eq:ExplicitIntegrand}.
This integrability property is crucial to the proofs of our main results due to a pervasive use of the isometry formula provided by Proposition~\ref{prop:ExtendedIsometry}.

The paper is organized as follows.
In Section~\ref{sec:Preliminaries} we give some preliminaries on point processes including the notions of Papangelou conditional intensity, classical stochastic intensity and an important relation between them.
In Section~\ref{sec:mainResults} we describe the main results of the paper and give their proofs in Section~\ref{sec:main}.
We also include an appendix, where we prove some technical lemmas and propositions.

\section{Preliminaries on point processes}
\label{sec:Preliminaries}

Let $E$ be a complete separable metric space and $\mathcal E$ the corresponding Borel $\sigma$-field.
Let $\Omega$ be the set of all integer-valued measures $\omega$ on $(\R_+\times E,\mathcal{B}(\R_+)\otimes \mathcal E)$, where $\mathcal B(\R_+)$ is the Borel $\sigma$-field on $\R_+$, such that $\omega(\{0\}\times E)=0$,
$\omega(\{t\}\times E)\leq 1$ for any $t\in\R_+$, and $\omega(K\times E)<\infty$ for any compact set $K\subset\R_+$
(in Remark~\ref{rem:JustificationOfOurChoiceOfOmega}, we shall explain in what sense our results apply under a more general definition of $\Omega$). We define
\[
	N(\omega):=\omega,\quad\omega\in\Omega
\]
and for a Borel set $B\in\mathcal{B}(\R_+)$ we shall consider the $\sigma$-field
\[
	\mathcal{F}_B:=\sigma\{N(A\times D):\,\,A\in\mathcal{B}(B),D\in\mathcal{E}\},
\]
where $\mathcal{B}(B)$ denotes the restriction of $\mathcal{B}(\R_+)$ to $B$.
For ease of notation we set $\mathcal{F}_t:=\mathcal{F}_{[0,t]}$ and $\mathcal{F}_{t^-}:=\mathcal{F}_{[0,t)}\equiv\bigvee_{0\leq s<t}\mathcal{F}_s$ for $t\in\R_+$.

We set $\mathcal{F}_\infty:=\bigvee_{t\in\R_+}\mathcal{F}_t$, let $\PP$ be a probability measure on $(\Omega,\mathcal F_\infty)$ and consider the canonical probability space $(\Omega,\mathcal F_\infty,\mathbb P)$.
The elements of $(\Omega,\mathcal F_\infty,\mathbb P)$ are known in the literature as simple and locally finite marked point processes on $\R_+$ with marks in $E$.
Throughout the paper we denote by $\mathbb E$ and $\Var$ the expectation and the variance operators with respect to $\mathbb P$, respectively.

By analogy with the un-marked setting, a mapping $X:\R_+\times\Omega\times E\to\R$ which is measurable with respect to the $\sigma$-fields $(\mathcal{B}(\R_+)\otimes\mathcal F_\infty\otimes\mathcal E,\mathcal{B}(\R))$,
where $\mathcal{B}(\R)$ is the Borel $\sigma$-field on $\R$, is called a stochastic process. Let $\mathcal{G}:=\{\mathcal{G}_t\}_{t\in\R_+}$ be a filtration such that $\mathcal{F}_\infty\supset\mathcal{G}_t\supseteq\mathcal{F}_t$,
$t\in\R_+$.
The $\mathcal{G}$-predictable $\sigma$-field on $\R_+\times\Omega$, denoted by $\mathcal{P}(\mathcal G)$, is the $\sigma$-field generated by the sets $(a,b]\times A$ with $A\in\mathcal{G}_a$, $a,b\in\R_+$.
Throughout the paper, for any $\omega\in\Omega$ and $(t,x)\in\mathbb{R}_+\times E$, we define $\omega-\varepsilon_{(t,x)}:=\omega$ if $(t,x)\notin\mathrm{Supp}(\omega)$,
where $\varepsilon_{(t,x)}$ is the Dirac measure at $(t,x)$ and $\mathrm{Supp}(\omega)$ denotes the support of $\omega$.
A stochastic process $X:\R_+\times\Omega\times E\to\R$ is said to be $\mathcal{G}$-predictable if it is measurable with respect to the $\sigma$-fields $(\mathcal{P}(\mathcal G)\otimes\mathcal{E},\mathcal{B}(\R))$.
We say that $X$ is predictable if it is measurable with respect to the $\sigma$-fields $(\mathcal{P}(\mathcal F)\otimes\mathcal{E},\mathcal{B}(\R))$,
where $\mathcal{F}:=\{\mathcal{F}_t\}_{t\in\R_+}$. For ease of notation, we set $X_{(t,x)}(\omega):=X(t,\omega,x)$ and for later purposes, we mention that if $X$ is predictable, then for fixed $t\in\R_+$ and $x\in E$
the random variable $X_{(t,x)}$ is measurable with respect to $\mathcal{F}_{t^-}$ (and therefore with respect to $\mathcal{F}_{t}$), and so $X_{(t,x)}(\omega)=X_{(t,x)}(\omega-\varepsilon_{(t,x)})$.
This claim follows by an obvious modification of the proof of Lemma~A3.3.I~p.~425 in \cite{daley}, see also Proposition~3.3 in \cite{last}.

We shall consider two different notions of conditional intensity for marked point processes: the Papangelou conditional intensity and the classical stochastic intensity.

\subsection{Marked point processes with a Papangelou conditional intensity}

Let $\nu$ denote a locally finite on $(E,\mathcal E)$. A non-negative stochastic process $\pi:\R_+\times\Omega\times E\to\R_+$ is said to be a Papangelou conditional intensity of $N$ with respect to $\d t\nu(\mathrm{d}x)$
if, for any non-negative stochastic process $X:\R_+\times\Omega\times E\to\R_+$,
\begin{equation}\label{eq:GNZ}
	\mathbb{E}\left[\int_{\R_+\times E}X_{(t,x)}(N-\varepsilon_{(t,x)})\,N(\d t\times\d x)\right]=\mathbb{E}\left[\int_{\R_+\times E}X_{(t,x)}(N)\,\pi_{(t,x)}(N)\,\d t\nu(\d x)\right].
\end{equation}
Intuitively, $\pi_{(t,x)}(\omega)\,\mathrm dt \nu(\d x)$ is the probability that the point process has a point in the infinitesimal region $\mathrm dt\d x$
given that it agrees with the configuration $\omega$ outside of $\mathrm dt\d x$. Point processes with a Papangelou conditional intensity are fully characterized in \cite{matthes} (see Section~3) and \cite{NZa} (see Theorem~2').

\begin{remark}
\label{rem:JustificationOfOurChoiceOfOmega}
	Some texts (e.g. \cite[Definitions~6.4.I.]{daley}) use a more general definition of a marked point process which allows for $N$ to have atoms $\varepsilon_{(t,x)},\varepsilon_{(t,y)}$ for the same $t\in\R_+$ and $x\neq y$.
	We considered a more restrictive set $\Omega$ in Section~\ref{sec:Preliminaries} since all our main results hold for marked point processes with a Papangelou conditional intensity.
	In such case, there exists a version of the more general marked point process which takes its values in the set $\Omega$, see the following Lemma~\ref{lem:DifferentDefinitionsOfMPP}.
\end{remark}

\begin{Lemma}
	\label{lem:DifferentDefinitionsOfMPP}
	Let $\Omega'$ be the set of all integer-valued measures $\omega$ on $(\R_+\times E,\mathcal{B}(\R_+)\otimes \mathcal E)$, such that $\omega(\{0\}\times E)=0$, $\omega(\{t\}\times\{x\})\leq 1$ for any $(t,x)\in\R_+\times E$, and $\omega(K\times E)<\infty$ for any compact set $K\subset\R_+$ (note that $\Omega\subset\Omega'$).
	Assume that the marked point process $N$ is defined on $\Omega'$ instead of $\Omega$ and define the quantities from Section~\ref{sec:Preliminaries} in an analogous manner.
	In particular,
        if $\mathbb P$ is
        a probability measure on $(\Omega',\mathcal F_\infty)$, where $\mathcal F_\infty$ is appropriately defined
        and $N$ has a Papangelou conditional intensity, then $N$ takes its values in $\Omega$ $\mathbb P$-almost surely, i.e. $\mathbb P(\Omega)=1$.
\end{Lemma}
We postpone the proof of this lemma to Section~\ref{subsec:ProofOfDifferentDefinitionsOfMPP} in the appendix in order to improve the flow of the
article.

Throughout this paper we shall often consider locally stable point processes, i.e. point processes $N$ with a Papangelou conditional intensity $\pi$ such that
\begin{equation}\label{eq:LocalStability}
	\pi_{(s,x)}(\omega)\le\beta(s,x),\quad\d s\nu(\d x)\mathbb{P}(\d\omega)\text{-almost everywhere}
\end{equation}
for a function $\beta:\R_+\times E\to\R$ which is integrable with respect to $\d s\d\nu$, on $[0,t]\times E$ for every $t\in\mathbb{R}_+$.
The local stability is known to be satisfied by a wide range of point processes (see e.g. \cite{lieshout}).

In \cite{last2}, the author provides a condition guaranteeing the existence of the predictable projection of a bounded
stochastic process. In the following proposition, we specialize this result to our setting while relaxing the
boundedness assumption.
\begin{Proposition}\label{cor:PredictableProjection}
	Assume that $N$ has a Papangelou conditional intensity $\pi$, i.e. \eqref{eq:GNZ} holds.
	Let $X:\R_+\times\Omega\times E\to\R$ be a stochastic process which is such that either $(i)$ $X\ge0$ or $(ii)$ for $\mathrm{d}t\nu(\mathrm{d}x)$-almost all $(t,x)\in\R_+\times E$, $X_{(t,x)}\in L^1(\Omega,\mathcal{F}_\infty,\mathbb P)$.
	Then there exists a predictable stochastic process $p(X):\R_+\times\Omega\times E\to\R$ called predictable projection, which is such that for $\mathrm{d}t\nu(\mathrm{d}x)$-almost all $(t,x)\in\R_+\times E$, we have
	\begin{equation*}
		p(X)_{(t,x)}=\mathbb E\bigl[X_{(t,x)}\big| \mathcal F_{t^-}\bigr],\quad\mathbb P\text{-almost surely.}
	\end{equation*}
	Additionally, under $(ii)$ the predictable projection $p(X)$ is such that, for $\mathrm{d}t\nu(\mathrm{d}x)$-almost all $(t,x)\in\R_+\times E$,
	\begin{equation*}
		p(X)_{(t,x)}<\infty,\quad\mathbb P\text{-almost surely.}
	\end{equation*}
\end{Proposition}
The proof is rather technical and therefore postponed to Section~\ref{subsec:ProofOfPredictableProjection} in the appendix.

Throughout this paper we will use the discrete Malliavin derivative of $F:\Omega\to\R$, defined as
\begin{equation*}
	D_{(t,x)} F(\omega):=F_{(t,x)}^+(\omega)-F(\omega),\quad (t,x)\in\R_+\times E,\ \omega\in\Omega,
\end{equation*}
where
\begin{equation*}
	F_{(t,x)}^+(\omega)
	:=F(\omega+\varepsilon_{(t,x)})
	=\begin{cases}
		F(\omega)&\text{ if }(t,x)\in\mathrm{Supp}(\omega),\\
		F(\omega+\varepsilon_{(t,x)})&\text{ if }(t,x)\notin\mathrm{Supp}(\omega).
	\end{cases}
\end{equation*}

Under suitable integrability conditions on $X$ and $\pi$, we shall consider the stochastic integral
\[
	\Delta(X):=\int_{\R_+\times E} X_{(t,x)}\,(N(\d t\times\d x)-\pi_{(t,x)}\,\d t\nu(\d x)),\quad\text{$\mathbb P$-almost surely}.
\]

If $X$ is predictable, then $X_{(t,x)}=X_{(t,x)}(N)=X_{(t,x)}(N-\varepsilon_{(t,x)})$, and so this integral can be rewritten as
\begin{equation*}
	\Delta(X)=\int_{\R_+\times E} X_{(t,x)}(N-\varepsilon_{(t,x)})\,(N(\d t\times\d x)-\pi_{(t,x)}\,\d t\nu(\d x)).
\end{equation*}
We conclude this paragraph
with the following integration by parts formula, see Corollary~3.1 in \cite{torrisi}.

\begin{Lemma}\label{le:IBP}
	Assume that $N$ has a Papangelou conditional intensity $\pi$. Then, for any predictable stochastic process $X:\mathbb R_+\times\Omega\times E\to\mathbb R$ and any random variable $G:\Omega\to\R$ such that
	\begin{equation}\label{eq:condint}
		\mathbb E\biggl[\int_{\R_+\times E}\bigl|GX_{(t,x)}\bigr|\pi_{(t,x)}\,\d t\nu(\d x)\biggr]<\infty\quad\text{and}\quad\mathbb E\biggl[\int_{\R_+\times E}\bigl|X_{(t,x)}D_{(t,x)}G\bigr|\pi_{(t,x)}\,\d t\nu(\d x)\biggr]<\infty,
	\end{equation}
	we have
	\begin{equation*}
		\mathbb E\bigl[G\Delta(X)\bigr]=\mathbb E\biggl[\int_{\R_+\times E}X_{(t,x)}\pi_{(t,x)}D_{(t,x)}G\,\d t\nu(\d x)\biggr].
	\end{equation*}
\end{Lemma}

\subsection{Marked point processes with a classical stochastic intensity}

A non-negative and $\mathcal G$-predictable stochastic process $\lambda:\R_+\times\Omega\times E\to\R_+$ such that for all
$t\in\mathbb{R}_+$ $\int_{[0,t]\times E}\lambda_{(s,x)}\,\d s\nu(\d x)<\infty$ $\mathbb P$-almost surely is said
to be a $\mathcal G$-stochastic intensity of $N$ (see \cite{bremaud}, \cite{daley} and \cite{lastbrandt}) if for any non-negative and $\mathcal G$-predictable stochastic process $X:\R_+\times\Omega\times E\to\R_+$ we have
\begin{equation}\label{eq:classicalStochasticIntensity}
	\mathbb{E}\biggl[\int_{\R_+\times E}X_{(t,x)}\,N(\d t\times\d x)\biggr]=\mathbb{E}\biggl[\int_{\R_+\times E}X_{(t,x)}\lambda_{(s,x)}\,\d s\nu(\d x)\biggr].
\end{equation}
Additionally, we call $\lambda$ a (classical) stochastic intensity if $\mathcal{G}:=\mathcal F$.
Roughly speaking, the quantity $\lambda_{(t,x)}(\omega)\,\mathrm dt\nu(\d x)$ is the probability that the point process has a point in the infinitesimal region
$\mathrm dt\d x$ given that it agrees with the configuration $\omega$ on $[0,t)\times E$.

Hereafter, we assume that $N$ has a $\mathcal G$-stochastic intensity $\lambda$ and, given two $\mathcal G$-predictable stochastic processes $X,Y:\R_+\times\Omega\times E\to\R$, we write $X\sim Y$ if $X$ and $Y$ are equal $\lambda_{(t,x)}(\omega)\,\d t\mathbb P(\d\omega)\nu(\d x)$-almost everywhere on $\R_+\times\Omega\times E$.

For $p\in[1,+\infty)$, we denote by $\mathcal{P}_p(\lambda)$ the family of equivalence classes (with respect to the equivalence relation $\sim$) formed by $\mathcal{G}$-predictable stochastic processes $X$ such that
\begin{equation*}
	\|X\|_{\mathcal P_p(\lambda)}^p:=\mathbb{E}\biggl[\int_{\R_+\times E}|X_{(t,x)}|^p\,\lambda_{(t,x)}\,\d t\nu(\d x)\biggr]<\infty
\end{equation*}
and, for ease of notation, we set $\mathcal{P}_{1,2}(\lambda):=\mathcal{P}_{1}(\lambda)\cap\mathcal{P}_{2}(\lambda)$.
Note that $\|\cdot\|_{\mathcal P_p(\lambda)}$ is a norm on $\mathcal{P}_p(\lambda)$.

For any $X\in\mathcal P_1(\lambda)$ we define the stochastic integral
\begin{equation*}
	\delta(X):=\int_{\R_+\times E}X_{(t,x)}\,(N(\d t\times\d x)-\lambda_{(t,x)}\,\d t\nu(\d x)),\quad\text{$\mathbb P$-almost surely,}
\end{equation*}
which is well defined $\mathbb P$-almost everywhere as the difference of two finite terms by \eqref{eq:classicalStochasticIntensity}.
The next proposition provides a fundamental isometry formula for marked point processes with a $\mathcal G$-stochastic intensity.

\begin{Proposition}[Theorem~3 of \cite{bremaudmassoulie} and p.~62 of \cite{watanabe}]
\label{bremaudmassoulie}
	Assume that $N$ has a $\mathcal G$-stochastic intensity $\lambda$. Then:\\
	$(i)$ $\mathbb E\bigl[\delta(X)\bigr]=0$, for any $X\in\mathcal{P}_1(\lambda)$;\\
	$(ii)$ $\mathbb E\bigl[\delta(X)\delta(Y)\bigr]=\mathbb E\Bigl[\int_{\R_+\times E}X_{(t,x)}Y_{(t,x)}\lambda_{(t,x)}\,\d t\nu(\d x)\Bigr]$, for any $X,Y\in\mathcal{P}_{1,2}(\lambda)$;\\
	$(iii)$ If $\{\ind_{[0,t]}(s)X_{(s,x)}\}_{(s,x)\in\R_+\times E},\{\ind_{[0,t]}(s)Y_{(s,x)}\}_{(s,x)\in\R_+\times E}\in\mathcal{P}_{1,2}(\lambda)$ for any $t\in\R_+$, then the stochastic process
	\[
		\delta(\ind_{[0,t]}X)\delta(\ind_{[0,t]}Y)-\int_{[0,t]\times E}X_{(s,x)}Y_{(s,x)}\lambda_{(s,x)}\,\d s\nu(\d x),\quad t\in\R_+
	\]
	is a $\mathcal{G}$-martingale.
\end{Proposition}

For later purposes, we extend the operator $\delta$ to $\mathcal P_1(\lambda)\cup\mathcal P_2(\lambda)$ and prove that the isometry formula of Proposition \ref{bremaudmassoulie}-$(ii)$ holds on $\mathcal P_2(\lambda)$ for the extension of $\delta$.
\begin{Proposition}\label{prop:ExtendedIsometry}
	Assume that $N$ has a $\mathcal G$-stochastic intensity $\lambda$ and that $\mathbb E\bigl[N([0,t]\times E)\bigr]<\infty$, for all $t\in\mathbb{R}_+$.
	The operator $\delta$ with domain $\mathcal P_1(\lambda)$ can be uniquely extended to an operator with domain $\mathcal P_1(\lambda)\cup\mathcal P_2(\lambda)$, which we still denote by $\delta$.
	Additionally, for any $X,Y\in\mathcal P_2(\lambda)$ we have
	\begin{equation*}
		\mathbb E\bigl[\delta(X)\delta(Y)\bigr]
		=\mathbb E\biggl[\int_{\R_+\times E}X_{(t,x)}Y_{(t,x)}\lambda_{(t,x)}\,\d t\nu(\d x)\biggr].
	\end{equation*}
\end{Proposition}
Although the proof of this proposition follows a standard Cauchy sequence argument along the lines of Section~II.2 in \cite{watanabe}, we have not found this precise result in the literature.
For this reason, we provide the proof in Section~\ref{subsec:ProofOfExtendedIsometry} of the appendix.

\subsection{A relation between the Papangelou conditional intensity and the classical stochastic intensity}

In the next lemma we show that the Papangelou conditional intensity of a marked point process $N$ determines its stochastic intensity.
\begin{Lemma}\label{le:PapStochlink}
	If $N$ has a Papangelou conditional intensity $\pi$ and
	\begin{equation}\label{eq:stocint}
		\int_{[0,t]\times E}p(\pi)_{(s,x)}\,\d s\nu(\d x)<\infty,\quad\text{$\mathbb P$-almost surely, for all $t\in\mathbb R_+$}
	\end{equation}
	then $N$ has a stochastic intensity $p(\pi)$.
\end{Lemma}
\noindent{\it Proof.}
Let $X$ be a non-negative and predictable stochastic process.
As recalled at the beginning of Section~\ref{sec:Preliminaries}, we have $X_{(t,x)}(\omega)=X_{(t,x)}(\omega-\varepsilon_{(t,x)})$ and $X_{(t,x)}$ is $\mathcal{F}_{t^-}$-measurable, for any $t\in\mathbb R_+$ and $x\in E$.
So by Fubini's theorem, standard properties of the conditional expectation and Proposition~\ref{cor:PredictableProjection}-$(i)$,
\begin{align*}
	\mathbb{E}\left[\int_{\mathbb R_+\times E} X_{(t,x)}\,N(\d t\times\d x)\right]&=
	\int_{\mathbb R_+\times E}\mathbb{E}\bigl[X_{(t,x)}\pi_{(t,x)}(N)\bigr]\,\d t\nu(\d x)\nonumber\\
	&=\int_{\mathbb R_+\times E}\mathbb{E}\bigl[X_{(t,x)}\mathbb{E}\bigl[\pi_{(t,x)}(N)\big| \mathcal{F}_{t^-}\bigr]\bigr]\,\d t\nu(\d x)\nonumber\\
	&=\mathbb{E}\biggl[\int_{\mathbb R_+\times E} X_{(t,x)}p(\pi)_{(t,x)}\,\d t\nu(\d x)\biggr].
\end{align*}
\\
\noindent$\square$

As already mentioned, the proofs of our main results are based on a Clark-Ocone formula for marked point processes, which, in turn, exploits the martingale representation theorem.
For this reason, we shall consider the $\mathbb P$-completed filtration $\overline{\mathcal F}\equiv\{\overline{\mathcal{F}}_t\}_{t\in\R_+}$, defined by
$\overline{\mathcal{F}}_t:=\mathcal{F}_t\vee\mathcal{N}=\sigma(\{A\cup H:\,\,A\in\mathcal{F}_t,\,H\in\mathcal N\})$, where $\mathcal{N}$ is the family of $\mathbb{P}$-null events of $\mathcal{F}_\infty$ (see \cite{bremaud} p.~309).

The next lemma, which we prove for the sake of completeness, guarantees that the notion of stochastic intensity is equivalent to that of $\overline{\mathcal F}$-stochastic intensity.
\begin{Lemma}\label{le:equiv}
	Let $\lambda:\R_+\times\Omega\times E\to\R_+$ be a predictable stochastic process.
	Then $N$ has a stochastic intensity $\lambda$ if and only if $N$ has a $\overline{\mathcal{F}}$-stochastic intensity $\lambda$.
\end{Lemma}
\noindent{\it Proof.}
We start by noting that any predictable stochastic process is $\overline{\mathcal{F}}$-predictable, and in particular $\lambda$ is both predictable and $\overline{\mathcal F}$-predictable.
The necessity follows directly from this fact.
For the sufficiency, taking in \eqref{eq:classicalStochasticIntensity} the process defined by $X_{(t,x)}:=\ind_{(a,b]}(t)\ind_A(\omega)\ind_L(x)$, $a,b\in\R_+$, $A\in\mathcal{F}_a$, $L\in\mathcal E$ we have
\[
	\mathbb E\bigl[N((a,b]\times L)\big| \mathcal{F}_a\bigr]
	=\mathbb E\biggl[\int_{(a,b]\times L}\lambda_{(t,x)}\d t\nu(\d x)\;\Big|\;\mathcal{F}_a\biggr],\quad\text{$\mathbb P$-almost surely.}
\]
We shall check later on that, for any 
non-negative or integrable random variable $Y$ and $t\in\mathbb R_+$
we have
\begin{equation}\label{eq:barnotbar}
	\mathbb E\bigl[Y\big| \mathcal{F}_t\bigr]
	=\mathbb E\bigl[Y\big| \overline{\mathcal{F}}_t\bigr],\quad\text{$\mathbb P$-almost surely.}
\end{equation}
As a consequence of \eqref{eq:barnotbar}, for any $a,b\in\R_+$ and $L\in\mathcal E$, we have
\begin{equation*}
	\mathbb E\bigl[N((a,b]\times L)\big| \overline{\mathcal{F}}_a\bigr]
	=\mathbb E\biggl[\int_{(a,b]\times L}\lambda_{(t,x)}\d t\nu(\d x)\;\Big|\;\overline{\mathcal{F}}_a\biggr],\quad\text{$\mathbb P$-almost surely.}
\end{equation*}
By this relation and a standard application of the monotone class theorem (see e.g. \cite{bremaud}, Theorem~T1 p.~260), we deduce
\begin{align*}
	\mathbb{E}\left[\int_{\mathbb R_+\times E} V_{(t,x)}\,N(\d t\times\d x)\right]
	=\mathbb{E}\biggl[\int_{\mathbb R_+\times E} V_{(t,x)}\lambda_{(t,x)}\,\d t\nu(\d x)\biggr],
\end{align*}
for any non-negative and $\overline{\mathcal{F}}$-predictable mapping $V$, and the proof is complete since $\lambda$ is $\overline{\mathcal F}$-predictable.

It remains to check \eqref{eq:barnotbar}.
For any $A\in\mathcal{F}_t$ and $H\in\mathcal N$,
\begin{align*}
	\mathbb E\bigl[\mathbb E\bigl[Y\big| \mathcal{F}_t\bigr]\ind_{A\cup H}\bigr]
	=\mathbb E\bigl[\mathbb E\bigl[Y\big| \mathcal{F}_t\bigr]\ind_A]+\mathbb E\bigl[\mathbb E\bigl[Y\big| \mathcal{F}_t\bigr]\ind_{H\setminus A}\bigr]
	&=\mathbb E\bigl[\mathbb E\bigl[Y\big| \mathcal{F}_t\bigr]\ind_A\bigr]\\
	&=\mathbb E\bigl[Y\ind_A\bigr]\\
	&=\mathbb E\bigl[Y\ind_{A\cup H}\bigr]-\mathbb E\bigl[Y\ind_{H\setminus A}\bigr]\\
	&=\mathbb E\bigl[Y\ind_{A\cup H}\bigr].
\end{align*}
Since $\mathbb E\bigl[Y\big| \mathcal{F}_t\bigr]$ is $\overline{\mathcal{F}}_t$-measurable, we conclude by the characterizing property of the conditional expectation.
\\
\noindent$\square$

\section{Main results}\label{sec:mainResults}

In this section we state our main achievements, i.e. a Poincar\'e inequality for square-integrable functionals of marked point processes with a Papangelou conditional intensity, a transportation cost inequality for functionals of marked point processes with a Papangelou conditional intensity, a transportation cost inequality for the law of a marked point process with a Papangelou conditional intensity,
and a variational representation formula for the Laplace transform of
(bounded from above) functionals of marked point processes with a Papangelou conditional intensity.
All these functional relations are obtained by applying a Clark-Ocone formula for square-integrable functionals of space-time point processes with a Papangelou conditional intensity, which generalizes
in various directions the corresponding formula in \cite{COFormula}.

Hereafter, we work under the convention $0/0:=0$.
Moreover, in order to be more precise, from now on we shall write $\mathcal{P}_p^{\mathcal G}(\lambda)$ in place of $\mathcal{P}_p(\lambda)$, $p\geq 1$, to stress that the stochastic processes therein are predictable with respect to some specific filtration $\mathcal G$.

\subsection{Poincar\'e inequality}

The following Poincar\'e inequality holds for functionals of marked point processes with a Papangelou conditional intensity.

\begin{Theorem}\label{prop:poincare}
	Assume
	\begin{equation}\label{eq:squareint}
		\mathbb E\bigl[N([0,t]\times E)^2\bigr]<\infty,\quad\text{for all $t\in\R_+$},
	\end{equation}
	and that $N$ has a Papangelou conditional intensity $\pi$ such that
	\begin{equation}\label{eq:IntegrabilityOfPi}
		\mathbb E\biggl[\biggl(\int_{[0,t]\times E}\pi_{(s,x)}\,\d s\nu(\d x)\biggr)^2\biggr]<\infty,\quad\text{for all $t\in\R_+$},
	\end{equation}
	and	
	\begin{equation}\label{eq:AssumptionPoincare}
		\gamma:=\left\|\int_{\R_+\times E}\frac{\Var\bigl[\pi_{(t,x)}\mid\mathcal F_{t^-}\bigr]}{\mathbb E\bigl[\pi_{(t,x)}\mid \mathcal F_{t^-}\bigr]}\,\mathrm dt\nu(\d x)\right\|_{L^\infty(\Omega,\mathcal{F}_\infty,\mathbb P)}<1,
	\end{equation}
	where $\Var\bigl[\pi_{(t,x)}\big| \mathcal F_{t^-}\bigr]:=\mathbb E\bigl[\pi_{(t,x)}^2\big| \mathcal F_{t^-}\bigr]-\mathbb E\bigl[\pi_{(t,x)}\big| \mathcal F_{t^-}\bigr]^2$.
	Then, for any $G\in L^2(\Omega,\mathcal{F}_\infty,\mathbb P)$
        we have
	\begin{equation}\label{eq:PoincareL2}
		\Var(G)\le\bigl(1-\sqrt\gamma\bigr)^{-2}\,\mathbb E\biggl[\int_{\R_+\times E}\pi_{(t,x)}\left|D_{(t,x)} G\right|^2\,\mathrm dt\nu(\d x)\biggr].
	\end{equation}
\end{Theorem}

\begin{remark}\label{re:poincare}
	Note that if $\mathbb E\bigl[\pi_{(t,x)}\big| \mathcal{F}_{t^-}\bigr]=0$ $\mathbb P$-almost surely, then $\mathbb P(\pi_{(t,x)}=0\mid\mathcal{F}_{t^-} )=1$ $\mathbb P$-almost surely.
	Therefore $\Var\bigl[\pi_{(t,x)}\big| \mathcal F_{t^-}\bigr]=0$ $\mathbb P$-almost surely and the ratio $\Var\bigl[\pi_{(t,x)}\big| \mathcal F_{t^-}\bigr]/\mathbb E\bigl[\pi_{(t,x)}\big| \mathcal F_{t^-}\bigr]$
        vanishes $\mathbb P$-almost surely by the convention $0/0:=0$.
\end{remark}

\begin{remark}\label{rem:locally}
	If $N$ is a Poisson process with locally integrable intensity function $\sigma(t,x)$ then $\pi\equiv\sigma$ and all the assumptions of Theorem~\ref{prop:poincare} hold.
	In particular, $\gamma=0$ and we recover the well-known Poincar\'e inequality for Poisson functionals (see e.g. \cite{last}, \cite{lastpenrosebook} and \cite{wu}).
\end{remark}

\begin{remark}\label{rem:AssumptionsAreSatisfiedIf}
	If $N$ is locally stable in the sense of \eqref{eq:LocalStability}, then assumptions \eqref{eq:squareint} and \eqref{eq:IntegrabilityOfPi} hold.
	Indeed, it is clear that \eqref{eq:LocalStability} implies \eqref{eq:IntegrabilityOfPi}.
	As far as the implication of \eqref{eq:squareint}, note that under assumption \eqref{eq:LocalStability} we have
	\begin{align}
		&\mathbb E\bigl[ N([0,t]\times E)^2\bigr] \nonumber\\
		&\qquad= \mathbb E\left[\int_{([0,t]\times E)^2}\pi_{(t,x)}(N)\pi_{(s,y)}(N+\varepsilon_{(t,x)})\,\d s \d t \nu(\d x)\nu(\d y)\right]
		+ \mathbb {E}\left[\int_{[0,t]\times E}\pi_{(t,x)}(N)\,\d t \nu(\d x)\right]\nonumber\\
		&\qquad\leq\int_{([0,t]\times E)^2}\beta(t,x)\beta(s,y)\,\d s \d t \nu(\d x) \nu(\d y)
		+\int_{[0,t]\times E}\beta(t,x)\,\d t\nu(\d x)<\infty,
	\end{align}
	where the first equality follows by the iterated Georgii-Nguyen-Zessin equation, see e.g. Proposition~15.5.II in \cite{daley2}.
	We also remark that \eqref{eq:IntegrabilityOfPi} is automatically verified if \eqref{eq:squareint} holds and the marked point process is attractive, in the sense that
	\begin{equation*}
		\pi_{(s,x)}(\omega)\le\pi_{(s,x)}(\eta),
	\end{equation*}
	for $\d s\nu(\mathrm{d}x)$-almost all $(s,x)\in\R_+\times E$ and all $\omega,\eta\in\Omega$ such that $\mathrm{Supp}(\omega)\subset\mathrm{Supp}(\eta)$ (see e.g. equation (3.7) in \cite{lieshout} for more explanations on this notion).
	Indeed, for all $t\in\R_+$ we have
	\begin{align*}
		\Bigl\|\int_{[0,t]\times E}\pi_{(s,x)}(N)\,\d s\nu(\d x)\Bigr\|_{L^2(\Omega,\mathcal{F}_\infty,\mathbb P)}
		&=
                \sqrt{
                  \mathbb E\biggl[\biggl(\int_{[0,t]\times E}\pi_{(s,x)}(N)\,\d s\nu(\d x)\biggr)^2\biggr]}
                \\
		&=
                \sqrt{
                  \mathbb E\biggl[\int_{[0,t]\times E}\int_{[0,t]\times E}\!\!\pi_{(s,x)}(N-\varepsilon_{(u,y)})\,\d s\nu(\d x)N(\d u\times\d y)\biggr]}
                \\
		&\le
                \sqrt{
                  \mathbb E\biggl[\int_{[0,t]\times E}N([0,t]\times E)\pi_{(s,x)}(N)\,\d s\nu(\d x)\biggr]}
                \\
        &=
                \sqrt{
                	\mathbb E\biggl[\int_{[0,t]\times E}(N-\varepsilon_{(s,x)})([0,t]\times E)\,N(\d s\times\d x)\biggr]}\\
		&=
                \sqrt{
                  \mathbb E\bigl[N([0,t]\times E)(N([0,t]\times E)-1)\bigr]
                }\\
		&<\infty,
	\end{align*}
	where we have applied \eqref{eq:GNZ} twice, and used the fact that by hypothesis $\mathbb E\bigl[N([0,t]\times E)^2\bigr]<\infty$.
\end{remark}

\begin{remark}\label{rem:locallybis}
	Let $N$ be a point process on $\R_+\times E$ with a Papangelou conditional intensity $\pi$ such that \eqref{eq:LocalStability} holds with dominating function $\beta$, for some measurable non-negative function $\alpha:\R_+\times E\to\R$,
	\[
		p(\pi)_{(t,x)}\geq\alpha(t,x)\quad\text{$\mathrm{d}t\nu(\mathrm{d}x)\mathbb{P}(\mathrm{d}\omega)$-almost everywhere}
	\]
	and
	\begin{equation}\label{eq:intless1}
		\int_{\R_+\times E}(\beta(t,x)-\alpha(t,x))\,\mathrm{d}t\nu(\mathrm{d}x)<1.
	\end{equation}
	Consequently, all assumptions
        of Theorem~\ref{prop:poincare} are satisfied.
	In particular, we note that
	\begin{align}
		0\leq\frac{\Var\bigl[\pi_{(t,x)}\big| \mathcal F_{t^-}\bigr]}{\mathbb E\bigl[\pi_{(t,x)}\big| \mathcal F_{t^-}\bigr]}&=\frac{\mathbb E\bigl[\pi_{(t,x)}^2\big| \mathcal F_{t^-}\bigr]}
		{\mathbb E\bigl[\pi_{(t,x)}\big| \mathcal F_{t^-}\bigr]}-\mathbb E\bigl[\pi_{(t,x)}\big| \mathcal F_{t^-}\bigr]\nonumber\\
		&\leq\beta(t,x)-\alpha(t,x),\nonumber
	\end{align}
	which, combined with \eqref{eq:intless1}, guarantees \eqref{eq:AssumptionPoincare}.
\end{remark}
\begin{remark}
	\label{rem:KL}
	As mentioned in the introduction, a Poincar\'e inequality for Gibbs point processes was proved in \cite{kondratiev} (see Corollaries 5.1 and 5.2 therein).
	Basically, Corollary 5.1 in \cite{kondratiev} (of which Corollary 5.2 is a small improvement) states and proves the following.
	If $N$ is a grand canonical Gibbs point process on $\R^d$ with activity parameter $z>0$ and non-negative pair potential $\phi$
	such that
	\[
	\delta:=z\int_{\R^d}(1-\mathrm{e}^{-\phi(x)})\,\mathrm{d}x<1,
	\]
	then, for any square-integrable functional $F$ of the point process, we have
	\[
	\mathbb{V}ar (F) \leq (1-\delta)^{-1}\mathbb{E}\biggl[\int_{\mathbb R^d} |D_x F|^2 \pi_x\,dx\biggr],
	\]
	where $D_x$ is the usual add one-cost operator,
	\[
	\pi_x(\omega):=z\exp(-\mathfrak{E}(x,\omega))
	\]
	is the Papangelou conditional intensity of $N$ and $\mathfrak{E}$ denotes the relative energy.
	On the one hand, the Poincar\'e inequality provided by Theorem~\ref{prop:poincare} can be applied e.g. to square-integrable functionals of renewal, non-linear Hawkes and Cox point processes,
	and these processes do not belong to the class of Gibbs point measures for which the inequality in \cite{kondratiev} applies.
	On the other hand, we were not able to apply our Poincar\'e inequality to the Gibbs measures considered in \cite{kondratiev}.
\end{remark}

The next three corollaries, whose proofs are given in Section~\ref{sec:main}, provide classes of non-Poissonian point processes which satisfy the Poincar\'e inequality \eqref{eq:PoincareL2}.

\begin{Corollary}[Renewal point processes]\label{cor:PoincareRenewal}
Let $N$ be a renewal point process on $[0,T]$, $T<\infty$, with points $T_0:=0<T_1<T_2<\cdots<T_{N([0,T])}$
and a spacing density $f$ such that $f$ is continuous on $[0,+\infty)$ and $f>0$ on $(a,C)$ for some $a\in[0,T]$ and $C\in(T,+\infty]$. Assume further that
		\begin{equation}
		\label{eq:ConditionOnSpacing}
			h(z):=\sup_{x\in[z,T]}\frac{f(x-z)}{f(x)}<\infty,
		\end{equation}
		for any $z\in[0,T]$, as well as
		\begin{equation}
		\label{eq:NewBound}
		\gamma:=\int_{0}^{T}\sup_{z\in[0,t]}f(z)\biggl[\frac{\overline F(T-t)^2}
		{\overline F(T-t+z)}
		+h(z)^2\bigl(\overline{F}(z)-\overline{F}(T-t+z)\bigr)-\frac{1}
		{\overline F(z)}\biggr]\,\mathrm{d}t
		<1,
		\end{equation}
		where $\overline{F}$ is the tail function of $f$.
		Then, for any square integrable functional $G$ of $N$,
		\begin{equation}
		\label{eq:PoincareRenewal}
		\Var(G)\le\bigl(1-\sqrt\gamma\bigr)^{-2}\,\mathbb E\biggl[\int_{0}^{T}\pi_{t}\left|D_{t} G\right|^2\,\mathrm dt\biggr],
		\end{equation}
		where the Papangelou conditional intensity $\pi$ of $N$ is given by
		\begin{equation*}
		\pi_t=
		\begin{cases}
		f(T_i-t)f(t-T_{i-1})\bigm/f(T_i-T_{i-1})&\text{ if }T_{i-1}\le t<T_i,\\
		f(t-T_{N([0,T])})\overline{F}(T-t)\bigm/\overline{F}(T-T_{N([0,T])})&\text{ if }T_{N([0,T])}\le t\le T.
		\end{cases}
		\end{equation*}
	\end{Corollary}
	
Next, we give some illustrating examples of renewal point processes which satisfy the assumptions of Corollary~\ref{cor:PoincareRenewal}.

\begin{Example}
	First, assume that the spacing density function $f$ is given by $f(x):=\lambda\exp(-\lambda x)$ for $\lambda>0$, i.e. $N$ is a Poisson process on $[0,T]$ with intensity $\lambda>0$.
	Then all assumptions of Corollary~\ref{cor:PoincareRenewal} are satisfied.
	In particular $h(z):=\exp(\lambda z)$, $\pi_t\equiv\lambda$ and $\gamma=0$.

	Second, consider the Weibull spacing density function
	\[
		f(x):=\beta x^{\beta-1}\exp\bigl(-x^\beta\bigr),\quad\beta>1,\ x\in\mathbb{R}_+.
	\]
	The corresponding tail function is $\overline{F}(x)=\exp\bigl(-x^\beta\bigr)$, $x\in\mathbb{R}_+$, and since $x\mapsto x^\beta$ is Lipschitz continuous on $[0,T]$, for any $z\in[0,T]$ we have
	\begin{equation*}
		h(z)
		=\sup_{x\in[z,T]}\biggl(1-\frac{z}{x}\biggr)^{\beta-1}\exp\bigl(x^\beta-(x-z)^\beta\bigr)
		\le\exp\bigl(\beta T^{\beta-1}z\bigr).
	\end{equation*}
	With the aim to check condition~\eqref{eq:NewBound}, we remark that again by Lipschitz continuity, we have
	\begin{align*}
	\gamma
	&\le\int_{0}^{T}\sup_{z\in[0,t]}\beta z^{\beta-1}\exp\bigl(-z^\beta\bigr)\biggl[\exp\bigl(-2(T-t)^\beta+(T-t+z)^\beta\bigr)\\
	&\qquad\qquad\qquad\qquad\qquad+\exp\bigl(2\beta T^{\beta-1}z\bigr)\Bigl[\exp\bigl(-z^\beta\bigr)-\exp\bigl(-(T-t+z)^\beta\bigr)\Bigr]-\exp\bigl(z^\beta\bigr)\biggr]\,\mathrm{d}t\\
	&\le\int_{0}^{T}\sup_{z\in[0,t]}\beta z^{\beta-1}\exp\bigl(-z^\beta\bigr)\biggl[\exp\bigl(2\beta T^{\beta-1}z-(T-t+z)^\beta\bigr)\\
	&\qquad\qquad\qquad\qquad\qquad+\exp\bigl(2\beta T^{\beta-1}z\bigr)\Bigl[\exp\bigl(-z^\beta\bigr)-\exp\bigl(-(T-t+z)^\beta\bigr)\Bigr]-\exp\bigl(z^\beta\bigr)\biggr]\,\mathrm{d}t\\
		&=\int_0^T\sup_{z\in[0,t]}\beta z^{\beta-1}\Bigl[\exp\bigl(2z\bigl((\beta-1)T^{\beta-1}+T^{\beta-1}-z^{\beta-1}\bigr)\bigr)-1\Bigr]\d t\\
	&\le\int_0^T\sup_{z\in[0,t]}\beta z^{\beta-1}\Bigl[\exp\bigl(2z\bigl((\beta-1)T^{\beta-1}+(\beta-1)T^{\beta-2}(T-z)\bigr)\bigr)-1\Bigr]\d t\\
		&\le\Bigl[\sup_{u\in[0,T]}\exp\bigl(2(\beta-1)T^{\beta-2}u(2T-u)\bigr)-1\Bigr]\int_0^T\sup_{z\in[0,t]}\beta z^{\beta-1}\d t\\
	&=T^\beta\bigl[\exp\bigl(2(\beta-1)T^\beta\bigr)-1\bigr]
	=:\gamma^*(\beta).
	\end{align*}
	Hence for fixed $T>0$ and $\beta_0$ satisfying the inequality $\gamma^*(\beta_0)<1$ (in particular, note that it suffices to take $\beta_0$ close to one) the corresponding renewal point process satisfies the Poincar\'e inequality \eqref{eq:PoincareRenewal}.
	
	Next, consider the generalized Pareto spacing density function
	\[
		f(x):=\lambda(1+\xi\lambda x)^{-(1+1/\xi)},\quad\lambda,\xi>0,\ x\in\mathbb{R}_+,
	\]
	whose tail function is $\overline{F}(x)=(1+\xi\lambda x)^{-1/\xi}$, $x\in\mathbb{R}_+$.
	For any $z\in[0,T]$ we have
	\begin{equation*}
		h(z)
		=\sup_{x\in[z,T]}\biggl(1+\frac{\xi\lambda z}{1+\xi\lambda(x-z)}\biggr)^{1+1/\xi}
		=(1+\xi\lambda z)^{1+1/\xi}.
	\end{equation*}
	We have
	\begin{align*}
	\gamma&=\int_{0}^{T}\sup_{z\in[0,t]}\lambda(1+\xi\lambda z)^{-(1+1/\xi)}\biggl[\frac{(1+\xi\lambda(T-t))^{-2/\xi}}{(1+\xi\lambda(T-t+z))^{-1/\xi}}\\
	&\qquad\qquad\qquad+(1+\xi\lambda z)^{2+2/\xi}\bigl((1+\xi\lambda z)^{-1/\xi}-(1+\xi\lambda(T-t+z))^{-1/\xi}\bigr)-\frac{1}{(1+\xi\lambda z)^{-1/\xi}}\biggr]\,\mathrm{d}t\\
		&=\int_{0}^{T}\sup_{z\in[0,t]}\lambda(1+\xi\lambda z)^{-1}\biggl[\biggl(\frac{1+\frac{\xi\lambda(T-t)}{1+\xi\lambda z}}{1+\xi\lambda(T-t)}\biggr)^{2/\xi}(1+\xi\lambda z)^{1/\xi}(1+\xi\lambda(T-t+z))^{-1/\xi}\\
	&\qquad\qquad\qquad+(1+\xi\lambda z)^{2}\bigl(1-(1+\xi\lambda z)^{1/\xi}(1+\xi\lambda(T-t+z))^{-1/\xi}\bigr)-1\biggr]\,\mathrm{d}t\\
	&\le\int_{0}^{T}\sup_{z\in[0,t]}\lambda(1+\xi\lambda z)^{-1}\biggl[(1+\xi\lambda z)^{1/\xi}(1+\xi\lambda(T-t+z))^{-1/\xi}\\
	&\qquad\qquad\qquad+(1+\xi\lambda z)^{2}\bigl(1-(1+\xi\lambda z)^{1/\xi}(1+\xi\lambda(T-t+z))^{-1/\xi}\bigr)-1\biggr]\,\mathrm{d}t\\
	&=\int_{0}^{T}\sup_{z\in[0,t]}\lambda(1+\xi\lambda z)^{-1}\Bigl(1-(1+\xi\lambda z)^{1/\xi}(1+\xi\lambda(T-t+z))^{-1/\xi}\Bigr)\Bigl((1+\xi\lambda z)^{2}-1\Bigr)\mathrm dt\\
	&\le\lambda\int_0^T\sup_{z\in [0,t]}\frac{\bigl((1+\xi\lambda t)^2-1\bigr)}{1+\xi\lambda z}\mathrm dt
	=\xi\lambda^2T^2\Bigl(1+\frac{\xi\lambda T}{3}\Bigr)
		=:\gamma^*(\xi).
	\end{align*}
	Hence for fixed $T,\lambda>0$ and $\xi_0$ satisfying the inequality $\gamma^*(\xi_0)<1$
(in particular, note that it suffices to take $\xi_0$ close to zero) the corresponding renewal point process satisfies the Poincar\'e inequality \eqref{eq:PoincareRenewal}.
\end{Example}

A nonlinear Hawkes process on $[0,T]$ with parameters $(h,\phi)$ is a point process $N$ on $[0,T]$ with stochastic intensity
\[
	\lambda_t:=\phi\left(\int_{(0,t)}h(t-s)N(\mathrm{d}s)\right),\quad t\in [0,T],
\]
where $\phi:\mathbb R\to\mathbb{R}_+$ and $h:[0,T]\to\R$ are two measurable functions, see e.g. \cite{bremaudmassouliestab}, \cite{bremaudmassoulie}, \cite{daley} and \cite{daley2}.

\begin{Corollary}[Nonlinear Hawkes processes]\label{cor:Hawkes}
	Assume that $N$ is a nonlinear Hawkes process on $[0,T]$, $T<\infty$, with parameters $(h,\phi)$ such that $h$ is non-negative and integrable on $[0,T]$,
and $\phi$ is Lipschitz continuous and non-increasing with $\phi(0)>0$.
	Additionally, assume
	\[
		\gamma:=\phi(0)\int_0^T\Bigl(\exp\Bigl(2\|\phi\|_{\mathrm{Lip}}\int_0^{\tau}h(z)\,\mathrm{d}z\Bigr)-1\Bigr)\,\mathrm{d}\tau<1,
	\]
	where $\|\phi\|_{\mathrm{Lip}}$ denotes the Lipschitz constant of $\phi$.
	Then, for any square-integrable functional $G$ of $N$,
	\begin{equation}
	\label{eq:PoincareHawkes}
	\Var(G)\le\bigl(1-\sqrt{\gamma}\bigr)^{-2}\,\mathbb E\biggl[\int_{0}^{T}\pi_{t}\left|D_{t} G\right|^2\,\mathrm dt\biggr],
	\end{equation}
	where the Papangelou conditional intensity $\pi$ of $N$ is given by
	\begin{equation}\label{eq:PapHaw}
		\pi_t:=\lambda_t\mathrm{E}_t(N)
	\end{equation}
	and
	\begin{align}
		\mathrm{E}_t(N):=&\exp\left(\int_{t}^T\Biggl[\phi\left(\int_{(0,s)}h(s-u)N(\mathrm{d}u)\right)-\phi\left(h(s-t)+\int_{(0,s)}h(s-u)N(\mathrm{d}u)\right)\Biggr]\,\mathrm{d}s\right)\nonumber\\
		&\times\exp\left(\int_{t}^{T}\Biggl[\log\phi\left(h(s-t)+\int_{(0,s)}h(s-u)N(\mathrm{d}u)\right)-\log\phi\left(\int_{(0,s)}h(s-u)N(\mathrm{d}u)\right)\Biggr]N(\mathrm{d}s)\right).\nonumber
	\end{align}
	In particular, we have
	\begin{equation}\label{eq:PapHawbis}
		\pi_t\leq\phi(0)\exp\left(\|\phi\|_{\mathrm{Lip}}\int_0^{T-t}h(s)\,\mathrm{d}s\right)
		\leq\phi(0)\exp\left(\|\phi\|_{\mathrm{Lip}}\|h\|_{L^1([0,T],\mathcal{B}([0,T]),\mathrm{d}t)}\right),
	\end{equation}
	$\d t\mathbb{P}(\d\omega)$-almost everywhere.
\end{Corollary}

Next, we give an example of nonlinear Hawkes processes which satisfies the assumptions of Corollary~\ref{cor:Hawkes}.

\begin{Example}
	Assume that $N$ is a nonlinear Hawkes process with parameters $\phi(x):=\alpha\min(\max(K-x,0),1)$ and $h:=\ind_{[0,z]}$ for some $\alpha,K,z>0$, $K$ integer.
	This is a notable example of nonlinear Hawkes process since $N([0,t])$, $t\in[0,T]$, is the total number of customers who have entered, in the time interval $[0,t]$, the
    Erlang loss system (or M/D/K/0 queue) with arrival rate $\alpha$, deterministic service time $z$ and number of servers equal to $K$, see \cite{bremaudmassouliestab} for details.
	An easy computation shows that $\gamma$ defined in Corollary~\ref{cor:Hawkes} is equal to
	\begin{equation*}
		\gamma(z):=\alpha\biggl[\frac{\mathrm e^{2\alpha\min(z,T)}}{2\alpha}-\frac{1}{2\alpha}-\min(z,T)+\bigl(\mathrm e^{2\alpha z}-1\bigr)(T-z)\ind_{\{z\le T\}}\biggr].
	\end{equation*}
	Therefore, for fixed $\alpha,T>0$, and $z_0$ satisfying the inequality $\gamma(z_0)<1$ (in particular, note that it suffices to take $z_0$ close to zero),
the corresponding nonlinear Hawkes process satisfies the Poincar\'e inequality \eqref{eq:PoincareHawkes}.
\end{Example}

We conclude this subsection by providing a class of Cox processes which satisfy the Poincar\'e inequality.

\begin{Corollary}[Cox processes]\label{cor:Cox}
	Assume $\nu(E)<\infty$ and let $N$ be a Cox process on $[0,T]\times E$, $T<\infty$, with stochastic intensity $\{\lambda_{(t,x)}\}_{(t,x)\in [0,T]\times E}$ such that, for some non-negative
	functions $\alpha,\beta\in L^1([0,T]\times E,\mathcal{B}([0,T])\otimes\mathcal{E},\mathrm{d}t\nu(\mathrm{d}x))$,
	\begin{equation}
	\label{eq:AssumptionOnLambdaCox}
	\alpha(t,x)\leq\lambda_{(t,x)}\leq\beta(t,x),\quad\text{$\mathrm{d}t\nu(\mathrm{d}x)\mathbb{P}(\mathrm{d}\omega)$-almost everywhere}
	\end{equation}
	and
	\begin{align}
	\gamma:=\|\beta\|_{L^1([0,T]\times E,\mathcal{B}([0,T])\otimes\mathcal{E},\mathrm{d}t\nu(\mathrm{d}x))}-\|\alpha\|_{L^1([0,T]\times E,\mathcal{B}([0,T])\otimes\mathcal{E},\mathrm{d}t\nu(\mathrm{d}x))}<1.\label{eq:intless2}
	\end{align}
	Then, for any square-integrable functional $G$ of $N$ we have
	\begin{equation*}
	\Var(G)\le\bigl(1-\sqrt\gamma\bigr)^{-2}\,\mathbb E\biggl[\int_{[0,T]\times E}\pi_{(t,x)}\left|D_{(t,x)} G\right|^2\,\mathrm dt\nu(\d x)\biggr],
	\end{equation*}
	where the Papangelou conditional intensity $\pi$ of $N$ is given for fixed $\omega\in\Omega$ by
	\begin{equation}\label{eq:PapCox}
	\pi_{(t,x)}(\omega):=
		\mathbb{E}[\mathcal{R}(\lambda,\omega)\lambda_{(t,x)}]
		\end{equation}
	where
	\begin{equation*}
	\mathcal{R}(\lambda,\omega):=\frac{\exp\left(-\int_{[0,T]\times E}\lambda_{(s,z)}\,\mathrm{d}s\nu(\mathrm{d}z)\right)\prod_{(s,z)\in\mathrm{Supp}(\omega)}\lambda_{(s,z)}}
	{\mathbb{E}\left[\exp\left(-\int_{[0,T]\times E}\lambda_{(s,z)}\,\mathrm{d}s\nu(\mathrm{d}z)\right)\prod_{(s,z)\in\mathrm{Supp}(\omega)}\lambda_{(s,z)}\right]}.
	\end{equation*}
		In particular, we have
	\begin{equation*}
	\pi_{(t,x)}(\omega)
	\leq\beta(t,x)
	\end{equation*}
	$\d t\mathbb{P}(\d\omega)\nu(\d x)$-almost everywhere.
\end{Corollary}

\subsection{Transportation cost inequalities}

Let $\chi$ be a Polish space equipped with its Borel $\sigma$-field $\mathcal{B}(\chi)$ and let
$d$ be a lower semi-continuous metric on $\chi$ (which does not necessarily generates the topology on $\chi$).
Letting $\sigma_1$ and $\sigma_2$ denote a couple of probability measures on $(\chi,\mathcal B(\chi))$, we define the transportation cost
\[
	\mathcal T_d(\sigma_1,\sigma_2):=\inf_\rho\int_{\chi\times\chi}d(x,y)\,\rho(\mathrm{d}x\times\mathrm{d}y)\in [0,\infty],
\]
where the infimum is taken over all the probability measures $\rho$ on $\chi\times\chi$ with first marginal $\sigma_1$ and second marginal $\sigma_2$.
We denote by $M_1(\chi,d)$ the set of all probability measures $\sigma$ on $(\chi,\mathcal B(\chi))$ such that $\int d(x_0,x)\,\sigma(\mathrm{d}x)<\infty$, for some $x_0\in\chi$ and we remark that for $\sigma_1,\sigma_2\in M_1(\chi,d)$, $\mathcal T_d(\sigma_1,\sigma_2)<\infty$.
The relative entropy of $\sigma_1$ with respect to $\sigma_2$ is defined by
\[
	H(\sigma_1\mid\sigma_2):=\int_\chi\log\left(\frac{\mathrm{d}\sigma_1}{\mathrm{d}\sigma_2}\right)\mathrm{d}\sigma_1
\]
if $\sigma_1$ is absolutely continuous with respect to $\sigma_2$ and $H(\sigma_1\mid\sigma_2):=+\infty$ otherwise (the reader is referred to \cite{villani} for more insight into the theory of optimal transportation).

In the following, we denote by
\[
	g^{\odot}(x):=\sup_{\theta\geq 0}(\theta x-g(\theta)),\quad x\in\mathbb{R}_+
\]
the monotone conjugate of a measurable function $g:\mathbb{R}_+\to [0,\infty]$.

\subsubsection{A transportation cost inequality for functionals of $N$}\label{sec:transpofunc}

In this subsection we take $\chi:=\mathbb R$, suppose that there exists a norm on $\chi$, say $\|\cdot\|_d$,
such that $d(x,y)=\|x-y\|_d$, $x,y\in\mathbb R$, and denote by $\mathcal{L}(X)$ the law of a real-valued random variable $X$. The following theorem holds.

\begin{Theorem}\label{prop:transport}
Suppose that $N$ satisfies \eqref{eq:LocalStability} for a measurable function $\beta$ and let $G$ be an integrable random variable such that
\[
	\|D_{(t,x)}G(\omega)\|_d\leq g_1(t,x)\quad\text{and}\quad\frac{\mathbb{E}\bigl[\|G-\mathbb E[G\mid\mathcal F_{t^-}]\|_d|\pi_{(t,x)}-p(\pi)_{(t,x)}|\mid \mathcal{F}_{t^-}\bigr](\omega)}{p(\pi)_{(t,x)}(\omega)}\leq g_2(t,x)
\]
for $\d t\d\mathbb P\d\nu$-almost all $(t,\omega,x)$ and some deterministic functions $g_1,g_2$ such that
\begin{equation}\label{eq:inth}
	\int_{\mathbb R_+\times E}|h(t,x)|^2\beta(t,x)\,\mathrm{d}t\nu(\d x)<\infty
\end{equation}
where $h:=g_1+g_2$. Then
\begin{equation}\label{eq:traspineq}
	c(\mathcal T_d(\sigma,\mathcal{L}(G)))\leq H(\sigma\mid\mathcal{L}(G)),\quad\text{for any $\sigma\in M_1(\R,d)$}
\end{equation}
where
\begin{equation}
\label{eq:DefinitionOfC}
	c(x):=\Lambda^\odot(x),\quad x\in\R_+
\end{equation}
and
\begin{equation}\label{eq:Lambda}
	\Lambda(\theta):=\int_{\R_+\times E}\bigl(\mathrm{e}^{\theta h(t,x)}-\theta h(t,x)-1\bigr)\beta(t,x)\,\d t\nu(\d x),\quad\theta\in\R_+.
\end{equation}
If additionally we assume that there exists $M>0$ such that $h(t,x)\le M$ for $\d t\d\nu$-almost all $(t,x)$, then
\begin{equation}
\label{eq:LowerBoundOfC1}
	c(x)\geq\widetilde{c}(x),\quad x\in\R_+,
\end{equation}
where
\begin{equation*}
	\widetilde{c}(x):=\frac{x+M^{-1}\int_{\R_+\times E}h(t,z)^2\beta(t,z)\,\d t\nu(\d z)}{M}\ln\biggl(1+\frac{x}{M^{-1}\int_{\R_+\times E}h(t,z)^2\beta(t,z)\,\d t\nu(\d z)}\biggr)-\frac{x}{M}.
\end{equation*}
\end{Theorem}

\begin{remark}\label{re:pi/ppi0}
	As noticed in Remark~\ref{re:poincare}, if $\mathbb E\bigl[\pi_{(t,x)}\big| \mathcal F_{t^-}\bigr]=0$ $\mathbb P$-almost surely, then $\mathbb P(\pi_{(t,x)}=0\mid\mathcal F_{t^-})=1$ $\mathbb P$-almost surely.
	Consequently, the ratio
	\[
		\frac{\mathbb{E}\bigl[\|G-\mathbb E[G\mid\mathcal F_{t^-}]\|_d|\pi_{(t,x)}-p(\pi)_{(t,x)}|\big| \mathcal{F}_{t^-}\bigr]}{p(\pi)_{(t,x)}}
	\]
	is always well-defined under the convention $0/0:=0$.
\end{remark}

\begin{remark}[Deviation inequality]\label{re:deviationbound}
	From the point of view of the applications, it is important to remark that a transportation cost inequality is often equivalent to a deviation bound.
	More precisely, in the context of Theorem \ref{prop:transport}, one has that the transportation cost inequality \eqref{eq:traspineq} is equivalent to the deviation bound
	\begin{equation}\label{eq:devineq1}
		\mathbb{P}\left(\frac{1}{n}\sum_{i=1}^{n}f(G_i)\geq\mathbb E\bigl[f(G)\bigr]+r\right)\leq\mathrm{e}^{-n c(r)},\quad\text{for any $n\geq 1$ and $r\in\mathbb{R}_+$}
	\end{equation}
	for any measurable function $f:\mathbb R\to\mathbb R$ which is Lipschitz continuous (with respect to $d$) with Lipschitz constant less than or equal to $1$, i.e. such that
	\[
		\sup_{x\neq y}\frac{|f(x)-f(y)|}{d(x,y)}\leq 1,
	\]
	where $\{G_n\}_{n\geq 1}$ is a sequence of independent random variables with the same law as $G$ (see \cite{gozlan} and Theorem~1.1(c) in \cite{MaWu}).
\end{remark}

To the best of our knowledge, the transportation cost inequality provided by Theorem~\ref{prop:transport} is new even in the Poisson case, which we state in a separate corollary.

\begin{Corollary}[Poisson processes]\label{cor:transportPoiss}
	Suppose that $N$ is a Poisson process on $\R_+\times E$ with mean measure $\beta(t,x)\,\d t\nu(\d x)$ and let $G$ be an integrable random variable such that
	\[
		\quad\|D_{(t,x)}G(\omega)\|_d\leq h(t,x)
	\]
	for $\d t\d\mathbb P\d\nu$-almost all $(t,\omega,x)$ and some deterministic function $h$ which satisfies \eqref{eq:inth}.
	Then the transportation cost inequality \eqref{eq:traspineq} holds. If additionally we assume that there exists $M>0$ such that $h(t,x)\le M$ for $\d t\d\nu$-almost all $(t,x)$,
	then \eqref{eq:LowerBoundOfC1} holds, and provides a more explicit bound on the corresponding deviation inequality \eqref{eq:devineq1}.
\end{Corollary}

In the following proposition, we specialize Theorem~\ref{prop:transport} to first order integrals.
The subsequent corollaries concern applications to renewal, nonlinear Hawkes and Cox point processes.

\begin{Proposition}[First order integrals]\label{prop:transportFirstOrder}
	Suppose that $N$ satisfies \eqref{eq:LocalStability} for a measurable function $\beta$,
	$d(x,y):=|x-y|$, take
	\begin{equation*}
	\label{eq:GFirstOrder}
		G:=\int_{\R_+\times E}g(s,y)\,N(\mathrm ds\times\mathrm dy),
	\end{equation*}
	for some measurable deterministic function $g$, and let $g_1,g_2$ be deterministic functions such that
	\begin{equation*}
		|g(t,x)|\le g_1(t,x)
	\end{equation*}
	and
	\begin{multline}
	\label{eq:g2Bounds}
		\sqrt2\biggl[\int_{[t,\infty)\times E}g(s,y)^2\beta(s,y)\,\mathrm ds\nu(\mathrm dy)
		+\biggl(\int_{[t,\infty)\times E}|g(s,y)|\beta(s,y)\,\mathrm ds\nu(\mathrm dy)\biggr)^2\biggr]^{1/2}\\
		\times\Biggl\|\frac{\sqrt{\mathbb V\mathrm{ar}\bigl[\pi_{(t,x)}\mid\mathcal F_{t^-}\bigr]}}{p(\pi)_{(t,x)}}\Biggr\|_{L^\infty(\Omega,\mathcal F_\infty,\mathbb P)}\le g_2(t,x),
	\end{multline}
	for $(t,x)\in\R_+\times E$.
	Assuming that $h:=g_1+g_2$ satisfies the corresponding assumption \eqref{eq:inth}, the transportation cost inequality \eqref{eq:traspineq} holds.
	If additionally we assume that there exists $M>0$ such that $h(t,x)\le M$ for $\d t\d\nu$-almost all $(t,x)$, then \eqref{eq:LowerBoundOfC1} holds, and provides a more explicit bound on the corresponding deviation inequality \eqref{eq:devineq1}.
\end{Proposition}

In the next corollaries, we provide classes of point processes which satisfy the assumptions of Proposition~\ref{prop:transportFirstOrder}.
We point out that we do not aim to optimize the assumptions, favoring instead clarity and conciseness.
The proofs are quite elementary, and provided in Section~\ref{subsec:proofOfApplicationsCorollaries} for the readers' convenience.

\begin{Corollary}[Renewal point processes]\label{cor:transportRenewal}
	Let $N$ be a renewal point process on $[0,T]$, $T<\infty$, with
		a spacing density $f$ such that $f$ is continuous on $[0,+\infty)$ and $f>0$ on $(a,C)$ for some $a\in[0,T]$ and $C\in(T,+\infty]$, and let $G$ be defined by
	\begin{equation}
	\label{eq:GFirstOrderSquareIntegrable}
		G:=\int_0^Tg(s)\,N(\mathrm ds),\quad g\in L^2([0,T]).
	\end{equation}
Assume further that there exists $\overline h>0$ such that
	\begin{equation}
	\label{eq:AssumptionOnHBar}
		\sup_{z\in[0,T]}\sup_{x\in[z,T]}\frac{f(x-z)}{f(x)}\le\overline h
	\end{equation}
	(note that this condition is always satisfied if $f>0$ on $[0,T]$.)
	After straightforward adjustments due to the unmarked setting, the transportation cost inequality \eqref{eq:traspineq} holds with $g_1\in L^2([0,T])$ such that
		$|g(t)|\leq g_1(t)$, $t\in [0,T]$,
		and $g_2\in L^2([0,T])$ such that
	\begin{equation}
	\label{eq:DefinG2Renewal}
		\sqrt{2\left(\overline h^2+\biggl(\int_T^Cf(x)\,\mathrm dx\biggr)^{-2}-1\right)}\biggl[\beta\int_t^T g(s)^2\,\mathrm ds
		+\beta^2\biggl(\int_t^T|g(s)|\,\mathrm ds\biggr)^2\biggr]^{1/2}
		\le g_2(t),\quad t\in [0,T]
	\end{equation}
	where
	\begin{equation}\label{eq:betaren}
		\beta:=\max\biggl(\frac{\sup_{x\in [0,T]}f(x)}{\int_T^C f(x)\,\mathrm{d}x},\sup_{x\in[0,C)}f(x)+\frac{\bigl(\sup_{x\in [0,T]}f(x)\bigr)^2}{\min_{x\in[a,T]}f(x)}\biggr).
	\end{equation}
	For any $M>0$ such that $\|g_1+g_2\|_\infty\leq M$, we have that \eqref{eq:LowerBoundOfC1} holds, which yields a more explicit bound on the corresponding deviation inequality \eqref{eq:devineq1}.
\end{Corollary}

\begin{Corollary}[Nonlinear Hawkes processes]\label{cor:transportHawkes}
	Assume that $N$ is a nonlinear Hawkes process on $[0,T]$, $T<\infty$, with parameters $(h,\phi)$ such that $h$ is non-negative and integrable on $[0,T]$,
	$\phi$ is Lipschitz continuous and non-increasing with $\phi(0)>0$, and let $G$ be defined by
		\eqref{eq:GFirstOrderSquareIntegrable}.
		After straightforward adjustments due to the unmarked setting, the transportation cost inequality \eqref{eq:traspineq} holds with $g_1\in L^2([0,T])$ such that
			$|g(t)|\leq g_1(t)$, $t\in [0,T]$,
		and $g_2\in L^2([0,T])$ such that
	\begin{multline}
	\label{eq:DefinG2Hawkes}
		g_2(t)\ge\sqrt 2\biggl[\phi(0)\int_t^T g(s)^2\exp\biggl(\|\phi\|_{\mathrm{Lip}}\int_0^{T-s}h(z)\,\mathrm dz\biggr)\,\mathrm{d}s
		\\
		+\biggl(\phi(0)\int_t^T|g(s)|\exp\biggl(\|\phi\|_{\mathrm{Lip}}\int_0^{T-s}h(z)\,\mathrm dz\biggr)\,\mathrm ds\biggr)^2\biggr]^{1/2}
		\sqrt{\exp\biggl(2\|\phi\|_{\mathrm{Lip}}\int_0^{T-t}h(z)\,\mathrm dz\biggr)-1},
	\end{multline}
	$t\in [0,T]$. For any $M>0$ such that $\|g_1+g_2\|_\infty\leq M$, we have that \eqref{eq:LowerBoundOfC1} holds, which yields a more explicit bound on the corresponding deviation inequality \eqref{eq:devineq1}.
\end{Corollary}

\begin{Corollary}[Cox processes]\label{cor:transportCox}
	Assume $\nu(E)<\infty$ and let $N$ be a Cox process on $[0,T]\times E$, $T<\infty$, with stochastic intensity $\{\lambda_{(t,x)}\}_{(t,x)\in [0,T]\times E}$ satisfying \eqref{eq:AssumptionOnLambdaCox}, for some non-negative functions $\alpha,\beta$ with $\beta\in L^1([0,T]\times E,\mathrm{d}t\nu(\mathrm{d}x))$ and $\beta^3\alpha^{-2}\in L^1([0,T]\times E,\mathrm dt\nu(\mathrm dx))$.
	If $G$ is defined by
	\begin{equation*}
	G:=\int_{[0,T]\times E}g(s,y)\,N(\mathrm ds\times\mathrm dy),\quad g\in L^2([0,T]\times E,\beta(t,x)\mathrm dt\nu(\mathrm dx)),
	\end{equation*}
	then the transportation cost inequality \eqref{eq:traspineq} holds with $g_1\in L^2([0,T]\times E,\beta(t,x)\mathrm dt\nu(\mathrm dx))$ such that
			$|g(t,x)|\leq g_1(t,x)$, $(t,x)\in [0,T]\times E$,
		and $g_2\in L^2([0,T]\times E,\beta(t,x)\mathrm dt\nu(\mathrm dx))$ such that
	\begin{multline*}
		\sqrt2\biggl[\int_{[t,T]\times E}g(s,y)^2\beta(s,y)\,\mathrm ds\nu(\mathrm dy)
		+\biggl(\int_{[t,T]\times E}|g(s,y)|\beta(s,y)\,\mathrm ds\nu(\mathrm dy)\biggr)^2\biggr]^{1/2}\\
		\times\sqrt{\frac{\beta(t,x)^2}{\alpha(t,x)^2}-1}\le g_2(t,x),
	\end{multline*}
	for $(t,x)\in[0,T]\times E$.
	For any $M>0$ such that $\|g_1+g_2\|_\infty\leq M$, we have that \eqref{eq:LowerBoundOfC1} holds, which yields a more explicit bound on the corresponding deviation inequality \eqref{eq:devineq1}.
\end{Corollary}

\subsubsection{A transportation cost inequality for the law of $N$}

In this subsection we take $\chi:=\Omega$ and equip this set with the vague convergence topology, that is, the coarsest topology such that the map $\omega\mapsto\int_{\mathbb R_+\times E} f(t,x)\,\omega(\d t\times\d x)$ is continuous, where $f:\R_+\times E\to\R$ is a continuous function with compact support.
It is well-known that this topology makes $\Omega$ a Polish space (see e.g. \cite{daley}).
In this subsection, we let $\varphi:\R_+\times E\to\R_+$ be a continuous function and define the following metric on $\Omega$:
\[
	d_\varphi(\omega,\omega'):=\int_{\R_+\times E}\varphi(t,x)|\omega-\omega'|(\d t\times\d x),\quad\omega,\omega'\in\Omega
\]
where, for $\kappa\in\Omega$, $|\kappa|:=\kappa^++\kappa^-$, and $\kappa^+$ and $\kappa^-$ denote respectively the positive and the negative parts of $\kappa$ in the Hahn-Jordan decomposition.
It is known that $d_\varphi$ is a lower semi-continuous metric on $\Omega$ (see Lemma 2.2 in \cite{MaWu}).
The following theorem holds.

\begin{Theorem}\label{thm:transportlawpp}
	Assume that $N$ satisfies \eqref{eq:LocalStability} for a measurable function $\beta$, that there exists a deterministic measurable function $\psi$ which verifies \eqref{eq:g2Bounds} with $\varphi$ in place of $g$ and $\psi$ in place of $g_2$ and that
	\begin{equation}\label{eq:ExpVarphi}
		\int_{\R_+\times E}|h_{\varphi}(t,x)|^2\beta(t,x)\,\d t \nu(\d x)<\infty,
	\end{equation}
	where $h_\varphi:=\varphi+\psi$.
	Then
	\begin{equation}\label{eq:transportomega}
		c(\mathcal{T}_{d_\varphi}(\mathbb Q,\mathbb P))\leq H(\mathbb Q\mid\mathbb P),\quad\text{for any $\mathbb Q\in M_1(\Omega,d_\varphi)$}
	\end{equation}
	where
	\begin{equation*}
		c(x):=\Lambda_\varphi^\odot(x),\quad x\in\mathbb{R}_+
	\end{equation*}
	and
	\begin{equation*}
		\Lambda_\varphi(\theta):=\int_{\R_+\times E}\bigl(\mathrm{e}^{\theta h_\varphi(t,x)}-\theta h_\varphi(t,x)-1\bigr)\beta(t,x)\,\d t\nu(\d x),\quad\theta\in\mathbb{R}_+.
	\end{equation*}
	In particular, when $N$ is a Poisson process with mean measure $\beta(t,x)\,\d t\nu(\d x)$ we recover the sharp transportation cost inequality in Remark~2.7 of \cite{MaWu}.
\end{Theorem}

\begin{remark}[Deviation inequality]\label{re:deviationboundomega}
	Here again, it turns out (see \cite{gozlan} and Theorem~1.1(c) in \cite{MaWu}) that, in the context of Theorem~\ref{thm:transportlawpp}, the transportation cost inequality \eqref{eq:transportomega}
is equivalent to the deviation bound
	\begin{equation}
	\label{eq:DeviationBoundN}
		\mathbb{P}\left(\frac{1}{n}\sum_{i=1}^{n}F(N_i)\geq\mathbb E\bigl[F(N)\bigr]+r\right)\leq\mathrm{e}^{-n c(r)},\quad\text{for any $n\geq 1$ and $r\in\mathbb{R}_+$}
	\end{equation}
	for any measurable function $F:\Omega\to\mathbb R$ which is Lipschitz continuous (with respect to $d_\varphi$) with Lipschitz constant less than or equal to $1$, i.e. such that
	\[
		\sup_{\omega\neq\omega'}\frac{|f(\omega)-f(\omega')|}{d_\varphi(\omega,\omega')}\leq 1,
	\]
	where $\{N_n\}_{n\geq 1}$, $N_n:\Omega\to\Omega$, is a sequence of independent marked point processes with the same law as $N$.
	As in the previous Subsection \ref{sec:transpofunc}, if additionally we assume
that there exists $M>0$ such that $h_\varphi(t,x)\le M$ for $\d t\d\nu$-almost all $(t,x)$,
then \eqref{eq:LowerBoundOfC1} with $h_\varphi$ in place of $h$ holds,
and provides a more explicit bound on the deviation inequality \eqref{eq:DeviationBoundN}.
\end{remark}

Arguing as in the proofs of Corollaries \ref{cor:transportRenewal}, \ref{cor:transportHawkes} and \ref{cor:transportCox},
one can show that Theorem~\ref{thm:transportlawpp} and the deviation bound in Remark~\ref{re:deviationboundomega} apply to renewal, nonlinear Hawkes and Cox point processes.
Similarly to the previous Subsection \ref{sec:transpofunc}, we
do not aim to optimize the assumptions, favoring instead clarity and conciseness.
\begin{Corollary}[Renewal point processes]\label{cor:transportRenewalForLawN}
	Let $N$ be a renewal point process on $[0,T]$, $T<\infty$, with a spacing density $f$ such that $f$ is continuous on $[0,+\infty)$ and $f>0$ on $(a,C)$ for some $a\in[0,T]$ and $C\in(T,+\infty]$.
	Assume further that
there exists $\overline h>0$ satisfying \eqref{eq:AssumptionOnHBar}.
	Then, after straightforward adjustments due to the unmarked setting, the transportation cost inequality \eqref{eq:transportomega} holds with
	$\psi\in L^2([0,T])$ such that
	\begin{equation*}
		\sqrt{2\left(\overline h^2+\biggl(\int_T^Cf(x)\,\mathrm dx\biggr)^{-2}-1\right)}\biggl[\beta\int_t^T\varphi(s)^2\,\mathrm ds
		+\beta^2\biggl(\int_t^T|\varphi(s)|\,\mathrm ds\biggr)^2\biggr]^{1/2}
		\le\psi(t),\quad t\in [0,T]
	\end{equation*}
where $\beta$ is defined by \eqref{eq:betaren}.
\end{Corollary}
\begin{Corollary}[Nonlinear Hawkes processes]\label{cor:transportHawkesForLawN}
	Assume that $N$ is a nonlinear Hawkes process on $[0,T]$, $T<\infty$, with parameters $(h,\phi)$ such that $h$ is non-negative and integrable on $[0,T]$,
	$\phi$ is Lipschitz continuous and non-increasing with $\phi(0)>0$.
	Then, after straightforward adjustments due to the unmarked setting, the transportation cost inequality \eqref{eq:transportomega} holds with
	$\psi\in L^2([0,T])$ such that
\begin{multline*}
		\psi(t)\ge\sqrt 2\biggl[\phi(0)\int_t^T\varphi(s)^2\exp\biggl(\|\phi\|_{\mathrm{Lip}}\int_0^{T-s}h(z)\,\mathrm dz\biggr)\,\mathrm{d}s
		\\
		+\biggl(\phi(0)\int_t^T|\varphi(s)|\exp\biggl(\|\phi\|_{\mathrm{Lip}}\int_0^{T-s}h(z)\,\mathrm dz\biggr)\,\mathrm ds\biggr)^2\biggr]^{1/2}
		\sqrt{\exp\biggl(2\|\phi\|_{\mathrm{Lip}}\int_0^{T-t}h(z)\,\mathrm dz\biggr)-1},
\end{multline*}
$t\in [0,T]$.
\end{Corollary}
\begin{Corollary}[Cox point processes]\label{cor:transportCoxForLawN}
	Assume $\nu(E)<\infty$ and let $N$ be a Cox process on $[0,T]\times E$, $T<\infty$, with stochastic intensity $\{\lambda_{(t,x)}\}_{(t,x)\in [0,T]\times E}$ satisfying \eqref{eq:AssumptionOnLambdaCox}, for some non-negative functions $\alpha,\beta$ with $\beta\in L^1([0,T]\times E,\mathrm{d}t\nu(\mathrm{d}x))$ and $\beta^3\alpha^{-2}\in L^1([0,T]\times E,\mathrm dt\nu(\mathrm dx))$.
	Assume further that $\varphi\in L^2([0,T]\times E,\beta(t,x)\mathrm{d}t\nu(\mathrm{d}x))$.
	Then, the transportation cost inequality \eqref{eq:transportomega} holds with
	$\psi\in L^2([0,T]\times E,\beta(t,x)\mathrm{d}t\nu(\mathrm{d}x))$ such that
\begin{multline*}
		\sqrt2\biggl[\int_{[t,T]\times E}\varphi(s,y)^2\beta(s,y)\,\mathrm ds\nu(\mathrm dy)
		+\biggl(\int_{[t,T]\times E}|\varphi(s,y)|\beta(s,y)\,\mathrm ds\nu(\mathrm dy)\biggr)^2\biggr]^{1/2}\\
		\times\sqrt{\frac{\beta(t,x)^2}{\alpha(t,x)^2}-1}\le\psi(t,x),
	\end{multline*}
	for $(t,x)\in[0,T]\times E$.
\end{Corollary}

\subsection{Variational representation of the Laplace transform}

The following Theorem~\ref{thm:laplace} generalizes to functionals of marked point processes with Papangelou conditional intensity
the variational representation formula for the Laplace transform of Poisson functionals given by Theorem~4.4 of \cite{zhang1}.
Note also that, in contrast to \cite{zhang1}, here the marked point process is defined on the whole half-line.

Hereafter, we suppose that $N$ satisfies \eqref{eq:squareint} and \eqref{eq:IntegrabilityOfPi}, and denote by $\mathcal{H}$ the subset of $\mathcal{P}_2^{\mathcal F}(p(\pi))$ formed by the real-valued processes $\phi\in L^\infty(\R_+\times\Omega\times E,\mathcal{B}(\R_+)\otimes\mathcal{F}_\infty\otimes\mathcal{E},\d t\d\mathbb P\d\nu)$ such that
\begin{equation}\label{hyp:phi1}
	\phi_{(t,x)}\geq c_\phi>-1,\quad\text{for some constant $c_\phi$, $\d t\d\mathbb P\nu(\d x)$-almost surely},
\end{equation}
\begin{equation}
\label{eq:BoundOnPhiFromMathcalH}
	\mathbb E\biggl[\biggl(\int_{\R_+\times E}|\phi_{(s,x)}|^2\,p(\pi)_{(s,x)}\,\d s\nu(\d x)\biggr)^2\biggr]<\infty
\end{equation}
and
\begin{equation}\label{eq:Ep2}
	\text{for any }T>0,\ (t,x)\mapsto\ind_{[0,T]}(t)\mathcal{E}_t(\phi)\in\mathcal P_2^{\mathcal F}(p(\pi)),
\end{equation}
where
\begin{multline*}
	\mathcal{E}_t(\phi):=\exp\biggl(\int_{[0,t]\times E}\log\bigl(1+\phi_{(s,x)}\bigr)(N(\d s\times\d x)-p(\pi)_{(s,x)}\,\d s\nu(\d x))\\
	+\int_{[0,t]\times E}\bigl(\log\bigl(1+\phi_{(s,x)}\bigr)-\phi_{(s,x)}\bigr)p(\pi)_{(s,x)}\,\d s\nu(\d x)\biggr),\quad t\in\R_+.
\end{multline*}

The following lemma ensures that $\{\mathcal{E}_t(\phi)\}_{t\in\R_+}$ is a square-integrable $\mathcal{F}$-martingale.
Its proof is postponed to Section~\ref{sec:main} (see Subsection~\ref{subsec:lap}).

\begin{Lemma}\label{le:mgproperty}
	Assume that $N$ satisfies \eqref{eq:squareint} and \eqref{eq:IntegrabilityOfPi}.
	Then, for any $\phi\in\mathcal{H}$, the stochastic process $\{\mathcal{E}_t(\phi)\}_{t\in\R_+}$ is a square-integrable $\mathcal{F}$-martingale.
\end{Lemma}

Under the assumptions of Lemma~\ref{le:mgproperty}, for $\phi\in\mathcal H$, we define a new probability measure $\mathbb P_\phi$ on $(\Omega,\mathcal{F}_\infty)$ by
\begin{equation}\label{eq:DefinitionPPhi}
	\frac{\d\mathbb P_{\phi}}{\d\mathbb P}\bigg|_{\mathcal F_t}=\mathcal{E}_t(\phi),\quad t\in\R_+.
\end{equation}

The following variational representation of the Laplace transform holds.

\begin{Theorem}\label{thm:laplace}
Suppose that $N$ satisfies \eqref{eq:squareint}, \eqref{eq:IntegrabilityOfPi},
\begin{equation}\label{eq:Kgrande}
K:=\int_0^\infty\biggl\|\int_Ep(\pi)_{(t,x)}\,\nu(\d x)\biggr\|_{L^\infty(\Omega,\mathcal F_\infty,\mathbb P)}\d t<\infty,
\end{equation}
and let $G$ be a random variable on $(\Omega,\mathcal{F}_\infty,\mathbb P)$ which is upper bounded.
Then
\[
	-\log\bigl(\mathbb E\bigl[\mathrm e^{-G}\bigr]\bigr)=\inf_{\phi\in\mathcal H}\mathbb E_\phi\bigl[G+L(\phi)\bigr],
\]
where $\mathbb E_\phi$ denotes the expectation under $\mathbb P_\phi$
and
\begin{equation*}
	L(\phi):=\int_{\R_+\times E}\bigl((1+\phi_{(s,x)})\log\bigl(1+\phi_{(s,x)}\bigr)-\phi_{(s,x)}\bigr)\,p(\pi)_{(s,x)}\,\d s\nu(\d x).
\end{equation*}
Additionally, if $G\in L^\infty(\Omega,\mathcal{F}_\infty,\mathbb P)$, then the infimum is uniquely attained at
\[
	\phi_{(t,x)}^{(F)}:=\frac{\varphi_{(t,x)}^{(F)}}{p(F)_{(t,x)}},\quad (t,x)\in\R_+\times E
\]
where $F:=\mathrm{e}^{-G}$ and
\begin{equation}\label{eq:ExplicitIntegrand}
	\varphi_{(t,x)}^{(F)}:=\frac{p(\pi F^+)_{(t,x)}-p(F)_{(t,x)}p(\pi)_{(t,x)}}{p(\pi)_{(t,x)}},\quad\text{$(t,x)\in\R_+\times E$.}
\end{equation}
\end{Theorem}
Here, we limit ourselves to note that $\varphi^{(F)}$ is well-defined and belongs to $\mathcal{P}_2^{\mathcal{F}}(p(\pi))$,
and refer the reader to Theorem~\ref{thm:clarkoperative} and Remark~\ref{rem:RemarkFollowingCOGeneralized} for details.

\begin{remark}\label{re:variationalrep}
Assumptions \eqref{eq:squareint}, \eqref{eq:IntegrabilityOfPi} and \eqref{eq:Kgrande} are satisfied if the inequality \eqref{eq:LocalStability} holds with
$\beta\in L^1(\R_+\times E,\mathcal{B}(\R_+)\otimes\mathcal{E},\d t\nu(\d x))$, see also Remark~\ref{rem:AssumptionsAreSatisfiedIf}.
As such, Theorem \ref{thm:laplace} applies to $(i)$ renewal point processes on $[0,T]$, $T<\infty$, with a spacing density $f$ such that $f$ is continuous
on $[0,+\infty)$ and $f>0$ on $(a,C)$ for some $a\in[0,T]$ and $C\in(T,+\infty]$; $(ii)$ nonlinear Hawkes processes on $[0,T]$, $T<\infty$, with parameters $(h,\phi)$
such that $h$ is non-negative and integrable on $[0,T]$, $\phi$ is Lipschitz continuous and non-increasing with $\phi(0)>0$; $(iii)$
if $\nu(E)<\infty$, Cox processes on $[0,T]\times E$, $T<\infty$, with stochastic intensity $\{\lambda_{(t,x)}\}_{(t,x)\in [0,T]\times E}$ such that
$\lambda_{(t,x)}\leq\beta(t,x)$, $(t,x)\in [0,T]\times E$,
with $\beta\in L^1([0,T]\times E,\mathrm{d}t\nu(\mathrm{d}x))$.
\end{remark}

\subsection{Clark-Ocone formula}\label{sec:Operative}

As already mentioned, the proofs of Theorems \ref{prop:poincare}, \ref{prop:transport}, \ref{thm:transportlawpp}
and \ref{thm:laplace} are based on a Clark-Ocone formula for marked point processes with Papangelou conditional intensity which generalizes the one derived in \cite{COFormula} (see Remark~\ref{re:genCO}).

\begin{Theorem}\label{thm:clarkoperative}
Suppose that $N$ satisfies \eqref{eq:squareint} and \eqref{eq:IntegrabilityOfPi}.
Then, for any $G\in L^2(\Omega,\mathcal{F}_\infty,\mathbb P)$,
\begin{equation}\label{eq:clarkoperative2}
	G=\mathbb{E}[G]+\int_{\R_+\times E}\varphi_{(t,x)}^{(G)}\,(N(\d t\times\d x)-p(\pi)_{(t,x)}\,\d t\nu(\d x)),\quad\text{$\mathbb{P}$-almost surely},
\end{equation}
where $\varphi^{(G)}$ is defined by \eqref{eq:ExplicitIntegrand} with $G$ in place of $F$.
Additionally, $\varphi^{(G)}\in\mathcal P_2^{\mathcal{F}}(p(\pi))$.
\end{Theorem}

\begin{remark}\label{re:genCO}
	The Clark-Ocone formula \eqref{eq:clarkoperative2} generalizes the corresponding formula in \cite{COFormula} to point processes on $\R_+$ with marks in $E$ and guarantees that the integrand $\varphi^{(G)}$ is square integrable with respect to $p(\pi)_{(t,x)}(\omega)\d t\mathbb P(\d\omega)\nu(\d x)$.
	        This integrability property is crucial in our proofs since it allows one to apply the isometry formula of Proposition \ref{prop:ExtendedIsometry}.
	We also remark that the proof of formula \eqref{eq:clarkoperative2} provided in this paper is shorter than the proof of the corresponding Clark-Ocone formula in \cite{COFormula}.
\end{remark}

\begin{remark}\label{rem:RemarkFollowingCOGeneralized}
	Under the assumptions of Theorem~\ref{thm:clarkoperative}, by Proposition~\ref{cor:PredictableProjection} we have that $p(G)_{(t,x)}$ exists and it is finite $\mathbb P$-almost surely.
	Additionally, we have $\pi_{(t,x)}(N)G_{(t,x)}^+(N)\in L^1(\Omega,\mathcal{F}_\infty,\mathbb P)$,
	for $\d t\nu(\d x)$-almost all $(t,x)$. Indeed, on any interval $[t_1,t_2]\subset\R_+$ we have
	\begin{align*}
		\int_{[t_1,t_2]\times E}\mathbb E\Bigl[\Bigl|\pi_{(t,x)}(N)G_{(t,x)}^+(N)\Bigr|\Bigr]\d t\nu(\d x)
		&=\mathbb E\left[\int_{\R_+\times E}\left|G(N)\right|\ind_{[t_1,t_2]}(t)\,N(\d t\times\d x)\right]\\
		&\le\|G\|_{L^2(\Omega,\mathcal{F}_\infty,\mathbb P)}\|N([t_1,t_2]\times E)\|_{L^2(\Omega,\mathcal{F}_\infty,\mathbb P)}
		<\infty.
	\end{align*}
	Here, the first equality follows from \eqref{eq:GNZ}.
	Therefore, again by Proposition~\ref{cor:PredictableProjection}, for $\d t\nu(\d x)$-almost all $(t,x)$, $p(\pi G^+)_{(t,x)}$ exists and it is finite $\mathbb P$-almost surely.
	We also note that the difference $p(\pi G^+)_{(t,x)}-p(G)_{(t,x)}p(\pi)_{(t,x)}$ is well-defined and finite $\mathbb P$-almost surely. Indeed, for any $T\geq 0$,
	\begin{equation}
	\label{eq:integralOfpPiFinite}
		\int_{[0,T]\times E}\mathbb E\bigl[p(\pi)_{(t,x)}\bigr]\,\d t\nu(\d x)=\mathbb E\left[\int_{[0,T]\times E}\pi_{(t,x)}\,\d t\nu(\d x)\right]=\mathbb E\bigl[N([0,T]\times E)\bigr]<\infty.
	\end{equation}
	So, for $\d t\nu(\d x)$-almost all $(t,x)$, $p(\pi)_{(t,x)}$ is finite $\mathbb P$-almost surely.
	As noticed in Remark~\ref{re:poincare} if $p(\pi)_{(t,x)}=0$ $\mathbb P$-almost surely then $\mathbb P(\pi_{(t,x)}=0\mid\mathcal{F}_{t^-} )=1$ $\mathbb P$-almost surely.
	Consequently, $p(\pi G^+)_{(t,x)}^2=\mathbb E\bigl[\pi_{(t,x)}G_{(t,x)}^+\big| \mathcal{F}_{t^-}\bigr]^2=0$ $\mathbb P$-almost surely.
	Thus, if $p(\pi)_{(t,x)}=0$ $\mathbb P$-almost surely then, by the convention $0/0:=0$, $\varphi_{(t,x)}^{(G)}=0$ $\mathbb P$-almost surely.
	Consequently, $\varphi^{(G)}$ is well-defined by \eqref{eq:ExplicitIntegrand}.
	Finally, we note that it is clearly predictable.
\end{remark}

\begin{remark}
	As mentioned in Remark~\ref{rem:AssumptionsAreSatisfiedIf}, one can exhibit a more explicit condition (i.e. the local stability of $N$) which guarantees that \eqref{eq:squareint} and \eqref{eq:IntegrabilityOfPi} are satisfied.
\end{remark}

\section{Proofs of the main results}\label{sec:main}

In this section we prove Theorems~\ref{prop:poincare}, \ref{prop:transport}, \ref{thm:transportlawpp}, \ref{thm:laplace},~\ref{thm:clarkoperative},
and Corollaries \ref{cor:PoincareRenewal}, \ref{cor:Hawkes}, \ref{cor:Cox}, \ref{cor:transportRenewal}, \ref{cor:transportHawkes} and \ref{cor:transportCox}.

\subsection{Proof of Theorem \ref{prop:poincare}}

By Lemma~\ref{le:PapStochlink}, Theorem~\ref{thm:clarkoperative} and Proposition~\ref{prop:ExtendedIsometry}, we have
\begin{align}
	\Var(G)
	&=\mathbb E\bigl[(G-\mathbb E[G])^2\bigr]\nonumber\\
	&=\mathbb E\left[\left(\int_{\R_+\times E}\varphi_{(t,x)}^{(G)}(N(\d t\times\d x)-p(\pi)_{(t,x)}\d t\nu(\d x))\right)^2\right]\nonumber\\
	&=\mathbb E\left[\int_{\R_+\times E}(\varphi_{(t,x)}^{(G)})^2p(\pi)_{(t,x)}\,\d t\nu(\d x)\right]\nonumber\\
	&=\mathbb E\left[\int_{\R_+\times E}\frac{\left(p(\pi G^+)_{(t,x)}-p(G)_{(t,x)}p(\pi)_{(t,x)}\right)^2}{p(\pi)_{(t,x)}}
	\ind_{\{p(\pi)_{(t,x)}>0\}}\,\mathrm{d}t\nu(\d x)\right]\nonumber\\
	&=\mathbb E\left[\int_{\R_+\times E}\frac{\biggl(\mathbb E\left[G^+_{(t,x)}\pi_{(t,x)}-G\,\mathbb E\left[\pi_{(t,x)}\big| \mathcal F_{t^-}\right]\big| \mathcal F_{t^-}\right]\biggr)^2}
	{\mathbb E\left[\pi_{(t,x)}\big| \mathcal F_{t^-}\right]}\ind_{\left\{\mathbb E\left[\pi_{(t,x)}\mid\mathcal F_{t^-}\right]>0\right\}}\,\mathrm{d}t\nu(\d x)\right]\label{eq:VarianceRepresentation}.
\end{align}
Note that
\begin{align}
&\mathbb E\left[G^+_{(t,x)}\pi_{(t,x)}-G\,\mathbb E\left[\pi_{(t,x)}\big| \mathcal F_{t^-}\right]\big| \mathcal F_{t^-}\right]\nonumber\\
&\,\,\,\,\,\,
\,\,\,\,\,\,
=\mathbb E\left[\pi_{(t,x)}\left(G^+_{(t,x)}-G\right)\big| \mathcal F_{t^-}\right]+
\mathbb E\left[(G-\mathbb E[G])\left(\pi_{(t,x)}-\mathbb E\left[\pi_{(t,x)}\big| \mathcal F_{t^-}\right]\right)\big| \mathcal F_{t^-}\right].\nonumber
\end{align}
Thus, for any $0<q<1-\gamma$, by \eqref{eq:VarianceRepresentation} and the convexity inequality
\begin{equation*}
(a+b)^2\le a^2/q+b^2/(1-q),\quad a,b\in\R,
\end{equation*}
we have
\begin{multline}
\label{eq:FirstBoundOnVarF}
\Var(G)
\le\frac1q\,\mathbb E\Biggl[\int_{\R_+\times E}\frac{\left(\mathbb E\left[\pi_{(t,x)}\left(G^+_{(t,x)}-G\right)\big| \mathcal F_{t^-}\right]\right)^2}
{\mathbb E\left[\pi_{(t,x)}\big| \mathcal F_{t^-}\right]}
\,\mathrm{d}t\nu(\d x)\Biggr]\\
+\frac1{1-q}\,\mathbb E\left[\int_{\R_+\times E}\frac{\left(\mathbb E\left[(G-\mathbb E[G])\left(\pi_{(t,x)}-\mathbb E\left[\pi_{(t,x)}\big| \mathcal F_{t^-}\right]\right)
	\big| \mathcal F_{t^-}\right]\right)^2}{\mathbb E\left[\pi_{(t,x)}\big| \mathcal F_{t^-}\right]}
\,\mathrm{d}t\nu(\d x)\right].
\end{multline}
By the Cauchy-Schwarz inequality
\begin{equation*}
\left(\mathbb E\left[\pi_{(t,x)}\Bigl(G^+_{(t,x)}-G\Bigr)\big| \mathcal F_{t^-}\right]\right)^2
\le\mathbb E\left[\pi_{(t,x)}\big| \mathcal F_{t^-}\right]\mathbb E\left[\pi_{(t,x)}\Bigl(G^+_{(t,x)}-G\Bigr)^2\big| \mathcal F_{t^-}\right]
\end{equation*}
and so
\begin{multline}
\label{eq:BoundOnFirstTerm}
\mathbb E\Biggl[\int_{\R_+\times E}\frac{\left(\mathbb E\Bigl[\pi_{(t,x)}\bigl(G^+_{(t,x)}-G\bigr)\big| \mathcal F_{t^-}\Bigr]\right)^2}{\mathbb E\left[\pi_{(t,x)}\big| \mathcal F_{t^-}\right]}
\,\mathrm{d}t\nu(\d x)\Biggr]
\le\mathbb E\left[\int_{\R_+\times E}\left|D_{(t,x)}G\right|^2\pi_{(t,x)}\,\mathrm dt\nu(\d x)\right].
\end{multline}
By the Cauchy-Schwarz inequality
\begin{align*}
\left(\mathbb E\left[(G-\mathbb E[G])\left(\pi_{(t,x)}-\mathbb E\left[\pi_{(t,x)}\big| \mathcal F_{t^-}\right]\right)\big| \mathcal F_{t^-}\right]\right)^2
&\leq\mathbb E\left[(G-\mathbb E[G])^2\left(\mathbb E\left[\pi_{(t,x)}^2\big| \mathcal F_{t^-}\right]-\mathbb E\left[\pi_{(t,x)}\big| \mathcal F_{t^-}\right]^2\right)\big| \mathcal{F}_{t^{-}}\right].
\end{align*}
So, by assumption \eqref{eq:AssumptionPoincare}
and the inequality
\[
(\mathbb E\bigl[\pi_{(t,x)}^2\big| \mathcal F_{t^-}\bigr]/\mathbb E\left[\pi_{(t,x)}\big| \mathcal F_{t^-}\right])-\mathbb E\left[\pi_{(t,x)}\big| \mathcal F_{t^-}\right]\ge 0
\]
which follows from Jensen's inequality, we have
\begin{align}
\mathbb E&\left[\int_{\R_+\times E}\frac{\left(\mathbb E\left[(G-\mathbb E[G])\left(\pi_{(t,x)}-\mathbb E\left[\pi_{(t,x)}\big| \mathcal F_{t^-}\right]\right)\big| \mathcal F_{t^-}\right]\right)^2}
{\mathbb E\left[\pi_{(t,x)}\big| \mathcal F_{t^-}\right]}
\,\mathrm{d}t\nu(\d x)\right]\nonumber\\
&\le\mathbb E\Biggl[(G-\mathbb E[G])^2\int_{\R_+\times E}\Biggl(\frac{\mathbb E\left[\pi_{(t,x)}^2\big| \mathcal F_{t^-}\right]}{\mathbb E\left[\pi_{(t,x)}\big| \mathcal F_{t^-}\right]}
-\mathbb E\left[\pi_{(t,x)}\big| \mathcal F_{t^-}\right]\Biggr)\,\mathrm{d}t\nu(\d x)\Biggr]
\le\gamma\,\Var(G)\label{eq:BoundOnSecondTerm}.
\end{align}
Combining \eqref{eq:FirstBoundOnVarF} with the bounds \eqref{eq:BoundOnFirstTerm} and \eqref{eq:BoundOnSecondTerm}, we obtain
\[
\Var(G)\le q^{-1}\mathbb E\left[\int_{\R_+\times E}\left|D_{(t,x)}G\right|^2\pi_{(t,x)}\,\mathrm dt\nu(\d x)\right]+\gamma(1-q)^{-1}\Var(G),
\]
i.e.
\begin{align}\label{eq:GeneralBoundOnVarG}
\Var(G)\le\frac{1/q}{1-\gamma/(1-q)}\,\mathbb E\left[\int_{\R_+\times E}\left|D_{(t,x)} G\right|^2\pi_{(t,x)}\,\mathrm{d}t\nu(\d x)\right],\quad\text{for any $0<q<1-\gamma$.}
\end{align}
Finally, we note that the choice $q=1-\sqrt\gamma$ minimizes the multiplicative constant appearing in \eqref{eq:GeneralBoundOnVarG},
and the proof is complete.

\subsection{Proofs of Corollaries \ref{cor:PoincareRenewal}, \ref{cor:Hawkes} and \ref{cor:Cox}}

	\noindent{\it Proof\,\,of\,\,Corollary\,\,\ref{cor:PoincareRenewal}.}
	By Corollary~2.9 in \cite{COFormula} $N$ has Papangelou conditional intensity $\{\pi_t\}_{t\in [0,T]}$ defined in the statement.
	By Proposition~2.10 in \cite{COFormula} $N$ is locally stable, and so by Remark~\ref{rem:AssumptionsAreSatisfiedIf} the corresponding assumptions \eqref{eq:squareint} and \eqref{eq:IntegrabilityOfPi} are satisfied.
	In order to verify the corresponding assumption \eqref{eq:AssumptionPoincare}, we compute the quantity
	\begin{equation}\label{eq:var}
		\int_{0}^{T}\frac{\Var[\pi_t\mid\mathcal{F}_{t^-}]}{\mathbb{E}\bigl[\pi_t\mid\mathcal{F}_{t^-}\bigr]}\,\mathrm{d}t=
		\int_{0}^{T}\biggl(\frac{\mathbb E\bigl[\pi_{t}^2\mid\mathcal F_{t^-}\bigr]}
		{\mathbb E\bigl[\pi_{t}\mid\mathcal F_{t^-}\bigr]}-\mathbb E\bigl[\pi_{t}\mid\mathcal F_{t^-}\bigr]\biggr)\mathrm{d}t.
	\end{equation}
	We recall that the points of $N$ are denoted by $T_0:=0<T_1<T_2<\cdots<T_{N([0,T])}$.
	For any $t\in[0,T]$, we have
	\begin{align}
	\mathbb E\bigl[\pi_t^2\mid\mathcal F_{t^-}\bigr]
	&=\mathbb E\biggl[\frac{f(t-T_{N([0,t))})^2f(T_{N([0,t))+1}-t)^2}{f(T_{N([0,t))+1}-T_{N([0,t))})^2}\ind_{\{N((t,T])>0\}}\nonumber\\
	&\qquad\qquad\qquad+\frac{f(t-T_{N([0,t))})^2\overline F(T-t)^2}
	{\overline F(T-T_{N([0,t))})^2}\ind_{\{N((t,T])=0\}}\;\Big|\;\mathcal F_{t^-}\biggr]\nonumber\\
	&=f(t-T_{N([0,t))})^2\mathbb E\biggl[\frac{f\bigl(T_{N([0,t))+1}-T_{N([0,t))}-(t-T_{N([0,t))})\bigr)^2}{f(T_{N([0,t))+1}-T_{N([0,t))})^2}\ind_{\{N((t,T])>0\}}\;\Big|\;\mathcal F_{t^-}\biggr]\nonumber\\
	&\qquad\qquad\qquad+\frac{f(t-T_{N([0,t))})^2\overline F(T-t)^2}
	{\overline F(T-T_{N([0,t))})^2}\mathbb P\bigl(N((t,T])=0\mid\mathcal F_{t^-}\bigr)\nonumber\\
	&\le f(t-T_{N([0,t))})^2h(t-T_{N([0,t))})^2\bigl(1-\mathbb P\bigl(N((t,T])=0\mid\mathcal F_{t^-}\bigr)\bigr)\nonumber\\
	&\qquad\qquad\qquad+\frac{f(t-T_{N([0,t))})^2\overline F(T-t)^2}
	{\overline F(T-T_{N([0,t))})^2}\mathbb P\bigl(N((t,T])=0\mid\mathcal F_{t^-}\bigr).\label{eq:PiTSquare}
	\end{align}
	Hence by \eqref{eq:var}, Lemma~\ref{le:PapStochlink} and Lemma~1 in \cite{bremaudmassouliestab},
	\begin{align}
	&\int_{0}^{T}\frac{\Var[\pi_t\mid\mathcal{F}_{t^-}]}{\mathbb{E}\bigl[\pi_t\mid\mathcal{F}_{t^-}\bigr]}\,\mathrm{d}t\nonumber\\
	&\qquad=\int_{0}^{T}\frac{f(t-T_{N([0,t))})}
	{\overline F(t-T_{N([0,t))})}\biggl[\frac{\overline F(t-T_{N([0,t))})^2\mathbb E\bigl[\pi_t^2\mid\mathcal F_{t^-}\bigr]}{f(t-T_{N([0,t))})^2}-1\biggr]\,\mathrm{d}t\nonumber\\
	&\qquad\le\int_{0}^{T}\frac{f(t-T_{N([0,t))})}
	{\overline F(t-T_{N([0,t))})}\biggl[\overline F(t-T_{N([0,t))})^2h(t-T_{N([0,t))})^2\bigl(1-\mathbb P\bigl(N((t,T])=0\mid\mathcal F_{t^-}\bigr)\bigr)\nonumber\\
	&\qquad\qquad\qquad\qquad\qquad\qquad\qquad+\frac{\overline F(t-T_{N([0,t))})^2\overline F(T-t)^2}
	{\overline F(T-t+(t-T_{N([0,t))}))^2}\mathbb P\bigl(N((t,T])=0\mid\mathcal F_{t^-}\bigr)-1\biggr]\,\mathrm{d}t\nonumber\\
	&\qquad=\int_{0}^{T}\frac{f(t-T_{N([0,t))})}
	{\overline F(t-T_{N([0,t))})}\biggl[\frac{\overline F(t-T_{N([0,t))})^2\overline F(T-t)^2}
	{\overline F(T-t+(t-T_{N([0,t))}))^2}\exp\biggl(-\int_t^T\frac{f(s-t+(t-T_{N([0,t))}))}{\overline F(s-t+(t-T_{N([0,t))}))}\,\mathrm ds\biggr)\nonumber\\
	&\qquad\qquad+\overline F(t-T_{N([0,t))})^2h(t-T_{N([0,t))})^2\biggl(1-\exp\biggl(-\int_t^T\frac{f(s-t+(t-T_{N([0,t))}))}{\overline F(s-t+(t-T_{N([0,t))}))}\,\mathrm ds\biggr)\biggr)-1\biggr]\,\mathrm{d}t\nonumber\\
	&\qquad\le\int_{0}^{T}\sup_{z\in[0,t]}\frac{f(z)}
	{\overline F(z)}\biggl[\frac{\overline F(z)^2\overline F(T-t)^2}
	{\overline F(T-t+z)^2}\exp\biggl(-\int_t^T\frac{f(s-t+z)}{\overline F(s-t+z)}\,\mathrm ds\biggr)\nonumber\\
	&\qquad\qquad\qquad\qquad\qquad\qquad\qquad\qquad+\overline F(z)^2h(z)^2\biggl(1-\exp\biggl(-\int_t^T\frac{f(s-t+z)}{\overline F(s-t+z)}\,\mathrm ds\biggr)\biggr)-1\biggr]\,\mathrm{d}t,\nonumber
	\end{align}
	$\mathbb P$-almost surely.
	It remains to notice that since $f(u)=-\mathrm d\overline F(u)\,/\,\mathrm du$, $u\ge0$, we have
	\begin{equation*}
		\exp\biggl(-\int_t^T\frac{f(s-t+z)}{\overline F(s-t+z)}\,\mathrm ds\biggr)
		=\frac{\overline F(T-t+z)}{\overline F(z)},
	\end{equation*}
	and so, by \eqref{eq:NewBound}, the corresponding assumption \eqref{eq:AssumptionPoincare} holds and the claim follows by Theorem~\ref{prop:poincare}.
	\\
	\noindent$\square$

\noindent{\it Proof\,\,of\,\,Corollary\,\,\ref{cor:Hawkes}.}
By Lemma~2.7 in \cite{COFormula}, $N$ has a Papangelou conditional intensity $\pi_t$ defined by \eqref{eq:PapHaw}.
The assumptions on the parameters $(h,\phi)$ guarantee $\lambda_t\leq\phi(0)$ and
\begin{equation}
\label{eq:BoundOnEt}
	\mathrm{E}_t(N)\leq\mathrm{e}^{\|\phi\|_{\mathrm{Lip}}\int_0^{T-t}h(s)\,\mathrm{d}s}\leq\mathrm{e}^{\|\phi\|_{\mathrm{Lip}}\|h\|_{L^1([0,T],\mathcal{B}([0,T]),\mathrm{d}t)}}.
\end{equation}
Consequently, the inequalities in \eqref{eq:PapHawbis} hold and by Remark~\ref{rem:AssumptionsAreSatisfiedIf} the corresponding conditions \eqref{eq:squareint} and \eqref{eq:IntegrabilityOfPi} are satisfied.
A straightforward computation gives
\begin{align*}
	\int_0^T\frac{\Var\bigl[\pi_{t}\big| \mathcal F_{t^-}\bigr]}{\mathbb E\bigl[\pi_{t}\big| \mathcal F_{t^-}\bigr]}\,\mathrm{d}t&=
	\int_0^T\lambda_t\left(\mathbb{E}\bigl[(\mathrm{E}_t(N))^2\big| \mathcal{F}_{t^-}\bigr]-1\right)\,\mathrm{d}t\\
	&\leq\int_0^T\lambda_t\left(\mathrm{e}^{2\|\phi\|_{\mathrm{Lip}}\int_0^{T-t}h(z)\,\mathrm{d}z}-1\right)\,\mathrm{d}t\\
	&\leq\gamma:=\phi(0)\int_0^T\left(\mathrm{e}^{2\|\phi\|_{\mathrm{Lip}}\int_0^{\tau}h(z)\,\mathrm{d}z}-1\right)\,\mathrm{d}\tau<1,
\end{align*}
and therefore the corresponding assumption \eqref{eq:AssumptionPoincare} holds. The claim follows by Theorem \ref{prop:poincare}.
\\
\noindent$\square$

\noindent{\it Proof\,\,of\,\,Corollary\,\,\ref{cor:Cox}.}
It is well-known (see e.g. \cite{lieshout}) that $N$ has Papangelou conditional intensity $\pi_{(t,x)}$
defined by \eqref{eq:PapCox}. The inequality $\pi_{(t,x)}(\omega)\leq\beta(t,x)$ easily follows by $\lambda_{(t,x)}\leq\beta(t,x)$ and so $N$
is locally stable. Since $\lambda_{(t,x)}\geq\alpha(t,x)$, condition~\eqref{eq:intless1} is implied by \eqref{eq:intless2}.
So, by Remark~\ref{rem:locallybis}, Theorem~\ref{prop:poincare} can be applied and the claim follows.
\\
\noindent$\square$

\subsection{Proof of Theorem \ref{prop:transport}}

The proof of Theorem \ref{prop:transport} is based on two preliminary propositions.
The first one consists in a result from \cite{gozlan} (see Theorem~3 therein, as well as Theorem~1.1 in \cite{MaWu}).
The second one,
whose proof is given at the end of this subsection,
provides a stochastic convex inequality for functionals of marked point processes and generalizes Theorem~4.1-$(ii)$ in \cite{kleinmaprivault}.
We recall that in general, $\chi$ denotes a Polish space equipped with its Borel $\sigma$-field $\mathcal{B}(\chi)$ and
$d$ denotes a lower semi-continuous metric on $\chi$
which does not necessarily generate the topology on $\chi$.

\begin{Proposition}\label{le:gozlan}
	Let $c:[0,\infty)\to [0,\infty]$ be a non-decreasing, left-continuous and convex function with $c(0)=0$, and let $\mu\in M_1(\chi,d)$.
	Then
	\[
	c(\mathcal T_d(\sigma,\mu))\leq H(\sigma\mid\mu),\quad\text{for any $\sigma\in M_1(\chi,d)$}
	\]
	if and only if, for any function $f:\chi\to\R$ which is measurable, bounded, Lipschitz continuous (with respect to the metric $d$) with Lipschitz constant less than or equal to $1$, i.e. $\sup_{x\neq y}\frac{|f(x)-f(y)|}{d(x,y)}\leq 1$, we have
	\[
		\int_\chi\mathrm{e}^{\theta\left(f(x)-\int_\chi f(z)\mu(\d z)\right)}\mu(\d x)\le\mathrm{e}^{c^\odot(\theta)},\quad\text{for any $\theta\in\R_+$.}
	\]
\end{Proposition}

\begin{Proposition}\label{lem:Theorem41Privault}
Let the notation of Theorem \ref{thm:clarkoperative} prevail. Assume that $N$ satisfies \eqref{eq:LocalStability} with a dominating function $\beta$, $G\in L^2(\Omega,\mathcal{F}_\infty,\mathbb P)$
and $\bigl|\varphi_{(t,x)}^{(G)}\bigr|\le h(t,x)$, $\d t\mathbb P(\d\omega)\nu(\d x)$-almost everywhere, for some deterministic function $h$ such that
\[
	\int_{\R_+\times E}|h(t,x)|^2\beta(t,x)\,\d t\nu (\d x)<\infty.
\]
Then, letting $\mathbb E'$ denote the expectation corresponding to a probability measure $\mathbb P'$
on $(\Omega,\mathcal F_\infty)$ under which $N$ is a Poisson process on $\R_+\times E$ with intensity function $\beta$,
\begin{equation}\label{g2.2}
\mathbb E\bigl[\phi(G-\mathbb E[G])\bigr]\le{\mathbb E'}\biggl[\phi\biggl(\int_{\R_+\times E}h(t,x)\bigl(N(\d t\times\d x)-\beta(t,x)\,\d t\nu(\d x)\bigr)\biggr)\biggr],
\end{equation}
for all twice continuously differentiable convex functions $\phi:\mathbb R\to\mathbb R$ such that $\phi'$ is convex.
\end{Proposition}

\noindent{\it Proof\,\,of\,\,Theorem\,\,\ref{prop:transport}.} We take $\chi=\mathbb R$ and let
$f:\R\to\R$ be bounded and Lipschitz continuous with Lipschitz constant less than or equal to $1$ (i.e $\sup_{x\neq y}|f(x)-f(y)|/d(x,y)\leq 1$).
We note that for $\mathrm dt\nu(\d x)\mathbb P(\d\omega)$-almost every $(t,x,\omega)$,
\begin{multline*}
	p(\pi f\circ G^+)_{(t,x)}-p(f\circ G)_{(t,x)}p(\pi)_{(t,x)}\\
	=\mathbb{E}[\pi_{(t,x)}D_{(t,x)}(f\circ G)\mid\mathcal{F}_{t^-}]
	+\mathbb{E}\bigl[\bigl(f\circ G-f(\mathbb{E}[G\mid\mathcal{F}_{t^-}])\bigr)\bigl(\pi_{(t,x)}-\mathbb E\bigl[\pi_{(t,x)}\big| \mathcal F_{t^-}\bigr]\bigr)\mid\mathcal{F}_{t^-}\bigr],
\end{multline*}
and so
\begin{align*}
	\left|\varphi_{(t,x)}^{(f\circ G)}(\omega)\right|
	&\le\frac{\mathbb E\Bigl[\pi_{(t,x)}\bigl|f\circ G^+_{(t,x)}-f\circ G\bigr|+\bigl|f\circ G-f(\mathbb{E}[G\mid\mathcal{F}_{t^-}])\bigr|\bigl|\pi_{(t,x)}-\mathbb E\bigl[\pi_{(t,x)}\mid \mathcal F_{t^-}\bigr]\bigr|\;\big|\;\mathcal F_{t^-}\Bigr](\omega)}{\mathbb E\left[\pi_{(t,x)}\mid\mathcal F_{t^-}\right](\omega)}\\
	&\le\frac{\mathbb E\left[\pi_{(t,x)}\|D_{(t,x)} G\|_d+\|G-\mathbb{E}[G\mid\mathcal{F}_{t^-}]\|_d\left|\pi_{(t,x)}-\mathbb E\left[\pi_{(t,x)}\mid\mathcal F_{t^-}\right]\right|\;\big|\;\mathcal F_{t^-}\right](\omega)}
	{\mathbb E\left[\pi_{(t,x)}\mid\mathcal F_{t^-}\right](\omega)}\\
	&\le h(t,x).
\end{align*}
So, for any $\theta\geq 0$, by Proposition~\ref{lem:Theorem41Privault} with $\phi(x):=\mathrm{e}^{\theta x}$ and $f(G)$ in place of $G$, we have
\begin{align*}
\mathbb E\Bigl[\mathrm{e}^{\theta(f(G)-\mathbb E\bigl[f(G)\bigr])}\Bigr]
&\le{\mathbb E'}\biggl[\exp\biggl(\theta\int_{\R_+\times E}h(t,x)\bigl(N(\d t\times\d x)-\beta(t,x)\,\d t\nu(\d x)\bigr)\biggr)\biggr]\\
&=\mathrm{e}^{\Lambda(\theta)},
\end{align*}
where $\Lambda$ is defined by \eqref{eq:Lambda}.
A straightforward computation shows that $\Lambda$ is non-negative, non-decreasing, left-continuous and convex, with $\Lambda(0)=0$.
Thus, by Proposition~1 in \cite{gozlan}
$\Lambda^{\odot\odot}=\Lambda$. The transportation cost inequality \eqref{eq:traspineq} then follows by Proposition \ref{le:gozlan}.
Now, assuming that $h$ is bounded by $M>0$, we have
\begin{equation*}
	\Lambda(\theta)
	\le\frac{\mathrm e^{\theta M}-\theta M-1}{M^2}\int_{\R_+\times E}h(t,z)^2\beta(t,z)\,\d t\nu(\d z),
\end{equation*}
and as in the proof of Theorem~2.6 in \cite{MaWu}, we get the inequality \eqref{eq:LowerBoundOfC1}.
\\
\noindent$\square$

\noindent{\it Proof\,\,of\,\,Proposition\,\,\ref{lem:Theorem41Privault}.}
This proof is inspired by that of Lemma~3.2 of \cite{wu1}.
Throughout this proof all the random quantities are defined on the product probability space $(\Omega^2,\mathcal{F}_\infty\otimes\mathcal{F}_\infty,\mathbb P\otimes\mathbb P')$,
and we let $N'(\omega,\omega')=\omega'$ and $N(\omega,\omega')=\omega$.
With an abuse of notation, we set $p(\pi)_{(t,x)}(\omega,\omega'):=p(\pi)_{(t,x)}(\omega)$, $\varphi^{(G)}_{(t,x)}(\omega,\omega'):=\varphi^{(G)}_{(t,x)}(\omega)$,
and we denote by $\widetilde{\mathbb E}$ the expectation with respect to $\mathbb P\otimes\mathbb P'$.

Let $\{M_t\}_{t\in\R_+}$ and $\{M_t^*\}_{t\in\R_+}$ be the stochastic processes defined, respectively, by
$$
	M_t(\omega,\omega'):=\int_{[0,t]\times E}\varphi_{(s,x)}^{(G)}(\omega)(\omega(\d s\times\d x )-p(\pi)_{(s,x)}(\omega)\,\d s\nu(\d x))
$$
and
\[	
	M_t^*(\omega,\omega'):=\int_{[t,\infty)\times E}h(s,x)(\omega'(\d s\times\d x )-\beta(s,x)\,\d s\nu(\d x)).
\]
Let $\{\mathcal{H}_t\}_{t\geq 0}$ and $\{\mathcal{H}_t^*\}_{t\geq 0}$ be, respectively, the forward and backward filtrations defined by
\[
	\mathcal{H}_t:=\mathcal F_t\otimes\mathcal{F}_\infty\quad\text{and}\quad
	\mathcal{H}_t^*:=\mathcal F_\infty\otimes\mathcal F_{[t,\infty)},\quad t\in\R_+.
\]
By Corollary~C4 p.~235 in \cite{bremaud} and standard properties of the conditional expectation, we have that $\{M_t\}_{t\geq 0}$ is an $\mathcal{H}^*$-adapted $\mathcal{H}$-martingale and $\{M_t^*\}_{t\geq 0}$ is an $\mathcal{H}$-adapted $\mathcal{H}^*$-backward martingale.
Letting $\varepsilon_{(t,u)}$ and $\varepsilon_t$ denote, respectively, the Dirac measure at $(t,u)\in\R_+\times\R$ and at $t\in\R_+$, we define the jump measures of $\{M_t\}_{t\in\R_+}$ and $\{M_t^*\}_{t\in\R_+}$ respectively by
\begin{equation*}
	\mu(\d s\times\d\tau)
	:=\sum_{t>0\st\Delta M_t\neq0}\varepsilon_{(t,\Delta M_t)}(\d s\times\d\tau)
	=\sum_{(t,x)\in\mathrm{Supp}(N)}\varepsilon_{\bigl(t,\varphi_{(t,x)}^{(G)}\bigr)}(\d s\times\d\tau)
\end{equation*}
and
\begin{equation*}
	\mu^*(\d s\times\d\tau)
	:=\sum_{(t,x)\in\mathrm{Supp}(N')}\varepsilon_{(t,h(t,x))}(\d s\times\d\tau),
\end{equation*}
where $\Delta M_t:=M_t-M_{t^-}$.
For any fixed $t\in\R_+$, denote by $\nu_t (\d\tau)$ the (random) image measure on $(\R,\mathcal{B}(\R))$ of $p(\pi)_{(t,x)}\nu(\d x)$ by the mapping $E\ni x\mapsto\varphi_{(t,x)}^{(G)}\in\R$, i.e. for any bounded and measurable $f:\R\to\R$,
\[
	\int_{\R}f(\tau)\,\nu_t(\d\tau):=\int_E f\Bigl(\varphi_{(t,x)}^{(G)}(\omega)\Bigr)p(\pi)_{(t,x)}(\omega)\,\nu(\d x),
\]
and similarly let $\nu_t^*$ be the measure on $(\R,\mathcal{B}(\R))$ defined by
\[
	\int_{\R}f(\tau)\,\nu_t^*(\d\tau):=\int_E f\bigl(h(t,x)\bigr)\beta(t,x)\,\nu(\d x).
\]
It turns out that $\nu_t(\d\tau)\,\d t$ is the dual $\mathcal H$-predictable projection of $\mu$ and $\nu_t^*(\d\tau)\,\d t$ is the dual $\mathcal H^*$-predictable projection of $\mu^*$.
Indeed focusing e.g. on $\mu$, again by Corollary~C4 p.~235 in \cite{bremaud} and standard properties of the conditional expectation, for any $t,\Delta t\ge0$ and $A\in\mathcal B(\R)$,
\begin{multline*}
	\widetilde{\mathbb E}\biggl[\mu([t,t+\Delta t]\times A)-\int_t^{t+\Delta t}\nu_s(A)\,\d s\big| \mathcal H_t\biggr]\\
	=\widetilde{\mathbb E}
	\biggl[\int_{[t,t+\Delta t]\times E}\ind_{\bigl\{\varphi_{(s,x)}^{(G)}\in A\bigr\}}(N(\d s\times\d x)-p(\pi)_{(s,x)}\d s\nu(\d x))\big| \mathcal H_t\biggr]
	=0.
\end{multline*}
Consequently, conditions $(3.1)$, $(3.2)$ and $(3.3)$ of \cite{kleinmaprivault} are verified.
We also note that condition $(3.4)$ of \cite{kleinmaprivault} is trivially satisfied with $H\equiv H^*\equiv 0$.
For any $t\in\R_+$, we define the following (random) measures on $(\R,\mathcal{B}(\R))$:
$$
	\widetilde{\nu}_t (\d\tau) := |\tau|^2 \nu_t (\d\tau)
	\quad\text{and}\quad
	\widetilde{\nu}^*_t (\d\tau) := |\tau|^2 \nu^*_t (\d\tau).
$$
For any $u\in\mathbb R$, we have
\begin{align*}
	\widetilde{\nu}_t ([u,\infty ))
	=\int_E\ind_{[u,\infty )}\bigl(\varphi_{(t,x)}^{(G)}\bigr)\bigl|\varphi_{(t,x)}^{(G)}\bigr|^2p(\pi)_{(t,x)}\,\nu (\d x)
	&\leq\int_E\ind_{[u,\infty )}\bigl(h(t,x)\bigr)\bigl|h(t,x)\bigr|^2\beta(t,x)\,\nu (\d x)\\
	&=\widetilde{\nu}^*_t ([u,\infty )).
\end{align*}
Furthermore, for any $u\in\R$,
\begin{equation*}
	\int_{0}^\infty\widetilde{\nu}^*_t ([u,\infty ))\,\d t
	\le\int_{\R_+\times E}\bigl|h(t,x)\bigr|^2\beta(t,x)\,\d t\nu (\d x)
	<\infty.
\end{equation*}
Therefore, for any $u\in\R$ and any $(\omega,\omega')\in\Omega^2$, $\widetilde{\nu}^*_t ([u,\infty ))<\infty$ $\d t$-almost everywhere, and so by Theorem~3.3 in \cite{kleinmaprivault} we have
\[
	\widetilde{\mathbb E}[\phi(M_t+M_t^*)]\le\widetilde{\mathbb E}[\phi(M^*_0)],
\]
for all $t\in\R_+$ and any function $\phi$ as in the statement.
By Theorem~\ref{thm:clarkoperative} and Proposition~\ref{bremaudmassoulie} we have
\begin{align*}
	\widetilde{\mathbb E}[(M_t+M_t^*+\mathbb E[G]-G)^2]
	&\le2\biggl(\widetilde{\mathbb E}\biggl[\biggl(\int_{[t,\infty)\times E}h(s,x)\bigl(N'(\d s\times\d x )-\beta(s,x)\,\d s\nu(\d x)\bigr)\biggr)^2\biggr]\\
	&\qquad\qquad+\widetilde{\mathbb E}\biggl[\biggl(\int_{[t,\infty)\times E}\varphi_{(s,x)}^{(G)}\bigl(N(\d s\times\d x )-p(\pi)_{(s,x)}\,\d s\nu(\d x)\bigr)\biggr)^2\biggr]\biggr)\\
	&=2\int_{[t,\infty)\times E}|h(s,x)|^2\beta(s, x)\,\d s\nu(\d x)\\
	&\qquad\qquad+2\mathbb E\biggl[\int_{[t,\infty)\times E}\bigl|\varphi_{(s,x)}^{(G)}\bigr|^2p(\pi)_{(s,x)}\,\d s\nu(\d x)\biggr]\\
	&\xrightarrow[t\to\infty]{}0.
\end{align*}
Thus, there exists a sequence $\{t_n\}_{n\geq 1}$ such that $M_{t_n}+M_{t_n}^*\to G-\mathbb{E}{G}$, $\widetilde{\mathbb{P}}$-almost surely, and therefore by Fatou's lemma
\begin{equation*}
	\widetilde{\mathbb E}[\phi(G-\mathbb E[G])]
	=\widetilde{\mathbb E}\Bigl[\liminf_{n\to\infty}\phi (M_{t_n}+M_{t_n}^*)\Bigr]
	\le\liminf_{n\to\infty}\widetilde{\mathbb E}[\phi (M_{t_n}+M_{t_n}^*)]
	\le\widetilde{\mathbb E}[\phi(M^*_0)],
\end{equation*}
which is exactly \eqref{g2.2}.
\\
\noindent$\square$

\subsection{Proofs of Proposition~\ref{prop:transportFirstOrder} and Corollaries \ref{cor:transportRenewal}, \ref{cor:transportHawkes}, \ref{cor:transportCox}}
\label{subsec:proofOfApplicationsCorollaries}

\noindent{\it Proof\,\,of\,\,Proposition\,\,\ref{prop:transportFirstOrder}.}
We shall apply Theorem~\ref{prop:transport}.
For any $(t,x)\in\R_+\times E$, we have
\begin{equation*}
	\bigl|D_{(t,x)}G\bigr|=|g(t,x)|\leq g_1(t,x),
\end{equation*}
and by the Cauchy-Schwarz inequality
\begin{multline*}
	\mathbb E\bigl[|G-\mathbb E[G\mid\mathcal F_{t^-}]||\pi_{(t,x)}-\mathbb E[\pi_{(t,x)}\mid\mathcal F_{t^-}]|\mid\mathcal F_{t^-}\bigr]\bigr|\\
	\le\mathbb E\bigl[(G-\mathbb E[G\mid\mathcal F_{t^-}])^2\mid\mathcal F_{t^-}\bigr]^{1/2}
	\mathbb E\bigl[(\pi_{(t,x)}-\mathbb E[\pi_{(t,x)}\mid\mathcal F_{t^-}])^2\mid\mathcal F_{t^-}\bigr]^{1/2}.
\end{multline*}
Therefore
\begin{multline*}
	\frac{\mathbb{E}\bigl[\bigl|G-\mathbb E[G\mid\mathcal F_{t^-}]\bigr||\pi_{(t,x)}-p(\pi)_{(t,x)}|\mid \mathcal{F}_{t^-}\bigr](\omega)}{p(\pi)_{(t,x)}(\omega)}\\
	\le\mathbb E\bigl[(G-\mathbb E[G\mid\mathcal F_{t^-}])^2\mid\mathcal F_{t^-}\bigr]^{1/2}\Biggl\|\frac{\sqrt{\mathbb V\mathrm{ar}\bigl[\pi_{(t,x)}\mid\mathcal F_{t^-}\bigr]}}{p(\pi)_{(t,x)}}\Biggr\|_{L^\infty(\Omega,\mathcal F_\infty,\mathbb P)}.
\end{multline*}
Additionally, by \eqref{eq:classicalStochasticIntensity}
\begin{equation*}
	\mathbb E[G\mid\mathcal F_{t^-}]
	=\int_{(0,t)\times E}g(s,y)\,N(\mathrm ds\times\mathrm dy)
	+\int_{[t,\infty)\times E}g(s,y)\mathbb E[\lambda_{(s,y)}\mid\mathcal F_{t^-}]\,\mathrm ds\nu(\mathrm dy),
\end{equation*}
and so by Proposition~\ref{bremaudmassoulie}{\it(ii)}
\begin{align}
	&\mathbb E\bigl[(G-\mathbb E[G\mid\mathcal F_{t^-}])^2\mid\mathcal F_{t^-}\bigr]\nonumber\\
	&\qquad=\mathbb E\biggl[\biggl(\int_{[t,\infty)\times E}g(s,y)\,N(\mathrm ds\times\mathrm dy)
	-\int_{[t,\infty)\times E}g(s,y)\mathbb E[\lambda_{(t,y)}\mid\mathcal F_{t^-}]\,\mathrm ds\nu(\mathrm dy)\biggr)^2\;\Big|\;\mathcal F_{t^-}\biggr]\nonumber\\
	&\qquad\le2\mathbb E\biggl[\biggl(\int_{[t,\infty)\times E}g(s,y)\,(N(\mathrm ds\times\mathrm dy)-\lambda_{(s,y)}\mathrm ds\nu(\mathrm dy)\biggr)^2\nonumber\\
	&\qquad\qquad\qquad\qquad\qquad\qquad\qquad+\biggl(\int_{[t,\infty)\times E}g(s,y)(\lambda_{(s,y)}-\mathbb E[\lambda_{(s,y)}\mid\mathcal F_{t^-}])\,\mathrm ds\nu(\mathrm dy)\biggr)^2\;\Big|\;\mathcal F_{t^-}\biggr]\nonumber\\
	&\qquad=2\biggl[\int_{[t,\infty)\times E}g(s,y)^2\mathbb E[\lambda_{(s,y)}\mid\mathcal F_{t^-}]\,\mathrm ds\nu(\mathrm dy)\nonumber\\
	&\qquad\qquad\qquad\qquad\qquad\qquad\qquad+\mathbb E\biggl[\biggl(\int_{[t,\infty)\times E}g(s,y)(\lambda_{(s,y)}-\mathbb E[\lambda_{(s,y)}\mid\mathcal F_{t^-}])\,\mathrm ds\nu(\mathrm dy)\biggr)^2\;\Big|\;\mathcal F_{t^-}\biggr]\biggr]\nonumber\\
	&\qquad\le2\biggl[\int_{[t,\infty)\times E}g(s,y)^2\beta(s,y)\,\mathrm ds\nu(\mathrm dy)
	+\biggl(\int_{[t,\infty)\times E}|g(s,y)|\beta(s,y)\,\mathrm ds\nu(\mathrm dy)\biggr)^2\biggr].\label{eq:BoundOnConditionalVarianceG}
\end{align}
The claim follows by Theorem~\ref{prop:transport}.
\\
\noindent$\square$

\noindent{\it Proof\,\,of\,\,Corollary\,\,\ref{cor:transportRenewal}.}
By the proof of Proposition 2.10 in \cite{COFormula} we have $\pi_t\leq\beta$, where $\beta$ is given by \eqref{eq:betaren}.
Additionally, by \eqref{eq:AssumptionOnHBar}, \eqref{eq:PiTSquare} and the form of the stochastic intensity for renewal processes,
\begin{align*}
	\frac{\sqrt{\mathbb V\mathrm{ar}\bigl[\pi_{t}\mid\mathcal F_{t^-}\bigr]}}{\lambda_t}
	=\sqrt{\mathbb E\biggl[\frac{\pi_t^2}{\lambda_t^2}\;\Big|\;\mathcal F_{t^-}\biggr]-1}
	&\le\sqrt{\overline h^2\overline F(t-T_{N([0,t))})^2+\frac{\overline F(T-t)^2\overline F(t-T_{N([0,t))})^2}{\overline F(T-T_{N([0,t))})^2}-1}\\
	&\le\sqrt{\overline h^2+\biggl(\int_C^Tf(x)\,\mathrm dx\biggr)^{-2}-1}.
\end{align*}
The claim follows by Proposition~\ref{prop:transportFirstOrder}.
\\
\noindent$\square$

\noindent{\it Proof\,\,of\,\,Corollary\,\,\ref{cor:transportHawkes}.}
By the first inequality in \eqref{eq:PapHawbis} we can take as dominating function of $\pi_t$
\begin{equation*}
	\beta(t)=\phi(0)\exp\biggl(\|\phi\|_{\mathrm{Lip}}\int_0^{T-t}h(z)\,\mathrm dz\biggr).
\end{equation*}
Additionally, by the first inequality in \eqref{eq:BoundOnEt},
\begin{equation*}
	\frac{\sqrt{\mathbb V\mathrm{ar}\bigl[\pi_{t}\mid\mathcal F_{t^-}\bigr]}}{\lambda_t}
	=\sqrt{\mathbb E\biggl[\frac{\pi_t^2}{\lambda_t^2}\;\Big|\;\mathcal F_{t^-}\biggr]-1}
	\le\sqrt{\exp\biggl(2\|\phi\|_{\mathrm{Lip}}\int_0^{T-t}h(z)\,\mathrm dz\biggr)-1}.
\end{equation*}
The claim follows by Proposition~\ref{prop:transportFirstOrder}.
\\
\noindent$\square$

\noindent{\it Proof\,\,of\,\,Corollary\,\,\ref{cor:transportCox}.} We already noticed in the proof of Corollary \ref{cor:Cox} that $\beta$ is a dominating function of
the Papangelou conditional intensity. The claim easily follows by Proposition~\ref{prop:transportFirstOrder} noticing that
\begin{equation*}
	\frac{\sqrt{\mathbb V\mathrm{ar}\bigl[\pi_{(t,x)}\mid\mathcal F_{t^-}\bigr]}}{\lambda_{(t,x)}}
	=\sqrt{\mathbb E\biggl[\frac{\pi_{(t,x)}^2}{\lambda_{(t,x)}^2}\;\Big|\;\mathcal F_{t^-}\biggr]-1}
	\le\sqrt{\frac{\beta(t,x)^2}{\alpha(t,x)^2}-1}.
\end{equation*}
\\
\noindent$\square$

\subsection{Proof of Theorem \ref{thm:transportlawpp}}

Since the proof is conceptually similar to that of Theorem \ref{prop:transport},
we only emphasize the main differences. We take $\chi=\Omega$, $d:=d_\varphi$ and let $F:(\Omega,d_\varphi)\to\R$ be bounded and Lipschitz continuous with Lipschitz constant less than or equal to one.
We have for $\mathrm dt\nu(\d x)\mathbb P(\d\omega)$-almost every $(t,x,\omega)$,
\begin{align*}
	\left|\varphi_{(t,x)}^{(F)}(\omega)\right|
	&\le\frac{\mathbb E\Bigl[\pi_{(t,x)}\bigl|F(N+\varepsilon_{(t,x)})-F(N)\bigr|+\bigl|F(N)-F(N|_{[0,t)\times E})\bigr|\bigl|\pi_{(t,x)}-\mathbb E\bigl[\pi_{(t,x)}\mid \mathcal F_{t^-}\bigr]
	\bigr|\mid\mathcal F_{t^-}\Bigr](\omega)}{\mathbb E\left[\pi_{(t,x)}\mid\mathcal F_{t^-}\right](\omega)}\\
	&\le\frac{\mathbb E\Bigl[\pi_{(t,x)}d_\varphi(N+\varepsilon_{(t,x)},N)+d_\varphi(N,N|_{[0,t)\times E})\bigl|\pi_{(t,x)}-\mathbb E\bigl[\pi_{(t,x)}\mid\mathcal F_{t^-}\bigr]\bigr|\mid \mathcal F_{t^-}\Bigr](\omega)}{\mathbb E\left[\pi_{(t,x)}\mid \mathcal F_{t^-}\right](\omega)}\\
	&=\varphi(t,x)+\frac{\mathbb E\Bigl[\left(\int_{[t,\infty)\times E}\varphi(s,y)N(\d s\times\d y)\right)
	\bigl|\pi_{(t,x)}-\mathbb E\bigl[\pi_{(t,x)}\mid\mathcal F_{t^-}\bigr]\bigr|\mid \mathcal F_{t^-}\Bigr](\omega)}{\mathbb E\left[\pi_{(t,x)}\mid \mathcal F_{t^-}\right](\omega)}\\
	&\le\varphi(t,x)+\Biggl\|\frac{\sqrt{\mathbb V\mathrm{ar}\bigl[\pi_{(t,x)}\mid\mathcal F_{t^-}\bigr]}}{p(\pi)_{(t,x)}}\Biggr\|_{L^\infty(\Omega,\mathcal F_\infty,\mathbb P)}\sqrt{\mathbb E\biggl[\biggl(\int_{[t,\infty)\times E}\varphi(s,y)N(\d s\times\d y)\biggr)^2\;\Big|\;\mathcal F_{t^-}\biggr](\omega)}.
\end{align*}
Proceeding similarly to the series of inequalities \eqref{eq:BoundOnConditionalVarianceG}, we obtain
\begin{multline*}
	\left|\varphi_{(t,x)}^{(F)}(\omega)\right|
	\le\varphi(t,x)+\sqrt 2\biggl[\int_{[t,\infty)\times E}\varphi(s,y)^2\beta(s,y)\,\mathrm ds\nu(\mathrm dy)
	+\biggl(\int_{[t,\infty)\times E}\varphi(s,y)\beta(s,y)\,\mathrm ds\nu(\mathrm dy)\biggr)^2\biggr]^{1/2}\\
	\times\Biggl\|\frac{\sqrt{\mathbb V\mathrm{ar}\bigl[\pi_{(t,x)}\mid\mathcal F_{t^-}\bigr]}}{p(\pi)_{(t,x)}}\Biggr\|_{L^\infty(\Omega,\mathcal F_\infty,\mathbb P)}.
\end{multline*}
So, for any $\theta\geq 0$, by Proposition \ref{lem:Theorem41Privault} with $\phi(x):=\mathrm{e}^{\theta x}$, we have
\begin{equation*}
	\mathbb E\Bigl[\mathrm{e}^{\theta(F-\mathbb E[F])}\Bigr]
	\le{\mathbb E'}\biggl[\exp\biggl(\theta\int_{\R_+\times E}h_\varphi(t,x)\bigl(N(\d t\times\d x)-\beta(t,x)\,\d t\nu(\d x)\bigr)\biggr)\biggr]
	=\mathrm{e}^{\Lambda_\varphi(\theta)}.
\end{equation*}
Note also that
\begin{multline*}
	\int_\Omega d_\varphi(\omega,\mathbf 0)\,\mathbb P(\mathrm d\omega)
	=\mathbb E\bigl[d_\varphi(N,\mathbf0)\bigr]
	=\mathbb E\biggl[\int_{\R_+\times E}\varphi(s,y)\,N(\mathrm ds\times\mathrm dy)\biggr]\\
	\le\int_{\R_+\times E}\varphi(s,y)\beta(s,y)\,\mathrm ds\nu(\mathrm dy)
	<\infty.
\end{multline*}
The remainder of the proof is similar to that of the final part of Theorem~\ref{prop:transport}.

\subsection{Proof of Theorem~\ref{thm:laplace}}\label{subsec:lap}

The proof of Theorem~\ref{thm:laplace} is based on two preliminary propositions, which extend to our setting Propositions~4.2~and~4.3 in \cite{zhang1}, respectively.

\begin{Proposition}\label{prop:1}
	Let the assumptions and notation of Theorem~\ref{thm:clarkoperative} prevail and let $G$ be a random variable on
	$(\Omega,\mathcal{F}_\infty,\mathbb P)$ which satisfies $0<c_0\leq G\leq c_1$ almost surely, for some positive constants $c_0,c_1>0$.
	Setting
	\begin{equation}\label{eq:fi}
		\phi_{(t,x)}^{(G)}:=\frac{\varphi_{(t,x)}^{(G)}}{p(G)_{(t,x)}},\quad t\in\R_+, x\in E,
	\end{equation}
	we have $\phi^{(G)}\in\mathcal{H}$.
	Additionally, for any $t\in\R_+$ we have
	\begin{equation}
	\label{eq:ConditionalExpectationAsDoobExponential}
		\mathbb E\bigl[G\big| \mathcal F_t\bigr]
		=\mathbb E[G]\cdot\mathcal E_t(\phi^{(G)}).
	\end{equation}
\end{Proposition}

\begin{Proposition}\label{prop:2}
	Let the assumptions and notation of Theorem~\ref{thm:clarkoperative} prevail and, for $\phi\in\mathcal H$, let $\mathbb P_\phi$ be defined by \eqref{eq:DefinitionPPhi}.
	For any predictable stochastic process $\psi$ such that
	\begin{equation}
	\label{eq:IntegrabilityAssumptionOnPsi}
		\mathbb E\biggl[\biggl(\int_{[0,T]\times E}|\psi_{(s,x)}|^2\,p(\pi)_{(s,x)}\,\d s\nu(\d x)\biggr)^2\biggr]<\infty,\quad\text{for any $T>0$}
	\end{equation}
	we have that, under $\mathbb P_\phi$,
	\[
		\biggl\{\int_{[0,t]\times E}\psi_{(s,x)}(N(\d s\times\d x)-p(\pi)_{(s,x)}\,\d s\nu(\d x))-\int_{[0,t]\times E}\psi_{(s,x)}\phi_{(s,x)}p(\pi)_{(s,x)}\,\d s\nu(\d x)\biggr\}_{t\in\R_+}
	\]
	is a square-integrable $\mathcal{F}$-martingale with null mean.
\end{Proposition}

These propositions are proved at the end of this subsection, and we start proving Lemma~\ref{le:mgproperty}.\\\\
\noindent{\it Proof\,\,of\,\,Lemma\,\,\ref{le:mgproperty}.} By It\^o's formula (see e.g. Theorem 5.1 p. 66 of \cite{watanabe}) we have
\begin{equation}\label{eq:SDEOfE}
	\mathcal{E}_t(\phi)=1+\int_{[0,t]\times E}\mathcal{E}_{s^-}(\phi)\phi_{(s,x)}(N(\d s\times\d x)-p(\pi)_{(s,x)}\d s\nu(\d x)),\quad t\in\R_+.
\end{equation}
Moreover, for any $T>0$ we have
\begin{equation}\label{eq:p12}
	\{\ind_{[0,T]}(t)\mathcal{E}_t(\phi)\phi_{(t,x)}\}_{(t,x)\in\R_+\times E}\in\mathcal{P}_{1,2}(p(\pi)).
\end{equation}
Indeed, letting $M>0$ denote a constant such that $|\phi_{(t,x)}|\le M$ $\d t\nu(\d x)\d\mathbb P$-almost everywhere, we have
\begin{align}
	\mathbb E\biggl[\int_{\R_+\times E}|\ind_{[0,T]}(t)\mathcal{E}_t(\phi)\phi_{(t,x)}|&p(\pi)_{(t,x)}\,\d t\nu(\d x)\biggr]\nonumber\\
	&\le M\mathbb E\biggl[\int_{[0,T]\times E}\mathcal{E}_t(\phi)p(\pi)_{(t,x)}\,\d t\nu(\d x)\biggr]\nonumber\\
	&\le M
        \sqrt{\mathbb E\biggl[\int_{[0,T]\times E}\mathcal{E}_t(\phi)^2p(\pi)_{(t,x)}\,\d t\nu(\d x)\biggr]}
        \sqrt{
          \mathbb E\bigl[N([0,T]\times E)\bigr]}
        \nonumber\\
	&<\infty\label{eq:Ifin}
\end{align}
and
\begin{align}
	\mathbb E\biggl[\int_{\R_+\times E}|\ind_{[0,T]}(t)\mathcal{E}_t(\phi)\phi_{(t,x)}|^2&p(\pi)_{(t,x)}\,\d t\nu(\d x)\biggr]\nonumber\\
	&\le M^2\mathbb E\biggl[\int_{[0,T]\times E}\mathcal{E}_t(\phi)^2p(\pi)_{(t,x)}\,\d t\nu(\d x)\biggr]\nonumber\\
	&<\infty,\label{eq:IIfin}
\end{align}
where \eqref{eq:Ifin} and \eqref{eq:IIfin} follow by \eqref{eq:Ep2} and the square integrability of $N([0,T]\times E)$.
By \eqref{eq:SDEOfE} and Proposition~\ref{bremaudmassoulie}-$(ii)$ we then have
\begin{equation}\label{eq:SecondMomentOfEt}
	\mathbb E\bigl[\mathcal E_t(\phi)^2\bigr]=1+\mathbb E\biggr[\int_{[0,t]\times E}\phi_{(s,x)}^2\mathcal E_s(\phi)^2p(\pi)_{(s,x)}\,\d s\nu(\d x)\biggr],
\end{equation}
which gives the square integrability of $\{\mathcal{E}_t(\phi)\}_{t\in\R_+}$.
The martingale property follows by \eqref{eq:SDEOfE} and Corollary~C4 p.~235 of \cite{bremaud}.\\
\noindent$\square$

\noindent{\it Proof\,\,of\,\,Theorem\,\,\ref{thm:laplace}.}
We divide the proof in three steps.
In the first step we identify $\d\mathbb P_{\phi}/\d\mathbb P$, in the second step we prove the claim when $G\in L^\infty(\Omega,\mathcal F_\infty,\mathbb P)$, in the third step we prove the variational representation in the more general case of functionals $G$ which are bounded from above.\\
\noindent$\it{Step\,\,1.}$\\
Letting $M>0$ denote a constant such that $|\phi_{(s,x)}|\le M$ $\d s\nu(\d x)\d\mathbb P$-almost everywhere, by \eqref{eq:SecondMomentOfEt} we have
\[
	\mathbb E\bigl[\mathcal E_t(\phi)^2\bigr]\leq 1+M^2\int_0^t\mathbb E\bigr[\mathcal E_s(\phi)^2\bigr]\biggl\|\int_Ep(\pi)_{(s,x)}\,\nu(\d x)\biggr\|_{L^\infty(\Omega,\mathcal F_\infty,\mathbb P)}\d s,
	\quad\text{$t\in\R_+$.}
\]
Let $T>0$ be arbitrarily fixed. Again by \eqref{eq:SecondMomentOfEt} we have that $s\mapsto\mathbb E\bigl[\mathcal E_s(\phi)^2\bigr]$ is non-decreasing and continuous on $[0,T]$ and by
\eqref{eq:Kgrande}
\[
s\mapsto\biggl\|\int_Ep(\pi)_{(s,x)}\,\nu(\d x)\biggr\|_{L^\infty(\Omega,\mathcal F_\infty,\mathbb P)}
\]
is integrable on $[0,T]$. Therefore, by Gr\"onwall's lemma
\begin{equation*}
	\mathbb E\bigl[\mathcal E_t(\phi)^2\bigr]
	\leq\exp\left(M^2\int_0^T\biggl\|\int_Ep(\pi)_{(s,x)}\,\nu(\d x)\biggr\|_{L^\infty(\Omega,\mathcal F_\infty,\mathbb P)}\,\d s\right)\leq\mathrm{e}^{M^2 K},\quad\text{for any $t\in\R_+$}
\end{equation*}
and so $\sup_{t\in\R_+}\mathbb E\bigl[\mathcal E_t(\phi)^2\bigr]\leq\mathrm e^{M^2K}$.
By this latter relation, \eqref{eq:SDEOfE}, Proposition~\ref{bremaudmassoulie}-$(ii)$ and \eqref{eq:Kgrande}, we have
\begin{align*}
	\mathbb E\bigl[\bigl(\mathcal E_{t+h}(\phi)-\mathcal E_t(\phi)\bigr)^2\bigr]
	&=\mathbb E\biggl[\int_{(t,t+h]\times E}\phi_{(s,x)}^2\mathcal E_s(\phi)^2p(\pi)_{(s,x)}\,\d s\nu(\d x)\biggr]\\
	&\le M^2\int_t^{t+h}\mathbb E\bigr[\mathcal E_s(\phi)^2\bigr]\biggl\|\int_Ep(\pi)_{(s,x)}\,\nu(\d x)\biggr\|_{L^\infty(\Omega,\mathcal F_\infty,\mathbb P)}\d s\\
	&\le M^2\mathrm{e}^{M^2 K}
	\int_t^{t+h}\biggl\|\int_Ep(\pi)_{(s,x)}\,\nu(\d x)\biggr\|_{L^\infty(\Omega,\mathcal F_\infty,\mathbb P)}\d s\\
	&\xrightarrow[t,h\to\infty]{}0,
\end{align*}
and thus $\mathcal E_t(\phi)$ converges in $L^2(\Omega,\mathcal F_\infty,\mathbb P)$ to a random variable $X$. Letting $\mathcal{E}_\infty(\phi)$ denote the $\mathbb P$-almost sure limit of $\mathcal{E}_t(\phi)$ as $t\to\infty$, we necessarily have $X=\mathcal{E}_\infty(\phi)$ almost surely, and so
\begin{equation}\label{eq:L2einfty}
	\mathcal E_\infty(\phi)\in L^2(\Omega,\mathcal F_\infty,\mathbb P).
\end{equation}
By the martingale property it follows that $\mathcal E_t(\phi)=\mathbb E\bigl[\mathcal E_n(\phi)\big| \mathcal F_t\bigr]$ for any $t\in\R_+$ and any integer $n>t$.
By the $L^2$-convergence of $\mathcal{E}_t(\phi)$ to $\mathcal{E}_\infty(\phi)$ we easily have that $\mathbb{E}\bigl[\mathcal{E}_n(\phi)\big| \mathcal{F}_t\bigr]$ converges to $\mathbb{E}\bigl[\mathcal{E}_\infty(\phi)\big| \mathcal{F}_t\bigr]$ in $L^1$ as $n\to\infty$.
This convergence holds almost surely for a suitable subsequence $\{n'\}$ and passing to the limit as $n'\to\infty$ in the equality $\mathcal E_t(\phi)=\mathbb E\bigl[\mathcal E_{n'}(\phi)\big| \mathcal F_t\bigr]$ we get $\mathcal E_t(\phi)=\mathbb E\bigl[\mathcal E_\infty(\phi)\big| \mathcal F_t\bigr]$.
By this relation we finally deduce $\d\mathbb P_{\phi}/\d\mathbb P=\mathcal E_\infty(\phi)$.\\
\noindent$\it{Step\,\,2.}$\\
Using the elementary inequality $|(1+x)\log(1+x)-x|\leq x^2/2$, $x>-1$, for any $\phi\in\mathcal H$, we have
\begin{align}
	\mathbb E_\phi[|L(\phi)|]
	&=\mathbb E\biggl[\mathcal E_\infty(\phi)\Bigl|\int_{\R_+\times E}\bigl((1+\phi_{(s,x)})\log\bigl(1+\phi_{(s,x)}\bigr)-\phi_{(s,x)}\bigr)\,p(\pi)_{(s,x)}\,\d s\nu(\d x)\Bigr|\biggr]\nonumber\\
	&\le\frac12\,\mathbb E\biggl[\mathcal E_\infty(\phi)\int_{\R_+\times E}|\phi_{(s,x)}|^2\,p(\pi)_{(s,x)}\,\d s\nu(\d x)\biggr]\nonumber\\
	&\le\frac12\,\|\mathcal E_\infty(\phi)\|_{L^2(\Omega,\mathcal F_\infty,\mathbb P)}\biggl\|\int_{\R_+\times E}|\phi_{(s,x)}|^2\,p(\pi)_{(s,x)}\,\d s\nu(\d x)\biggr\|_{L^2(\Omega,\mathcal F_\infty,\mathbb P)}\nonumber\\
	&<\infty,\label{eq:LPhiIsInL1}
\end{align}
where the finiteness of the $L^2$-norms follows by \eqref{eq:BoundOnPhiFromMathcalH} and \eqref{eq:L2einfty}. Let $\phi\in\mathcal{H}$ and $T>0$ be arbitrarily fixed,
and set $\psi_{(s,x)}:=\log (1+\phi_{(s,x)})$, $(s,x)\in [0,T]\times E$. Clearly $\psi$ is predictable.
Additionally, we note that there exists $C>0$ such that $\bigl|\log\bigl(1+\phi_{(s,x)}\bigr)\bigr|\le C|\phi_{(s,x)}|$, $(s, x)\in[0,T]\times E$ (since there exists $c_\phi>-1$ such that $\phi_{(s,x)}\ge c_\phi$).
Therefore, by
\eqref{eq:BoundOnPhiFromMathcalH} we have \eqref{eq:IntegrabilityAssumptionOnPsi}.
Note that
\begin{align}
	\mathbb E_\phi\bigl[\log\bigl(\mathcal E_\infty(\phi)\bigr)\bigr]
	&=\mathbb E_\phi\biggl[\int_{\R_+\times E}\log\bigl(1+\phi_{(s,x)}\bigr)(N(\d s\times\d x)-p(\pi)_{(s,x)}\d s\nu(\d x))\nonumber\\
	&\,\,\,\,\,\,
	\,\,\,\,\,\,
	\,\,\,\,\,\,
	\,\,\,\,\,\,		
	+\int_{\R_+\times E}\bigl(\log\bigl(1+\phi_{(s,x)}\bigr)-\phi_{(s,x)}\bigr)\,p(\pi)_{(s,x)}\,\d s\nu(\d x)\biggr]\nonumber\\
	&=\mathbb E_\phi[L(\phi)]+\mathbb E_\phi\biggl[\int_{\R_+\times E}\log\bigl(1+\phi_{(s,x)}\bigr)(N(\d s\times\d x)-p(\pi)_{(s,x)}\d s\nu(\d x))\nonumber\\
	&\,\,\,\,\,\,
	\,\,\,\,\,\,
	\,\,\,\,\,\,
	\,\,\,\,\,\,
	-\int_{\R_+\times E}\phi_{(s,x)}\log\bigl(1+\phi_{(s,x)}\bigr)\,p(\pi)_{(s,x)}\,\d s\nu(\d x)\biggr]\nonumber\\
	&=\mathbb E_\phi[L(\phi)],
\label{eq:EfiLfi}
\end{align}
where we have used that the martingale provided by Proposition~\ref{prop:2} has null mean.
By Jensen's inequality and this relation
\begin{align}
	-\log\bigl(\mathbb E\bigl[\mathrm e^{-G}\bigr]\bigr)&=-\log\Bigl(\mathbb E_\phi\bigl[\exp\bigl(-G-\log\bigl(\mathcal E_\infty(\phi)\bigr)\bigr)\bigr]\Bigr)\nonumber\\
	&\le\mathbb E_\phi\bigl[G+\log\bigl(\mathcal E_\infty(\phi)\bigr)\bigr]=\mathbb E_\phi[G+L(\phi)]<\infty,\label{eq:LastEqualityInVariational}
\end{align}
and so
\begin{equation}\label{eq:Iineq}
	-\log\bigl(\mathbb E\bigl[\mathrm e^{-G}\bigr]\bigr)
	\le\inf_{\phi\in\mathcal H}\mathbb E_\phi\bigl[G+L(\phi)\bigr].
\end{equation}
Setting $F:=\mathrm{e}^{-G}$, we have
\[
	0<\mathrm{e}^{-\|G\|_{L^{\infty}(\Omega,\mathcal{F}_\infty,\mathbb P)}}\leq F\leq\mathrm{e}^{\|G\|_{L^\infty(\Omega,\mathcal{F}_\infty,\mathbb P)}},\quad\text{$\mathbb P$-almost surely.}
\]
Therefore, by Proposition~\ref{prop:1}
\[
	\mathcal{E}_t(\phi^{(F)})=\frac{\mathbb E\bigl[F\big| \mathcal{F}_t\bigr]}{\mathbb E[F]},\quad\forall t\in\R_+
\]
and letting $t$ go to infinity we deduce
\[
	\mathcal{E}_\infty(\phi^{(F)})=\frac{\mathbb E\bigl[F\big| \mathcal{F}_\infty\bigr]}{\mathbb E[F]}=\frac{F}{\mathbb E[F]},
\]
where $\phi_{(t,x)}^{(F)}:=\varphi_{(t,x)}^{(F)}/p(F)_{(t,x)}\in\mathcal{H}$, $(t,x)\in\R_+\times E$.
Therefore
\begin{align}
	-\log\bigl(\mathbb E\bigl[\mathrm e^{-G}\bigr]\bigr)&=-\mathbb E_{\phi^{(F)}}\bigl[\log\bigl(F\mathcal E_\infty(\phi^{(F)})^{-1}\bigr)\bigr]\nonumber\\
	&=\mathbb E_{\phi^{(F)}}\bigl[G+\log\bigl(\mathcal E_\infty(\phi^{(F)})\bigr)\bigr]\nonumber\\
	&=\mathbb E_{\phi^{(F)}}\bigl[G+L(\phi^{(F)})\bigr],\label{eq:comb}
\end{align}
where the latter equality follows by \eqref{eq:EfiLfi}.
Combining \eqref{eq:comb} and \eqref{eq:Iineq} we deduce
\[
	-\log\bigl(\mathbb E\bigl[\mathrm e^{-G}\bigr]\bigr)
	=\inf_{\phi\in\mathcal H}\mathbb E_\phi\bigl[G+L(\phi)\bigr]
\]
and the infimum is attained at $\phi^{(F)}$.
It remains to show that the infimum is uniquely attained at $\phi^{(F)}$, i.e. if $\phi\in\mathcal H$ is a stochastic process at which the infimum is attained then necessarily $\phi_{(t,x)}(\omega)=\phi_{(t,x)}^{(F)}(\omega)$ for $p(\pi)_{(t,x)}(\omega)\d t\nu(\d x)\d\mathbb{P}(\omega)$-almost all $(t,x,\omega)$.
So let $\phi\in\mathcal{H}$ be a process at which the infimum is attained.
Then Jensen's inequality \eqref{eq:LastEqualityInVariational}
holds as an equality and therefore we have
\[
	\exp\bigl(-G-\log\bigl(\mathcal E_\infty(\phi)\bigr)\bigr)=\mathbb{E}[F],\quad\text{$\mathbb{P}_\phi$-almost surely.}
\]
Similarly
\[
	\exp\bigl(-G-\log\bigl(\mathcal E_\infty(\phi^{(F)})\bigr)\bigr)=\mathbb{E}[F],\quad\text{$\mathbb{P}_\phi$-almost surely.}
\]
Therefore $\mathcal E_\infty(\phi)=\mathcal E_\infty(\phi^{(F)})$ $\mathbb P_\phi$-almost surely.
Since the probability measures $\mathbb P_\phi$ and $\mathbb P$ are equivalent it follows that $\mathcal E_\infty(\phi)=\mathcal E_\infty(\phi^{(F)})$ $\mathbb P$-almost surely.
Consequently, for any $t\in\R_+$ we have
$\mathcal E_t(\phi)=\mathcal E_t(\phi^{(F)})$ $\mathbb P$-almost everywhere,
and so by \eqref{eq:SDEOfE} we find
\begin{equation*}
	\int_{[0,t]\times E}\mathcal{E}_{s^-}(\phi)\bigl(\phi_{(s,x)}-\phi_{(s,x)}^{(F)}\bigr)(N(\d s\times\d x)-p(\pi)_{(s,x)}\d s\nu(\d x))=0.
\end{equation*}
Taking the expectation of the square of this quantity, by Proposition \ref{prop:ExtendedIsometry} we have
\begin{equation*}
	\mathbb E\biggl[\int_{[0,t]\times E}\mathcal{E}_{s^-}(\phi)^2\bigl(\phi_{(s,x)}-\phi_{(s,x)}^{(F)}\bigr)^2\,p(\pi)_{(s,x)}\d s\nu(\d x)\biggr]=0.
\end{equation*}
Since $\mathcal{E}_{s^-}(\phi)>0$ $\d s\d\mathbb P$-almost surely, we get $\phi=\phi^{(F)}$ in $\mathcal P_2(p(\pi))$, and the proof is complete.\\
\noindent$\it{Step\,\,3.}$\\
For any integer $n\geq 1$, define $G^{(n)}:=\max\{G,-n\}$.
Since $G$ is upper bounded we have $G^{(n)}\in L^\infty(\Omega,\mathcal F_\infty,\mathbb P)$.
Therefore, by Step~2 we have
\begin{equation*}
	-\log\bigl(\mathbb E\bigl[\mathrm e^{-G^{(n)}}\bigr]\bigr)=\inf_{\phi\in\mathcal{H}}\mathbb E_\phi\bigl[G^{(n)}+L(\phi)\bigr]
	\geq\inf_{\phi\in\mathcal H}\mathbb E_\phi\bigl[G+L(\phi)\bigr].
\end{equation*}
In addition, since the sequence $\{\mathrm{e}^{-G^{(n)}}\}_{n\geq 1}$ is non-decreasing we get
\[
	\lim_{n\to\infty}\mathbb E\bigl[\mathrm e^{-G^{(n)}}\bigr]=\mathbb E\bigl[\mathrm e^{-G}\bigr]
\]
by the monotone convergence theorem, hence
taking the limit as $n\to\infty$ we obtain
\begin{equation*}
	-\log\bigl(\mathbb E\bigl[\mathrm e^{-G}\bigr]\bigr)
	\geq\inf_{\phi\in\mathcal H}\mathbb E_\phi\bigl[G+L(\phi)\bigr].
\end{equation*}
For the reversed inequality we note that, for any $\psi\in\mathcal{H}$,
\begin{equation}
	-\log\bigl(\mathbb E\bigl[\mathrm e^{-G}\bigr]\bigr)=\lim_{n\to\infty}\inf_{\phi\in\mathcal H}\mathbb E_\phi\bigl[G^{(n)}+L(\phi)\bigr]
	\leq\lim_{n\to\infty}\mathbb E_\psi\bigl[G^{(n)}+L(\psi)\bigr]
	=\mathbb E_\psi\bigl[G+L(\psi)\bigr]\label{eq:reverse}
\end{equation}
where the latter equality follows by the monotone convergence theorem since the sequence $\{G^{(n)}\}_{n\geq 1}$ is non-increasing in $n$ and each $G^{(n)}$ is (for $n$ large enough) bounded above by a same constant.
Taking the infimum on $\mathcal H$ in \eqref{eq:reverse} yields the reversed inequality and the proof is complete.
\\
\noindent$\square$

\noindent{\it Proof\,\,of\,\,Proposition\,\,\ref{prop:1}.}
Since $N$ has stochastic intensity $p(\pi)$, the stochastic process
\[
	\left\{\int_{[0,t]\times E}\varphi_{(s,x)}^{(G)}\,(N(\d s\times\d x)-p(\pi)_{(s,x)}\,\d s\nu(\d x))\right\}_{t\in\R_+}
\]
is an $\mathcal{F}$-martingale, and by Theorem~\ref{thm:clarkoperative} we have
\begin{equation}\label{eq:clarkoperative2bis}
	\mathbb E\bigl[G\big| \mathcal{F}_t\bigr]=\mathbb{E}[G]+\int_{[0,t]\times E}\varphi_{(s,x)}^{(G)}\,(N(\d s\times\d x)-p(\pi)_{(s,x)}\,\d s\nu(\d x)),\quad\text{$\mathbb{P}$-almost surely.}
\end{equation}
Letting $\phi$ be defined by \eqref{eq:fi} and suppressing the dependence of $\phi$ on $G$ for ease of notation we have, since $0<c_0\leq G\leq c_1$,
\begin{align}
	\phi_{(t,x)}=\frac{\varphi_{(t,x)}^{(G)}}{p(G)_{(t,x)}}&=\frac{p(\pi G^+)_{(t,x)}-p(G)_{(t,x)}p(\pi)_{(t,x)}}{p(G)_{(t,x)} p(\pi)_{(t,x)}}\nonumber\\
	&=\frac{p(\pi G^+)_{(t,x)}}{p(G)_{(t,x)} p(\pi)_{(t,x)}}-1\nonumber\\
	&\geq\frac{c_0}{c_1}-1>-1,\quad\text{$\mathbb P$-almost surely},
        \nonumber
\end{align}
and $\phi$ is predictable and bounded with
$\sup_{(t,x)\in\R_+\times E}|\phi_{(t,x)}|\leq (c_1/c_0)+1$, $\mathbb P$-almost surely.
Set
\[
	X_t:=\int_{[0,t]\times E}\phi_{(s,x)}(N(\d s\times\d x)-p(\pi)_{(s,x)}\,\d s\nu(\d x)),\quad t\in\R_+.
\]
By Proposition~\ref{bremaudmassoulie}-$(ii)$ we have
\begin{equation}\label{eq:efi2fin}
	\sup_{t\in\R_+}\mathbb{E}\bigl[X_t^2\bigr]=\mathbb{E}\biggl[\int_{\R_+\times E}|\phi_{(s,x)}|^{2}p(\pi)_{(s,x)}\,\d s\nu(\d x)\biggr]<\infty,
\end{equation}
where the finiteness of the latter quantity follows noticing that $|\phi_{(t,x)}|\le|\varphi_{(t,x)}^{(G)}|/c_0$, $t\in\R_+$, $x\in E$, and $\varphi^{(G)}\in\mathcal P_2(p(\pi))$.
So $\{X_t\}_{t\in\R_+}$ is a square-integrable $\mathcal{F}$-martingale.
By the definition of $\phi$, the relation $p(G)_t=\mathbb E\bigl[G\big| \mathcal{F}_{t^-}\bigr]$ $\mathbb P$-almost surely and \eqref{eq:clarkoperative2bis}, we have
\begin{equation}\label{eq:exp}
	\mathbb E\bigl[G\big| \mathcal{F}_{t}\bigr]
	=\mathbb{E}[G]+\int_0^t\mathbb{E}\bigl[G\big| \mathcal{F}_{s^-}\bigr]\,\d X_s,\quad\text{$\mathbb{P}$-almost surely}.
\end{equation}
Note that
\[
	(c_1/c_0)+1\geq\Delta X_s:=X_s-X_{s^-}\geq (c_0/c_1)-1>-1,
\]
and so $X$ is of finite variation (and c\`adl\`ag).
Therefore, its quadratic variation process is given by
\begin{equation}
\label{eq:QuadraticVariationOfXComputation}
	[X,X]_t=\sum_{0<s\leq t}(\Delta X_s)^2=\int_{[0,t]\times E}\phi_{(s,x)}^2N(\d s\times\d x)
\end{equation}
and thus the path-by-path continuous part of $[X,X]$ is equal to zero. By \eqref{eq:exp} and Theorem~37 p.~84 in \cite{protter}, we have
\begin{align}
	\frac{\mathbb E\bigl[G\big| \mathcal{F}_t\bigr]}{\mathbb{E}[G]}&=\exp(X_t)\prod_{0<s\leq t}(1+\Delta X_s)\exp(-\Delta X_s)\nonumber\\
	&=\exp\biggl(
        X_t+\sum_{0<s\leq t}(\log(1+\Delta X_s)-\Delta X_s)\biggr).
        \label{Eq:1}
\end{align}
We note that
\begin{align}
	\sum_{0<s\leq t}(\log(1+\Delta X_s)-\Delta X_s)&=\int_{[0,t]\times E}(\log(1+\phi_{(s,x)})-\phi_{(s,x)})N(\d s\times\d x)\nonumber\\
	&=\int_{[0,t]\times E}(\log(1+\phi_{(s,x)})-\phi_{(s,x)})(N(\d s\times\d x)-p(\pi)_{(s,x)}\,\d s\nu(\d x))\nonumber\\
	&\,\,\,\,\,\,
	\,\,\,\,\,\,
	+\int_{[0,t]\times E}(\log(1+\phi_{(s,x)})-\phi_{(s,x)})p(\pi)_{(s,x)}\,\d s\nu(\d x)\nonumber\\
	&=-X_t+\log\mathcal{E}_t(\phi).\nonumber
\end{align}
Substituting this expression into \eqref{Eq:1} we deduce
\[
	\mathbb E\bigl[G\big| \mathcal{F}_t\bigr]=\mathbb{E}[G]\mathcal{E}_t(\phi),\quad t\in\R_+.
\]
In particular, for any $T>0$, we have $(t,x)\mapsto\ind_{[0,T]}(t)\mathcal{E}_t(\phi)\in\mathcal P_2(p(\pi))$ since
\begin{equation*}
	\mathbb E\biggl[\int_{[0,T]\times E}\mathcal{E}_s(\phi)^2p(\pi)_{(s,x)}\,\d s\nu(\d x)\biggr]
	\le\Bigl(\frac{c_1}{c_0}\Bigr)^2\mathbb E\bigl[N([0,T]\times E)\bigr]
<\infty.
\end{equation*}
Finally, we prove that \eqref{eq:BoundOnPhiFromMathcalH} holds (so that $\phi\in\mathcal H$).
First, set
\[
Y_t:=\int_{[0,t]\times E}\varphi_{(s,x)}^{(G)}(N(\d s\times\d x)-p(\pi)_{(s,x)}\,\d s\nu(\d x)),\quad t\in\R_+,
\]
and note that (by the same computation as in \eqref{eq:QuadraticVariationOfXComputation}) its quadratic variation process is
\begin{equation*}
	[Y,Y]_t=\int_{[0,t]\times E}\bigl|\varphi_{(s,x)}^{(G)}\bigr|^2N(\d s\times\d x),\quad t\in\R_+.
\end{equation*}
By Burkholder-Davis-Gundy's and Doob's inequalities, there exists a positive constant $C>0$ such that for all $t\ge0$ we have
\begin{equation*}
	\mathbb E\biggl[\biggl(\int_{[0,t]\times E}\bigl|\varphi_{(s,x)}^{(G)}\bigr|^2\,N(\d s\times\d x)\biggl)^2\biggr]
	=\mathbb E\bigl[[Y,Y]_t^2\bigr]
	\le C\mathbb E\biggl[\Bigl(\sup_{s\in[0,t]}|Y_s|\Bigr)^4\biggr]
	\le C\Bigl(\frac43\Bigr)^4\mathbb E\bigl[|Y_t|^4\bigr].
\end{equation*}
By \eqref{eq:clarkoperative2bis} the right-most term of this relation is equal to $C(4/3)^4\mathbb E\bigl[|\mathbb E\bigl[G\big| \mathcal F_t\bigr]-\mathbb E[G]|^4\bigr]$,
which is in turn less than or equal to a positive constant, say $C'>0$, which is independent of $t$.
Hence for any $t\in\R_+$ we have
\begin{align}
	&\mathbb E\biggl[\biggl(\int_{[0,t]\times E}|\phi_{(s,x)}|^2\,p(\pi)_{(s,x)}\,\d s\nu(\d x)\biggr)^2\biggr]\nonumber\\
	&\,\,\,
	\leq
	2\,\mathbb E\biggl[\biggl(\int_{[0,t]\times E}|\phi_{(s,x)}|^2\,N(\d s\times\d x)\biggr)^2\biggr]
	+2\,\mathbb E\biggl[\biggl(\int_{[0,t]\times E}|\phi_{(s,x)}|^2(N(\d s\times\d x)-p(\pi)_{(s,x)}\,\d s\nu(\d x))\biggr)^2\biggr]\nonumber\\
	&\,\,\,
	\le
	\frac{2}{c_0^4}\,\mathbb E\biggl[\biggl(\int_{[0,t]\times E}|\varphi_{(s,x)}^{(G)}|^2\,N(\d s\times\d x)\biggr)^2\biggr]
	+2\,\mathbb E\biggl[\biggl(\int_{[0,t]\times E}|\phi_{(s,x)}|^2(N(\d s\times\d x)-p(\pi)_{(s,x)}\,\d s\nu(\d x))\biggr)^2\biggr]\nonumber\\
	&\,\,\,
	\le C''
	+2\,\mathbb E\biggl[\biggl(\int_{[0,t]\times E}|\phi_{(s,x)}|^2(N(\d s\times\d x)-p(\pi)_{(s,x)}\,\d s\nu(\d x))\biggr)^2\biggr]\nonumber\\
	&\,\,\,
	=C''
	+2\,\mathbb E\biggl[\int_{[0,t]\times E}|\phi_{(s,x)}|^4p(\pi)_{(s,x)}\,\d s\nu(\d x)\biggr]\nonumber\\
	&\,\,\,
	\le C''
	+2\Bigl(\frac{c_1}{c_0}+1\Bigr)^2\,\mathbb E\biggl[\int_{\R_+\times E}|\phi_{(s,x)}|^2p(\pi)_{(s,x)}\,\d s\nu(\d x)\biggr]\nonumber\\
	&\,\,\,
	<\infty,\nonumber
\end{align}
where $C'':=(2C')/c_0^4>0$
and we have applied Proposition~\ref{bremaudmassoulie} and \eqref{eq:efi2fin}.
Taking the limit as $t$ goes to infinity in the above relations
finally yields \eqref{eq:BoundOnPhiFromMathcalH}.
\\
\noindent$\square$

\noindent{\it Proof\,\,of\,\,Proposition\,\,\ref{prop:2}.}
For $t\in\R_+$ we put
\[
	\delta(\ind_{[0,t]}\mathcal{E}(\phi)\phi):=\int_{[0,t]\times E}\mathcal{E}_{s^-}(\phi)\phi_{(s,x)}(N(\d s\times\d x)-p(\pi)_{(s,x)}\d s\nu(\d x)),
\]
\[
\delta(\ind_{[0,t]}\psi):=\int_{[0,t]\times E}\psi_{(s,x)}(N(\d s\times\d x)-p(\pi)_{(s,x)}\d s\nu(\d x)),
\]
and note that by Proposition~\ref{bremaudmassoulie}-$(iii)$ we have that
\begin{equation*}
	\delta(\ind_{[0,t]}\psi)\delta(\ind_{[0,t]}\mathcal{E}(\phi)\phi)-\int_{[0,t]\times E}\psi_{(s,x)}\mathcal{E}_s(\phi)\phi_{(s,x)}p(\pi)_{(s,x)}\d s\nu(\d x),\quad t\in\R_+,\quad\text{is an $\mathcal F$-martingale}
\end{equation*}
since \eqref{eq:IntegrabilityAssumptionOnPsi} implies $(t,x)\mapsto\ind_{[0,T]}(t)\psi_{(t,x)}\in\mathcal{P}_2(p(\pi))$ for any $T>0$, along a similar computation as in \eqref{eq:Ifin} we have $(t,x)\mapsto\ind_{[0,T]}(t)\psi_{(t,x)}\in\mathcal{P}_1(p(\pi))$ for any $T>0$, and \eqref{eq:Ifin} and \eqref{eq:IIfin} guarantee \eqref{eq:p12}.

Using the ``angle bracket'' notation,
see p.~53  in \cite{watanabe} and pp.~122-123 in \cite{protter},
we have
\begin{equation*}
	\langle Z,\mathcal{E}(\phi)\rangle_t=\langle Z,\mathcal{E}(\phi)-1\rangle_t
	=\int_{[0,t]\times E}\mathcal{E}_s(\phi)\phi_{(s,x)}\psi_{(s,x)}p(\pi)_{(s,x)}\,\d s\nu(\d x),\quad t\in\R_+,
\end{equation*}
where we put $Z_t:=\delta(\ind_{[0,t]}\psi)$ and used \eqref{eq:SDEOfE}, i.e. $\mathcal{E}_t(\phi)-1=\delta(\ind_{[0,t]}\mathcal{E}(\phi)\phi)$.
Since $\{Z_t\}_{t\in\R_+}$ is an $\mathcal{F}$-martingale under $\mathbb P$, with $Z_0=0$, and $\mathcal{E}_t(\phi)$ is the Radon-Nikodym derivative of $\mathbb P_\phi$ with respect to $\mathbb P$ on $\mathcal{F}_t$, by the Meyer-Girsanov theorem, see Theorem 36 p. 133 in \cite{protter}, we have
\begin{equation}
\label{eq:ProcessUnderPPhi}
	Z_t-\int_{[0,t]\times E}\mathcal{E}_{s^-}(\phi)^{-1}\,\d \langle Z,\mathcal{E}(\phi)\rangle_s=
	Z_t-\int_{[0,t]\times E}\phi_{(s,x)}\psi_{(s,x)}p(\pi)_{(s,x)}\,\d s\nu(\d x),\quad t\in\R_+
\end{equation}
is a local $\mathcal{F}$-martingale under $\mathbb P_\phi$.

In the following, $\mathbb E_\phi$ denotes the expectation with respect to $\mathbb P_\phi$.
In order to prove that the process defined in \eqref{eq:ProcessUnderPPhi} is an $\mathcal F$-martingale under $\mathbb P_\phi$, it suffices to prove that the stochastic process $\{(1+\phi_{(t,x)})p(\pi)_{(t,x)}\}_{(t,x)\in\R_+\times E}$ is a stochastic intensity of $N$ under $\mathbb P_\phi$ i.e., for any non-negative and predictable stochastic process $X:\R_+\times\Omega\times E\to\R_+$ we have
\begin{equation}
\label{eq:StochasticIntensityUnderPPhi}
	\mathbb{E}_\phi\biggl[\int_{\R_+\times E}X_{(t,x)}\,N(\d t\times\d x)\biggr]
	=\mathbb{E}_\phi\biggl[\int_{\R_+\times E}X_{(t,x)}(1+\phi_{(s,x)})p(\pi)_{(s,x)}\,\d s\nu(\d x)\biggr],
\end{equation}
and the integrability condition
\begin{equation}
\label{eq:IntegrabilityOfPsi1}
	\mathbb{E}_\phi\biggl[\int_{[0,T]\times E}|\psi_{(t,x)}|(1+\phi_{(s,x)})p(\pi)_{(s,x)}\,\d s\nu(\d x)\biggr]<\infty,\quad T>0.
\end{equation}
Indeed, the martingale property follows by Corollary C4 p. 235 in \cite{bremaud}.
In order to prove that the process defined in \eqref{eq:ProcessUnderPPhi} is square integrable under $\mathbb P_\phi$,
it suffices to prove
\begin{equation}
\label{eq:IntegrabilityOfPsi}
	\mathbb{E}_\phi\biggl[\int_{[0,T]\times E}|\psi_{(t,x)}|^2(1+\phi_{(s,x)})p(\pi)_{(s,x)}\,\d s\nu(\d x)\biggr]<\infty,\quad T>0.
\end{equation}
The square integrability then follows by Proposition \ref{bremaudmassoulie}-$(ii)$. We start by proving \eqref{eq:StochasticIntensityUnderPPhi}.
Setting $X_{(t,x)}:=\ind_{(a,b]}(t)\ind_A(\omega)\ind_L(x)$, $a,b\in\R_+$, $A\in\mathcal{F}_a$, $L\in\mathcal E$, and reasoning exactly as in the first part of the proof we have that the process defined by \eqref{eq:ProcessUnderPPhi} with $X$ in place of $\psi$ is a local $\mathcal F$-martingale under $\mathbb P_\phi$,
            as $X\in\mathcal P_2(p(\pi))$ since
            $\mathbb E\bigl[N((a,b]\times L)\bigr]<\infty$.
Therefore there exists a sequence of $\mathcal{F}$-stopping times $\{T_n\}_{n\ge0}$ increasing to infinity such that
\begin{equation*}
	\mathbb{E}_\phi\biggl[\int_{[0,\min\{T,T_n\}]\times E}X_{(t,x)}\,N(\d t\times\d x)\biggr]
	=\mathbb{E}_\phi\biggl[\int_{[0,\min\{T,T_n\}]\times E}X_{(t,x)}(1+\phi_{(s,x)})p(\pi)_{(s,x)}\,\d s\nu(\d x)\biggr],
\end{equation*}
for any $T\in\R_+$ and $n\geq 0$.
Letting $n$ and $T$ go to infinity in the above equation, by the monotone convergence theorem we obtain \eqref{eq:StochasticIntensityUnderPPhi} for simple predictable stochastic processes.

The result follows for a general non-negative predictable stochastic process by a standard application of the monotone class theorem, see e.g. \cite{bremaud}, Theorem~T1 p.~260.
Finally we prove \eqref{eq:IntegrabilityOfPsi1} and \eqref{eq:IntegrabilityOfPsi}.
Letting $M>0$ denote a constant such that $|\phi_{(t,x)}|\le M$ $\d t\nu(\d x)\d\mathbb P$-almost everywhere, for any $T>0$, since the probability measures $\mathbb P$ and $\mathbb P_\phi$ are equivalent we also have $|\phi_{(t,x)}|\le M$ $\d t\nu(\d x)\d\mathbb P_\phi$-almost everywhere, for any $T>0$.
Therefore,
\begin{align*}
	\mathbb{E}_\phi\biggl[&\int_{[0,T]\times E}|\psi_{(t,x)}|^2(1+\phi_{(s,x)})p(\pi)_{(s,x)}\,\d s\nu(\d x)\biggr]\\
	&\le(1+M)\,\mathbb{E}\biggl[\mathcal E_T(\phi)\int_{[0,T]\times E}|\psi_{(t,x)}|^2p(\pi)_{(s,x)}\,\d s\nu(\d x)\biggr]\\
	&\le(1+M)\,\|\mathcal E_T(\phi)\|_{L^2(\Omega,\mathcal F_\infty,\mathbb P)}\,
        \sqrt{
          \mathbb{E}\biggl[\biggl(\int_{[0,T]\times E}|\psi_{(t,x)}|^2p(\pi)_{(s,x)}\,\d s\nu(\d x)\biggr)^2\biggr]}
        \\
	&<\infty,
\end{align*}
and
\begin{align*}
	&\mathbb{E}_\phi\biggl[\int_{[0,T]\times E}|\psi_{(t,x)}|(1+\phi_{(s,x)})p(\pi)_{(s,x)}\,\d s\nu(\d x)\biggr]\\
  &\;\le
  \sqrt{(1+M)
    \mathbb{E}_\phi\biggl[\int_{[0,T]\times E}|\psi_{(t,x)}|^2(1+\phi_{(s,x)})p(\pi)_{(s,x)}\,\d s\nu(\d x)\biggr]}
  \sqrt{
    \mathbb E_\phi[N([0,T]\times E)]}
  \\
  &\;=
  \sqrt{
    (1+M)
    \mathbb{E}_\phi\biggl[\int_{[0,T]\times E}|\psi_{(t,x)}|^2(1+\phi_{(s,x)})p(\pi)_{(s,x)}\,\d s\nu(\d x)\biggr]}
  \sqrt{
    \mathbb E\bigl[\mathcal E_T(\phi)N([0,T]\times E)\bigr]}
  \\
  &\;\le
  \sqrt{(1+M)
    \mathbb{E}_\phi\biggl[\int_{[0,T]\times E}\!|\psi_{(t,x)}|^2(1+\phi_{(s,x)})p(\pi)_{(s,x)}\,\d s\nu(\d x)\biggr]}
  \mathbb E\bigl[\mathcal E_T(\phi)^2\bigr]^{1/4}\mathbb E\bigl[N([0,T]\times E)^2\bigr]^{1/4}\\
	&\;<\infty,
\end{align*}
which concludes the proof.
\\
\noindent$\square$

\subsection{Proof of Theorem \ref{thm:clarkoperative}}

The proof of the Clark-Ocone formula relies on
the following propositions, which will be shown at the end of this subsection.
The first one provides an $\overline{\mathcal{F}}$-predictable representation formula for square-integrable functionals of a marked point process with a stochastic intensity.
The second one gives a formula which allows us to transform the expectation of the product between a square-integrable functional of a marked point process with a stochastic intensity and an integral with respect to the compensated marked point process into the expectation of an integral with respect to the measure $\mathrm{d}t\nu(\mathrm{d}x)$.

\begin{Proposition}\label{Thm:clark}
	Assume that $N$ has a $\overline{\mathcal F}$-stochastic intensity $\lambda$.
	Then, for any $G\in L^2(\Omega,\mathcal{F}_\infty,\mathbb P)$, there exists $u^{(G)}\in\mathcal{P}_2^{\overline{\mathcal F}}(\lambda)$ such that, for all $t\in\R_+$,
	\begin{equation}\label{eq:clarkrepresentation}
		\mathbb E\bigl[G\big| \overline{\mathcal{F}}_t\bigr]=\mathbb{E}[G]+\int_{[0,t]\times E}u_{(s,x)}^{(G)}\,(N(\d s\times\d x)-\lambda_{(s,x)}\,\d s\nu(\d x)),\quad\text{$\mathbb{P}$-almost surely}.
	\end{equation}
	In particular,
	\begin{equation}\label{eq:clarkrepresentationTInfinity}
		G=\mathbb{E}[G]+\int_{\R_+\times E}u_{(s,x)}^{(G)}\,(N(\d s\times\d x)-\lambda_{(s,x)}\,\d s\nu(\d x)),\quad\text{$\mathbb{P}$-almost surely}.
	\end{equation}
\end{Proposition}
We remark that the integrand $u^{(G)}\in\mathcal{P}_2^{\overline{\mathcal F}}(\lambda)$ is not made explicit, in contrast with the Clark-Ocone formula.

\begin{Proposition}\label{Thm:duality}
	Assume that $N$ has a $\overline{\mathcal F}$-stochastic intensity $\lambda$ and that $\mathbb E\bigl[N([0,t]\times E)\bigr]<\infty$, for all $t\in\mathbb{R}_+$.
	Let $G\in L^2(\Omega,\mathcal{F}_\infty,\mathbb P)$ and $u\in\mathcal{P}_{1,2}^{\overline{\mathcal F}}(\lambda)$. Then we have
	\begin{equation*}
		\mathbb{E}\bigl[G\delta(u)\bigr]=\mathbb{E}\biggl[\int_{\R_+\times E}u_{(t,x)}^{(G)}u_{(t,x)}\,\lambda_{(t,x)}\,\d t\nu(\d x)\biggr].
	\end{equation*}
\end{Proposition}

\noindent{\it Proof\,\,of\,\,Theorem\,\,\ref{thm:clarkoperative}.}
We divide the proof into two steps.
In the first step we derive
a predictable representation of $G$, in the second step we identify the integrand of the predictable representation.\\
\noindent{\it Step\,\,1.}\\
By \eqref{eq:integralOfpPiFinite} we have that for all $t\in\mathbb{R}_+$, $\int_{[0,t]\times E}p(\pi)_{(s,x)}\,\d s\nu(\d x)<\infty$ $\mathbb P$-almost surely, so by Lemma~\ref{le:PapStochlink}, $p(\pi)$ is a stochastic intensity of $N$.
Hence by Lemma~\ref{le:equiv}, $p(\pi)$ is a $\overline{\mathcal F}$-stochastic intensity of $N$.
Therefore, by Proposition~\ref{Thm:clark} we have
\begin{equation*}
	G=\mathbb{E}[G]+\int_{\R_+\times E}u_{(s,x)}^{(G)}\,(N(\d s\times\d x)-p(\pi)_{(s,x)}\,\d s\nu(\d x)),\quad\text{$\mathbb{P}$-almost surely},
\end{equation*}
for some $u^{(G)}\in\mathcal{P}_2^{\overline{\mathcal{F}}}(p(\pi))$.
By Proposition~\ref{cor:PredictableProjection} there exists an $\mathcal F$-predictable stochastic process $p(u^{(G)})$ such that $p(u^{(G)})_{(t,x)}=\mathbb{E}\bigl[u_{(t,x)}^{(G)}\big| \mathcal{F}_{t^{-}}\bigr]$ $\mathbb P$-almost surely.
By relation \eqref{eq:barnotbar} we deduce
\begin{equation}\label{eq:puG}
	p(u^{(G)})_{(t,x)}=\mathbb{E}\bigl[u_{(t,x)}^{(G)}\big| \overline{\mathcal{F}}_{t^{-}}\bigr]=u_{(t,x)}^{(G)},\quad\text{$\mathbb P$-almost surely},
\end{equation}
where the latter equality follows by the $\overline{\mathcal{F}}$-predictability of $u^{(G)}$ which guarantees that $u_{(t,x)}^{(G)}$ is $\overline{\mathcal{F}}_{t^{-}}$-measurable (as already noticed in Section~\ref{sec:Preliminaries}, this follows by an obvious modification of the proof of Lemma~A3.3.I p.~425 in \cite{daley}).
So finally, we deduce the predictable representation
\begin{equation}\label{eq:clarkpred}
	G=\mathbb{E}[G]+\int_{\R_+\times E}p(u^{(G)})_{(s,x)}\,(N(\d s\times\d x)-p(\pi)_{(s,x)}\,\d s\nu(\d x)),\quad\text{$\mathbb{P}$-almost surely}.
\end{equation}
\noindent{\it Step\,\,2.}\\
Let $u:\R_+\times\Omega\times E\to\R$ be a predictable stochastic process and, for any $n\geq 0$, define the (predictable) stochastic process $u^{(n)}$ by $u_{(s,x)}^{(n)}:=u_{(s,x)}\ind_{\{u_{(s,x)}\in[-n,n]\}}$.
For any $n\geq 0$ and $t\in\R_+$, by the square integrability of $G$ and $N([0,t]\times E)$ and \eqref{eq:IntegrabilityOfPi}, we have
\begin{align*}
	\mathbb E\biggl[\int_{\R_+\times E}|u_{(s,x)}^{(n)}\ind_{[0,t]}(s)G|\pi_{(s,x)}\,\d s\nu(\d x)\biggr]
	&\le n\mathbb E\biggl[|G|\int_{[0,t]\times E}\pi_{(s,x)}\,\d s\nu(\d x)\biggr]\\
	&\le n\|G\|_{L^2(\Omega,\mathcal{F}_\infty,\mathbb P)}\Bigl\|\int_{[0,t]\times E}\pi_{(s,x)}\,\d s\nu(\d x)\Bigr\|_{L^2(\Omega,\mathcal{F}_\infty,\mathbb P)}<\infty
\end{align*}
and
\begin{align*}
	\mathbb E\biggl[\int_{\R_+\times E}|u_{(s,x)}^{(n)}\ind_{[0,t]}(s)G^+_{(s,x)}|\pi_{(s,x)}\,\d s \nu(\d x)\biggr]
	&\le n\mathbb E\biggl[\int_{[0,t]\times E}|G^+_{(s,x)}|\pi_{(s,x)}\,\d s\nu(\d x)\biggr]\\
	&=n\mathbb E\biggl[\int_{[0,t]\times E}|G(N+\varepsilon_{(s,x)})|\pi_{(s,x)}(N)\,\d s\nu(\d x)\biggr]\\
	&=n\mathbb E\biggl[\int_{[0,t]\times E}|G(N)|N(\d s\times\d x)\biggr]\\
	&=n\mathbb E\bigl[|G|N([0,t]\times E)\bigr]\\
	&\le n\|G\|_{L^2(\Omega,\mathcal{F}_\infty,\mathbb P)}\|N([0,t]\times E)\|_{L^2(\Omega,\mathcal{F}_\infty,\mathbb P)}<\infty,
\end{align*}
thus both conditions in \eqref{eq:condint} are verified for the predictable mapping $X_{(s,x)}^{(n)}:=u_{(s,x)}^{(n)}\ind_{[0,t]}(s)$ and the random variable $F:=G$.
The stochastic process $X^{(n)}$ defined above is clearly in $\mathcal P_{1,2}^{\mathcal F}(p(\pi))\subset\mathcal P_{1,2}^{\overline{\mathcal F}}(p(\pi))$, since it is bounded and $\mathbb E\bigl[N([0,t]\times E)\bigr]<\infty$. Hence by Lemma~\ref{le:IBP} (with $X:=X^{(n)}$), Proposition~\ref{Thm:duality} (with $X:=X^{(n)}$) and \eqref{eq:puG}, we have
\begin{align}
    &\mathbb E\left[\int_{\R_+\times E}X_{(s,x)}^{(n)}\pi_{(s,x)}D_{(s,x)}G\,\d s\nu(\d x)\right]\nonumber\\
	&=\mathbb E\Bigl[G\Delta(X^{(n)})\Bigr]\nonumber\\
  &=\mathbb E\left[G\int_{[0,t]\times E}u_{(s,x)}^{(n)}(N(\d s\times\d x)-p(\pi)_{(s,x)}\,\d s\nu(\d x))\right]
    +\mathbb E\left[G\int_{[0,t]\times E}u_{(s,x)}^{(n)}(p(\pi)_{(s,x)}-\pi_{(s,x)})\,\d s\nu(\d x)\right]\nonumber\\
	&=\mathbb E\left[\int_{[0,t]\times E}u_{(s,x)}^{(n)}u_{(s,x)}^{(G)}p(\pi)_{(s,x)}\,\d s\nu(\d x)\right]+
	\mathbb E\left[G\int_{[0,t]\times E}u_{(s,x)}^{(n)}(p(\pi)_{(s,x)}-\pi_{(s,x)})\,\d s\nu(\d x)\right]\nonumber\\
    &=\mathbb E\left[\int_{[0,t]\times E}u_{(s,x)}^{(n)}p(u^{(G)})_{(s,x)}p(\pi)_{(s,x)}\,\d s\nu(\d x)\right]+
	\mathbb E\biggl[G\int_{[0,t]\times E}u_{(s,x)}^{(n)}(p(\pi)_{(s,x)}-\pi_{(s,x)})\,\d s\nu(\d x)\biggr].\nonumber
\end{align}
Therefore we have
\begin{align}
	0&=\mathbb E\biggl[\int_{[0,t]\times E}u_{(s,x)}^{(n)}\bigl(\pi_{(s,x)}D_{(s,x)}G-p(u^{(G)})_{(s,x)}p(\pi)_{(s,x)}-G p(\pi)_{(s,x)}+G\pi_{(s,x)}\bigr)\,\d s\nu(\d x)\biggr]\nonumber\\
    &=\mathbb E\biggl[\int_{[0,t]\times E}u_{(s,x)}^{(n)}
    \bigl(\mathbb E\bigl[\pi_{(s,x)}G^+_{(s,x)}\big| \mathcal F_{s^-}\bigr]-p(u^{(G)})_{(s,x)}p(\pi)_{(s,x)}-\mathbb E\bigl[G\big| \mathcal F_{s^-}\bigr]p(\pi)_{(s,x)}\bigr)\,\d s\nu(\d x)\biggr].\nonumber
\end{align}
As discussed in Remark~\ref{rem:RemarkFollowingCOGeneralized} the predictable projections $p(\pi G^+)$ and $p(G)$ exist and so
\begin{equation*}
	\mathbb E\biggl[\int_{[0,t]\times E}u_{(s,x)}^{(n)}\bigl(p(\pi G^+)_{(s,x)}-p(u^{(G)})_{(s,x)}p(\pi)_{(s,x)}-p(G)_{(s,x)}p(\pi)_{(s,x)}\bigr)\,\d s\nu(\d x)\biggr]=0,
\end{equation*}
for all $n\ge0$ and $t\in\R_+$.
As shown in Remark~\ref{rem:RemarkFollowingCOGeneralized} the stochastic process $\varphi^{(G)}$ in \eqref{eq:ExplicitIntegrand} is well-defined and predictable.
So the above relation can be rewritten as
\begin{equation}\label{eq:normordensity}
	\mathbb E\biggl[\int_0^t\int_Eu_{(s,x)}^{(n)}\bigl(\varphi_{(s,x)}^{(G)}-p(u^{(G)})_{(s,x)}\bigr)p(\pi)_{(s,x)}\,\d s\nu(\d x)\biggr]=0,\quad\forall n\ge0,\ \forall t\in\mathbb{R}_+.
\end{equation}
Since $u$ is an arbitrary predictable stochastic process, by choosing $u=\varphi^{(G)}-p(u^{(G)})$, equation~\eqref{eq:normordensity} reads as
\begin{equation*}
	\mathbb E\left[\int_0^t\int_E\bigl(\varphi_{(s,x)}^{(G)}-p(u^{(G)})_{(s,x)}\bigr)^2\ind_{\{\varphi_{(s,x)}^{(G)}-p(u^{(G)})_{(s,x)}\in[-n,n]\}}
    p(\pi)_{(s,x)}\,\d s\nu(\d x)\d s\right]=0,
\end{equation*}
for all $n\ge0$ and $t\in\R_+$.
By the monotone convergence theorem, letting $n$ tend to infinity, we get
\[
	\mathbb E\biggl[\int_0^t\int_E\bigl(\varphi_{(s,x)}^{(G)}-p(u^{(G)})_{(s,x)}\bigr)^2 p(\pi)_{(s,x)}\,\d s\nu(\d x)\biggr]=0,\quad\text{for all $t\in\mathbb{R}_+$}
\]	which implies that $p(u^{(G)})_{(s,x)}(\omega)=\varphi_{(s,x)}^{(G)}(\omega)$ $\mathbb P(\mathrm{d}\omega)p(\pi)_{(s,x)}(\omega)\,\d s\nu(\d x)$-almost everywhere on $(\R_+\times\Omega\times E,\mathcal{B}(\R_+)\otimes\mathcal{F}_\infty\otimes\mathcal E)$.
Combining this with \eqref{eq:clarkpred} we finally deduce \eqref{eq:clarkoperative2}.
Since $u^{(G)}$ is square integrable with respect to $\mathbb P(\mathrm{d}\omega)p(\pi)_{(s,x)}(\omega)\,\d s\nu(\d x)$ by Jensen's inequality we easily have that $p(u^{(G)})$ is square integrable with respect to $\mathbb P(\mathrm{d}\omega)p(\pi)_{(s,x)}(\omega)\,\d s\nu(\d x)$ and so $\varphi^{(G)}\in\mathcal P_2^{\mathcal F}(p(\pi))$.
\noindent$\square$

\noindent$\it{Proof\,\,of\,\,Proposition\,\,\ref{Thm:clark}.}$
Define the $\overline{\mathcal F}$-martingale $G_t:=\mathbb{E}\bigl[G\big| \overline{\mathcal{F}}_t\bigr]$, $t\in\R_+$.
Since $G$ is square integrable, by Jensen's inequality we have $\sup_{t\in\R_+}\mathbb{E}\bigl[G_t^2\bigr]\le\mathbb E\bigl[G^2\bigr]<\infty$. Then $\{G_t\}_{t\geq 0}$ is square integrable and \eqref{eq:clarkrepresentation} follows by Theorem~T8~p.~239 of \cite{bremaud} with the completed filtration, noticing that $\mathbb E\bigl[G\big| \overline{\mathcal{F}}_0\bigr]=\mathbb E\bigl[G\big| \mathcal{F}_0\bigr]=\mathbb E[G]$.
The first equality of this latter relation follows by \eqref{eq:barnotbar} and the second inequality follows by $\mathcal{F}_0=\{\emptyset,\Omega\}$.
\\
\noindent$\square$

\noindent{\it Proof\,\,of\,\,Proposition\,\,\ref{Thm:duality}.}
By Proposition~\ref{bremaudmassoulie}-$(i)$ we have $\mathbb{E}\bigl[\delta(u)\bigr]=0$.
By Proposition~\ref{Thm:clark} we have $G-\mathbb{E}[G]=\delta(u^{(G)})$, for some $u^{(G)}\in\mathcal{P}_2^{\overline{\mathcal F}}(\lambda)$.
Therefore
\begin{equation*}
	\mathbb{E}\bigl[G\delta(u)\bigr]
	=\mathbb{E}\bigl[(G-\mathbb{E}[G])\delta(u)\bigr]=\mathbb{E}\left[\delta(u^{(G)})\delta(u)\right]
	=\mathbb{E}\left[\int_{\R_+\times E}u_{(t,x)}^{(G)}u_{(t,x)}\,\lambda_{(t,x)}\,\d t\nu(\d x)\right],
\end{equation*}
where the latter equality follows by Proposition \ref{prop:ExtendedIsometry}.
\\
\noindent$\square$

\section{Appendix}

\subsection{Proof of Lemma~\ref{lem:DifferentDefinitionsOfMPP}}
\label{subsec:ProofOfDifferentDefinitionsOfMPP}

We have
\begin{align*}
	1-\mathbb P(\Omega)
	&=\mathbb P(\{\omega\in\Omega'\st\exists (t,x)\in\mathrm{Supp}(\omega),\ \omega(\{t\}\times E)\ge2\})\\
		&=\mathbb E\bigl[\ind_{\{\exists (t,x)\in \mathrm{Supp}(N),\ N(\{t\}\times E)\ge2\}}\bigr]\\
	&\le\mathbb E\biggl[\int\ind_{\{N(\{t\}\times E)\ge2\}}\,N(\d t\times\d x)\biggr]\\
	&=\mathbb E\biggl[\int_{\R_+\times E}\!\!\ind_{\{(N+\varepsilon_{(t,x)})(\{t\}\times E)\ge2\}}\pi_{(t,x)}\,\d t\nu(\d x)\biggr]\\
	&=\mathbb E\biggl[\int_{\R_+\times E}\ind_{\{N(\{t\}\times E)\ge1\}}\pi_{(t,x)}\,\d t\nu(\d x)\biggr]\\
	&=0,
\end{align*}
where the last equality follows since for any $\omega\in\Omega'$, we have that $\ind_{\{\omega(\{t\}\times E)\ge1\}}=0$ $\d t$-almost surely.
Indeed, for any $\omega\in\Omega'$ and any $T>0$, we have $\omega([0,T]\times E)<\infty$, therefore the cardinal of $\{t\in[0,T]\st\omega(\{t\}\times E)\ge1\}$ is finite, and thus
\begin{equation*}
	\int_0^T\ind_{\{t\in[0,T]\st\omega(\{t\}\times E)\ge1\}}\,\d t=0.
\end{equation*}

\subsection{Proof of Proposition~\ref{cor:PredictableProjection}}
\label{subsec:ProofOfPredictableProjection}

The proof uses the following lemma which guarantees the existence of a predictable version of a bounded stochastic process $X$.

\begin{Lemma}\label{le:predproj}
	Assume that $N$ has a Papangelou conditional intensity $\pi$, i.e. \eqref{eq:GNZ} holds, and let $X\in L^\infty(\R_+\times\Omega\times E,\mathcal{B}(\R_+)\otimes\mathcal{F}_\infty\times\mathcal{E},\mathrm{d}t\nu(\mathrm{d}x)\mathbb{P}(\mathrm{d}\omega))$ be a real-valued stochastic process.
	Then there exists a predictable stochastic process $p(X):\R_+\times\Omega\times E\to\R$ such that, for all $(t,x)\in\R_+\times E$, we have
	\begin{equation*}
		p(X)_{(t,x)}=\mathbb E\bigl[X_{(t,x)}\big| \mathcal F_{t^-}\bigr],\quad\mathbb P\text{-almost surely.}
	\end{equation*}
\end{Lemma}

\noindent{\it Proof\,\,of\,\,Proposition\,\,\ref{cor:PredictableProjection}.}
If $X$ is assumed to be non-negative, we set $X^{(n)}:=\min(X,n)$, $n\ge0$.
By Lemma~\ref{le:predproj}, for any $n$ and $(t,x)$,
\[
	p(X_{(t,x)}^{(n)})=\mathbb{E}\bigl[X_{(t,x)}^{(n)}\big| \mathcal{F}_{t^-}\bigr],\quad\mathbb{P}\text{-almost surely.}
\]
By the monotone convergence theorem
\[
	\lim_{n\to\infty}p(X_{(t,x)}^{(n)})= \mathbb{E}\bigl[X_{(t,x)}\big| \mathcal{F}_{t^-}\bigr],\quad\mathbb{P}\text{-almost surely.}
\]
Since $p(X^{(n)})$ is predictable for each $n$ its limit exists and is predictable, completing the proof under the assumption that $X$ is non-negative.
	
If one assumes that, for $\mathrm{d}t\nu(\mathrm{d}x)$-almost all $(t,x)\in\R_+\times E$, $X_{(t,x)}\in L^1(\Omega,\mathcal{F}_\infty,\mathbb P)$, then we write $X=X^+-X^-$, where $X^+:=\max(X,0)$ and $X^-:=-\min(X,0)$.
Applying the first part of the proposition to $X^+$ and $X^-$ we have that there exist two predictable stochastic processes $p(X^+)$ and $p(X^-)$ such that
\begin{equation*}
	p(X^+)_{(t,x)}=\mathbb E\bigl[X_{(t,x)}^+\big| \mathcal F_{t^-}\bigr]\quad\text{and}\quad
	p(X^-)_{(t,x)}=\mathbb E\bigl[X_{(t,x)}^-\big| \mathcal F_{t^-}\bigr],\quad\mathbb P\text{-almost surely.}
\end{equation*}
By taking the expectation of these two equalities, one has, for $\mathrm{d}t\nu(\mathrm{d}x)$-almost all $(t,x)$,
\[
	p(X^+)_{(t,x)}<\infty\quad\text{and}\quad p(X^-)_{(t,x)}<\infty,
\]
and therefore
\begin{equation}
\label{eq:ComputingPredictableProjection}
	\mathbb E\bigl[X_{(t,x)}\big| \mathcal F_{t^-}\bigr]
	=\mathbb E\bigl[X_{(t,x)}^+\big| \mathcal F_{t^-}\bigr]-\mathbb E\bigl[X_{(t,x)}^-\big| \mathcal F_{t^-}\bigr]
	=p(X^+)_{(t,x)}-p(X^-)_{(t,x)},\quad\mathbb P\text{-almost surely.}
\end{equation}
	The claim follows by noticing that the right-hand side of \eqref{eq:ComputingPredictableProjection} is a predictable stochastic process.
\\
\noindent$\square$

\noindent{\it Proof\,\,of\,\,Lemma\,\,\ref{le:predproj}.}
The idea is to apply the existence part of Theorem~3.3 in \cite{last2} and the last displayed formula on p.~368 again in \cite{last2}.
Following the notation of \cite{last2}, we consider the locally compact second countable Hausdorff space	 $\mathbf X:=\R_+\times E$ and the DC-semiring $\mathcal{S}:=\{(a,b]\times B\st a,b\in\R_+,\ a<b,\ B\in\mathcal E,\ B\text{ relatively compact}\}$.
Moreover, we consider the system $\Gamma:=\{\Gamma_{(t,x)},\Gamma_S:\,\,(t,x)\in\R_+\times E,\,S\in\cal S\}$ defined by	 $\Gamma_{(t,x)}:=[0,t)\times E$ and $\Gamma_{(a,b]\times B}:=[0,a]\times E$. It is readily checked that it satisfies conditions (2.3) and (2.4) in \cite{last2}.
In this setting, for $(t,x)\in\R_+\times E$, $a,b\in\R_+$: $a<b$, $B\in\cal E$ a relatively compact set, the $\sigma$-fields $\mathcal{F}((t,x))$ and $\mathcal{F}((a,b]\times B)$ defined on p. 364 of \cite{last2} are given by
\[
	\mathcal F((t,x)):=\sigma(N|_{[0,t)\times E})=\mathcal{F}_{t^-}\quad\text{and}\quad\mathcal{F}((a,b]\times B):=\sigma(N|_{[0,a]\times E})=\mathcal{F}_a,
\]
where $N|_A$ denotes the restriction of the random measure $N$ to $A\in\mathcal{B}(\R_+)\otimes\mathcal{E}$.
One may easily check that the predictable $\sigma$-algebra $\mathcal{P}$ in \cite{last2} coincides with the $\sigma$-field $\mathcal{P}(\mathcal{F})\otimes\mathcal{E}$, and that the point process $N$ satisfies condition $(2.5)$ on p.~364 of \cite{last2}.
We note that if condition $\Sigma(\Lambda)$ on p.~367 of \cite{last2} holds, then the claim of the lemma follows by the existence part of Theorem~3.3 in \cite{last2} and the last displayed formula on p.~368 again in \cite{last2}.
In our setting, condition $\Sigma(\Lambda)$ on p.~367 of \cite{last2}
reads
\[
	\Sigma(\Lambda):\,\,\mathbb{P}(N((a,b]\times E)=0\mid\mathcal{F}_a)>0\quad\text{$\mathbb{P}$-a.s., $a,b\in\R_+$: $a<b$}.
\]
Since $N$ has a Papangelou conditional intensity it satisfies condition $(\Sigma)$ from Remark~2.5$(c)$ in \cite{georgii} which, in our setting, reads
\[
	(\Sigma):\,\,\mathbb{P}(N((a,b]\times E)=0\mid\mathcal{F}_a\vee\mathcal{F}_{(b,\infty)})>0\quad\text{$\mathbb{P}$-a.s., $a,b\in\R_+$: $a<b$},
\]
where $\mathcal{F}_{(b,\infty)}:=\sigma(N|_{(b,\infty)\times E})$.
Condition~$\Sigma(\Lambda)$ easily follows by $(\Sigma)$ and the properties of the conditional expectation, indeed
$$
	\mathbb{P}(N((a,b \hskip0.03cm ]\times E)=0\mid\mathcal{F}_a)
	=\mathbb{E}\bigl[\mathbb{P}(N((a,b \hskip0.03cm ]\times E)=0\big| \mathcal{F}_a\vee\mathcal{F}_{(b,\infty)})\bigr]>0
$$
$\mathbb{P}$-a.s., for any $a,b\in\R_+$ such that $a<b$.
\\
\noindent$\square$

\subsection{Proof of Proposition \ref{prop:ExtendedIsometry}}
\label{subsec:ProofOfExtendedIsometry}

The proof uses the following lemma which, under a mild integrability condition on $N$, guarantees that $\mathcal{P}_{1,2}(\lambda)$ is dense in $\mathcal{P}_2(\lambda)$.
Although its proof is quite standard, we include it for the sake of completeness.

\begin{Lemma}\label{lem:Density}
	Assume that $N$ has $\mathcal{G}$-stochastic intensity $\lambda$ and that $\mathbb{E}\bigl[N([0,t]\times E)\bigr]<\infty$ for any $t\in\R_+$.
	Then $\mathcal{P}_{1,2}(\lambda)$ is dense in $\mathcal{P}_2(\lambda)$.
\end{Lemma}

\noindent{\it Proof\,\,of\,\,Proposition\,\,\ref{prop:ExtendedIsometry}.}
Let $X\in\mathcal{P}_2(\lambda)$ and let $(X^{(n)})_{n\ge0}\subset\mathcal{P}_{1,2}(\lambda)$ be a sequence given by Lemma~\ref{lem:Density} which converges to $X$ in $\mathcal{P}_2(\lambda)$.
By the isometry formula Proposition~\ref{bremaudmassoulie}-$(ii)$
we have
\begin{align*}
	\mathbb E\bigl[(\delta(X^{(n)})-\delta(X^{(n+m)}))^2\bigr]
	&=\mathbb E\bigl[\delta(X^{(n)}-X^{(n+m)})^2\bigr]\\
	&=\mathbb E\biggl[\int_{\R_+\times E}(X_{(t,x)}^{(n)}-X_{(t,x)}^{(n+m)})^2\lambda_{(t,x)}\,\d t\nu(\d x)\biggr]\\
	&=\|X^{(n)}-X^{(n+m)}\|_{\mathcal P_2(\lambda)}\\
	&\xrightarrow[m,n\to\infty]{}0.
\end{align*}
Since the space $L^2(\Omega,\mathcal{F}_\infty,\mathbb P)$ is complete the sequence $(\delta(X^{(n)}))_{n\ge0}$ converges in $L^2(\Omega,\mathcal{F}_\infty,\mathbb P)$.
Its limit does not depend on the sequence $(X^{(n)})_{n\ge0}$, indeed if $(Y^{(n)})_{n\ge0}\subset\mathcal{P}_{1,2}(\lambda)$ is another sequence which converges to $X$ in $\mathcal{P}_2(\lambda)$, again by Proposition~\ref{bremaudmassoulie}-$(ii)$ we have
\begin{align*}
	\mathbb E\bigl[(\delta(X^{(n)})-\delta(Y^{(n)}))^2\bigr]
	&=\mathbb E\biggl[\int_{\R_+\times E}(X_{(t,x)}^{(n)}-Y_{(t,x)}^{(n)})^2\lambda_{(t,x)}\,\d t\nu(\d x)\biggr]\\
	&=\|X^{(n)}-Y^{(n)}\|_{\mathcal P_2(\lambda)}\\
	&\xrightarrow[n\to\infty]{}0.
\end{align*}
For an arbitrary $X\in\mathcal{P}_{2}(\lambda)$, we denote by $\delta(X)$ the limit in $L^2(\Omega,\mathcal{F}_\infty,\mathbb P)$ of $\delta(X^{(n)})$ for some sequence $(X^{(n)})_{n\geq 0}\subset\mathcal{P}_{1,2}(\lambda)$ converging to $X$.
In order to conclude, for $X,Y\in\mathcal{P}_2(\lambda)$ we let $(X^{(n)})_{n\ge0},(Y^{(n)})_{n\ge0}$ be the corresponding sequences of $\mathcal{P}_{1,2}(\lambda)$ tending to $X$ and $Y$ respectively.
By Proposition~\ref{bremaudmassoulie}-$(ii)$ we have
\begin{equation*}
	\mathbb E\bigl[\delta(X^{(n)})\delta(Y^{(n)})\bigr]=\mathbb E\biggl[\int_{\R_+\times E}X_{(t,x)}^{(n)}Y_{(t,x)}^{(n)}\lambda_{(t,x)}\,\d t\nu(\d x)\biggr],
\end{equation*}
and taking the limit in $n$ on both sides yields the conclusion.
\\
\noindent$\square$

\noindent{\it Proof\,\,of\,\,Lemma\,\,\ref{lem:Density}.}
Take $X\in\mathcal{P}_2(\lambda)$ and define $X^{(n)}$ by $X^{(n)}_{(t,x)}:=\max(\min(X_{(t,x)},n),-n)\ind_{[0,n]}(t)$ for each $n\ge0$.
It is readily seen that $X^{(n)}$ is $\mathcal G$-predictable and $|X_{(t,x)}^{(n)}|\leq n\ind_{[0,n]}(t)$ for any $t\in\R_+$, $x\in E$ and $n\in\mathbb N$.
Therefore, by \eqref{eq:classicalStochasticIntensity} and the finiteness of the first moment of $N([0,n]\times E)$, for $p\in\{1,2\}$,
\begin{equation*}
	\mathbb{E}\biggl[\int_{\R_+\times E}|X_{(t,x)}^{(n)}|^p\,\lambda_{(t,x)}\,\d t\nu(\d x)\biggr]
	\leq n^p\,\mathbb E\biggl[\int_0^n\int_E\lambda_{(t,x)}\,\d t\nu(\d x)\biggr]
	=n^p\,\mathbb{E}\bigl[N([0,n]\times E)\bigr]
	<\infty,
\end{equation*}
and so $X^{(n)}\in\mathcal{P}_{1,2}(\lambda)$. Additionally we have
\begin{align}
	\mathbb{E}\biggl[\int_{\R_+\times E}|X_{(t,x)}^{(n)}-X_{(t,x)}|^2\,\lambda_{(t,x)}\,\d t\nu(\d x)\biggr]&
	=\mathbb{E}\biggl[\int_{\R_+\times E}|X_{(t,x)}^{(n)}-X_{(t,x)}|^2\ind_{\{|X_{(t,x)}|\leq n\}}\,\lambda_{(t,x)}\,\d t\nu(\d x)\biggr]\nonumber\\
	&\qquad+\mathbb{E}\biggl[\int_{\R_+\times E}|X_{(t,x)}^{(n)}-X_{(t,x)}|^2\ind_{\{|X_{(t,x)}|>n\}}\,\lambda_{(t,x)}\,\d t\nu(\d x)\biggr]\nonumber\\
	&\leq\mathbb{E}\biggl[\int_{n}^\infty\int_E|X_{(t,x)}|^2\ind_{\{|X_{(t,x)}|\leq n\}}\,\lambda_{(t,x)}\,\d t\nu(\d x)\biggr]\nonumber\\
	&\qquad+2\mathbb{E}\biggl[\int_{\R_+\times E}|X_{(t,x)}^{(n)}|^2\ind_{\{|X_{(t,x)}|>n\}}\,\lambda_{(t,x)}\,\d t\nu(\d x)\biggr]\nonumber\\
	&\qquad+2\mathbb{E}\biggl[\int_{\R_+\times E}|X_{(t,x)}|^2\ind_{\{|X_{(t,x)}|>n\}}\,\lambda_{(t,x)}\,\d t\nu(\d x)\biggr]\nonumber\\
	&\leq\mathbb{E}\biggl[\int_{n}^\infty\int_E|X_{(t,x)}|^2\ind_{\{|X_{(t,x)}|\leq n\}}\,\lambda_{(t,x)}\,\d t\nu(\d x)\biggr]\nonumber\\
	&\qquad+2\mathbb{E}\biggl[\int_{\R_+\times E}n^2\ind_{\{|X_{(t,x)}|>n\}}\,\lambda_{(t,x)}\,\d t\nu(\d x)\biggr]\nonumber\\
	&\qquad+2\mathbb{E}\biggl[\int_{\R_+\times E}|X_{(t,x)}|^2\ind_{\{|X_{(t,x)}|>n\}}\,\lambda_{(t,x)}\,\d t\nu(\d x)\biggr]\nonumber\\
	&\leq\mathbb{E}\biggl[\int_{n}^\infty\int_E|X_{(t,x)}|^2\ind_{\{|X_{(t,x)}|\leq n\}}\,\lambda_{(t,x)}\,\d t\nu(\d x)\biggr]\nonumber\\
	&\qquad+4\mathbb{E}\biggl[\int_{\R_+\times E}|X_{(t,x)}|^2\ind_{\{|X_{(t,x)}|>n\}}\,\lambda_{(t,x)}\,\d t\nu(\d x)\biggr],\nonumber
\end{align}
and this latter term tends to zero as $n\to\infty$ since $X\in\mathcal{P}_2(\lambda)$.
Hence $\mathcal{P}_{1,2}(\lambda)$ is dense in $\mathcal{P}_2(\lambda)$.
\\
\noindent$\square$

\end{document}